%% file: main.tex
\documentclass[preprint,11pt]{elsarticle}

\usepackage{geometry}
\usepackage{xcolor}
\usepackage{siunitx}
\usepackage[export]{adjustbox}
\usepackage{tabularx}
\usepackage{graphicx}
\usepackage{paralist}
\usepackage{url}
\usepackage{wrapfig}
\usepackage{amsmath}
\usepackage{amssymb}
\usepackage{placeins}
\usepackage{tikz, standalone}
\usepackage{caption}
\usepackage{subcaption}
\usepackage{ifthen,tikz}
\usepackage{pgfplots}
\pgfplotsset{compat=newest}
\usepackage{comment}
\usepackage{enumerate}
\usepackage{algorithm}
\usepackage{algpseudocode}
\usepackage{eucal}
\usepackage{enumitem}
\usepackage{booktabs}
\usepackage{nicefrac}

\usetikzlibrary{shapes,positioning,intersections,quotes, calc, angles}
\usepackage{color}

\newcommand{\MarginPar}[1]{\marginpar{
\vskip-\baselineskip
\raggedright\tiny\sffamily
\hrule\smallskip{\color{red}#1}\par\smallskip\hrule}}

\newcommand{\commentout}[1]{}

\usepackage{soul} % Use \st{} for strikeout, \ul{} for underline

\newcommand{\sfrac}[2]{\mathchoice
  {\kern0em\raise.5ex\hbox{\the\scriptfont0 #1}\kern-.15em/
   \kern-.15em\lower.25ex\hbox{\the\scriptfont0 #2}}
  {\kern0em\raise.5ex\hbox{\the\scriptfont0 #1}\kern-.15em/
   \kern-.15em\lower.25ex\hbox{\the\scriptfont0 #2}}
  {\kern0em\raise.5ex\hbox{\the\scriptscriptfont0 #1}\kern-.2em/
   \kern-.15em\lower.25ex\hbox{\the\scriptscriptfont0 #2}}
  {#1\!/#2}}

\newcommand{\half}{{\sfrac{1}{2}}}
\newcommand{\nph}{{n+\sfrac{1}{2}}}

\newcommand{\Diag}{\mathrm{Diag}}
\newcommand{\vfrac}{\Lambda}
\newcommand{\cons}{\mathrm{c}}
\newcommand{\ncons}{\mathrm{nc}}
\newcommand{\vol}{V}
\newcommand{\delM}{\delta M}
\newcommand{\nto}{u\rightarrow c}
\newcommand{\otn}{c\rightarrow u}
\newcommand{\gtv}{g\rightarrow v}
\newcommand{\vtg}{v\rightarrow g}
\newcommand{\NFRD}{\mathbf{N}}

\DeclareMathSymbol{\shortminus}{\mathbin}{AMSa}{"39}

\def\half   {\frac{1}{2}}

\def\U       {{U}}

\def\nph    {{n+\sfrac{1}{2}}}

\def\iph    {i+\sfrac{1}{2}}

\def\imh    {i-\sfrac{1}{2}}
\def\iph    {i+\sfrac{1}{2}}

\def\iph    {{i+\sfrac{1}{2},j,k}}
\def\ijph   {{i,j+\sfrac{1}{2},k}}
\def\ijkph  {{i,j,k+\sfrac{1}{2}}}
\def\imh    {{i-\sfrac{1}{2},j,k}}
\def\ijmh   {{i,j-\sfrac{1}{2},k}}
\def\ijkmh  {{i,j,k-\sfrac{1}{2}}}
\def\ivec   {{i,j,k}}

\def\IVEC   {{I,J,K}}

\journal{Journal of Computational Physics}

\begin{document}

\begin{frontmatter}

\title{A New Re-redistribution Scheme for Weighted State Redistribution with Adaptive Mesh Refinement}

\author[inst1]{I. Barrio Sanchez}

\affiliation[inst1]{organization={University of Pittsburgh},
            addressline={4200 Fifth Ave}, 
            city={Pittsburgh},
            postcode={15260}, 
            state={PA},
            country={USA}}

\author[inst2]{A. S. Almgren}

\affiliation[inst2]{organization={Lawrence Berkeley National Laboratory},
            addressline={1 Cyclotron Rd}, 
            city={Berkeley},
            postcode={94720}, 
            state={CA},
            country={USA}}

\author[inst2]{J. B. Bell}
\author[inst3]{M. T. Henry de Frahan}
\affiliation[inst3]{organization={National Renewable Energy Laboratory},
            addressline={15013 Denver W Pkwy}, 
            city={Golden},
            postcode={80401}, 
            state={CO},
            country={USA}}
\author[inst2]{W. Zhang}

%
%\date{}
%\maketitle

\begin{abstract}
\input{00_abstract/abstract}
\end{abstract}




\end{frontmatter}

% Section 1 - Introduction
\section{Introduction}
\input{01_introduction/background}

% Section 2 - Preliminaries
\section{Mathematical preliminaries}\label{sec:math-prelim}
\input{02_algorithms/prelim_intro}

\subsection{Notation}\label{sec:math}
\input{02_algorithms/notation}

\subsection{Single level algorithm}
\input{02_algorithms/finite_volume_intro}

%\subsection{Small Cell Problem}
%\input{02_algorithms/redistribution}

%\subsection{Flux Redistribution Scheme} \label{sec:FRD}
%\input{02_algorithms/FRD_intro}

\subsection{Multilevel algorithm}\label{sec:AMR}
\input{02_algorithms/multilevel_algorithm}

\subsection{Synchronization}\label{sec:refluxing}
\input{02_algorithms/synchronization}

%\subsection{Flux Re-redistribution} \label{sec:FRRD}
%\input{02_algorithms/flux_re_redistribution}

\section{Weighted state redistribution scheme} \label{sec:SRD}
\input{02_algorithms/SRD_intro}

\input{02_algorithms/srd_matrix}

% Section 4 - Incflo
%\section{Time-Stepping for Conservation Law}

%\subsection{Embedded Boundary Representation}
%\subsubsection{Simple Geometry Toolbox}
%\input{04_incflo/eb_representation/simple_geometry}

%\subsubsection{Complex Geometries}
%\input{04_incflo/eb_representation/complex_geometry}

% Section 5 - Results
\section{Numerical results}\label{sec:num}
In this section, we present numerical results demonstrating the properties of the proposed state re-redistribution method. We use two shock tube problems, with different initial conditions and grid refinement strategies to demonstrate the mass conservation properties (Section~\ref{sec:num_res_conservation}) and to provide a validation case (Section~\ref{sec:sod-validation}). A shock interacting with a cylinder provides additional validation in the presence of multiple AMR levels (Section~\ref{sec:cyl-validation}). These cases are run with the CAMR code~\cite{url:CAMR}. Finally, we demonstrate the use of state re-redistribution with a more realistic three-dimensional geometry and present some performance metrics showing the cost savings from only partially refining the embedded boundary (Section~\ref{sec:pelec-ic}).  This final problem is run with the compressible reacting flow solver, PeleC~\cite{PeleC_IJHPCA,PeleViz,Sitaraman2021}, part of the Pele suite of solvers. 
\input{04_results/results_2d}
\input{04_results/pele-sod}
\input{04_results/pelec-challenge}

% Section 6 - Conclusions
\section{Conclusions}\label{sec:conclusion}
\input{conclusions}

% Section 6 - Funding
\section{Funding}
The work of JBB was supported by the U.S. Department of Energy, Office of Science, Office of Advanced Scientific Computing Research, Applied Mathematics Program under contract No. DE-AC02-05CH11231.  The work of ASA was supported in part by the U.S. Department of Energy (DOE) Office of Advanced Scientific Computing Research (ASCR) via the Scientific Discovery through Advanced Computing (SciDAC) program FASTMath Institute.
The work of IBS and WZ was supported by the Exascale Computing Project (17-SC-20-SC), a collaborative effort of the U.S. Department of Energy Office of Science and the National Nuclear Security Administration.
This work was authored in part by the National Renewable Energy Laboratory, operated by Alliance for Sustainable Energy, LLC, for the U.S. Department of Energy (DOE) under Contract No. DE-AC36-08GO28308.  Funding was provided by U.S. Department of Energy Office of Science and National Nuclear Security Administration. This research was supported by the Exascale Computing Project (17-SC-20-SC), a collaborative effort of the U.S. Department of Energy Office of Science and the National Nuclear Security Administration. A portion of the research was performed using computational resources sponsored by the Department of Energy’s Office of Energy Efficiency and Renewable Energy and located at the National Renewable Energy Laboratory. The views expressed in the article do not necessarily represent the views of the DOE or the U.S. Government. The U.S. Government retains and the publisher, by accepting the article for publication, acknowledges that the U.S. Government retains a nonexclusive, paid-up, irrevocable, worldwide license to publish or reproduce the published form of this work, or allow others to do so, for U.S. Government purposes.

\appendix

\section{Flux redistribution scheme} \label{sec:FRD}
\input{02_algorithms/FRD_intro}

\section{Flux re-redistribution} \label{sec:FRRD}
\input{02_algorithms/flux_re_redistribution}

\commentout{
\section{Matrix form of flux redistribution}\label{sec:MFRD}
\input{02_algorithms/frd_matrix}

\section{Matrix form of flux re-redistribution}\label{sec:MFRRD}
\input{02_algorithms/frd_matrix_rere}

}

\section{Piecewise-linear Godunov with embedded boundaries}\label{sec:EBGdnv}
\input{02_algorithms/ebgdnv}

\newpage
% References
%\bibliographystyle{plain}
\bibliographystyle{elsarticle-num} 
\bibliography{biblio/asa.bib,%
              biblio/exascale.bib,%
              biblio/drangarajan.bib,%
              biblio/pele.bib,%
              biblio/refs.bib}

\end{document}

%% file: 00_abstract/abstract.tex
State redistribution (SRD) is a recently developed technique for stabilizing cut cells that result from finite-volume embedded boundary methods.  SRD has been successfully applied to a variety of compressible and incompressible flow problems.  When used in conjunction with adaptive mesh refinement (AMR), additional steps are needed to preserve the accuracy and conservation properties of the solution if the embedded boundary is not restricted to a single level of the mesh hierarchy. In this work, we extend the weighted state redistribution algorithm to cases where cut cells live at or near a coarse-fine interface within the domain.
The resulting algorithm maintains conservation and is demonstrated on several two- and three-dimensional example problems.

%% file: 01_introduction/background.tex
% General Intro
\commentout{
Fluid flows in and around irregular geometries are ubiquitous in the physical world and the numerical simulation of these flows, from the biological flows inside the human body to the high-speed flows around aircraft control surfaces, is at the heart of numerous research questions. For simulations of fluid flows that use an underlying mesh, there are two broad approaches to address irregular geometries. The first is a body-fitted approach in which the mesh itself is fit to the irregular geometry; this approach is used with both structured and unstructured meshes \cite{iaccarino2003}.
Body-fitted approaches can yield highly accurate solutions 

but can become difficult to use with moving geometries due to the required remeshing.
In this paper, we focus on using a non-body fitted approach.

Non-body fitted approaches typically use an underlying Cartesian mesh and represent the irregular geometries using an implicit function or an overlapping mesh. The way in which the Cartesian mesh is then coupled with the irregular geometry leads to the different families of methods. In the Immersed Boundary (IB) Method, pioneered by Peskin \cite{peskin1972flow}, the Cartesian mesh is coupled with the irregular geometry by introducing forcing terms at the physical location of the interface between the fluid and the irregular geometry. This forcing term can be viewed as force density that imposes the boundary condition between the fluid and the irregular geometry. In order to impose the correct boundary conditions and avoid modifications to the computational stencil, some means of interpolation from the location of the irregular geometry to the nearest point(s) on the Cartesian grid is used. 
%The introduction of these interpolation schemes can then lead to a 
 In \cite{gronskis2016simple}, the authors perform a thorough comparison of different interpolation schemes and present a novel scheme for irregular boundaries that is compatible with domain decomposition. The authors show that their new method, while not strictly conservative, is computationally efficient and accurate by solving a variety of canonical problems in two-dimensions.  

While the IB method has been successfully used in a wide variety of settings (see, e.g., \cite{lin2007fixed, taira2007immersed}) since its introduction, it does suffer from drawbacks. IB schemes often suffer from loss of conservation and volume leakage.  In addition, the forcing conditions used to couple the fluid flow and irregular geometry are problem-specific, requiring the need to tune the method based on physical and material characteristics. Additionally, in practice the forcing terms are often applied in a band around the interface between the fluid and geometry to avoid stiffening the equations and inducing non-physical oscillations \cite{iaccarino2003}. This band approach can be made conservative, but leads to smearing of the solution near the interface and reduced accuracy when the boundary is closely coupled with the surrounding fluid. Other approaches, such as the direct-forcing scheme in \cite{mohd1997combined} and \cite{fadlun2000combined}, avoid smearing the forcing term, but are not conservative. For a more detailed review of IB methods, including recent advances in addressing conservation and leakage, see \cite{huang2019recent} and the references therein.
%\ann{This seems redundant -- you introduced a review-type article in previous paragraph, now giving a different one here.  This paragraph and previous need to be better integrated} \matt{Let's talk about this. The comparison noted above looked at different interpolation schemes, whereas the review paper here is truly a review of the IB methods.}

%{\color{red} what about loss of conservation / "leakage" in IB methods?}

%In order to preserve the numerical solution of flow problems near the interface, a
A number of sharp interface methods have also been developed. The immersed interface method (IIM) of LeVeque and Li \cite{iim1_lai2001remark, iim2_li2001immersed, iim3_leveque1994immersed} was the earliest implementation of a sharp interface method applied to two-dimensional incompressible flow with a moving geometry. The IIM uses a standard finite difference scheme at all of the regular grid points away from the interface. At the interface, the IIM uses second order Taylor expansions on both sides of the interface and a jump condition across the interface. The IIM achieves second order accuracy, but its implementation is complex and difficult to extend to three dimensions.  

The embedded boundary / cut-cell approach
%is a simpler method that couples the Cartesian mesh with the irregular geometry by creating 
intersects the solid geometry with the mesh and defines separate, non-overlapping fluid and solid regions.
%at the physical location of the interface. These intersections effectively cut the cells along the interface between the fluid and irregular geometry without any restrictions on how the intersections are placed. 
The process of mesh generation for the cut-cell approach is robust since it simply involves computing intersections, and does not impose constraints on the resulting mesh,
%and easily automated, 
but can also result in cut cells near the interface that have arbitrarily small volumes. For explicit time-stepping methods in the context of a standard finite volume discretization, this results in a prohibitively restrictive requirement for the time step \cite{almgren1997cartesian}. This is referred to as the "small cell problem" and is one of the major challenges when using the embedded boundary / cut-cell approach. 
%\deepak{perhaps worth adding some context for what "easily automated" means.} \matt{this is a good point} \ann{let's just remove that comment? see my rewording above} \matt{Looks good.}

In general, there are two approaches to alleviating the problems that arise from the small cell problem: cell merging and dynamic stabilizing procedures. Cell merging is a geometric manipulation that removes small cells by merging them with neighboring cells, either regular or cut. The dynamic stabilizing procedures address the problem by stabilizing cut cells in a post-processing step that replaces the values at the cut cell with values that maintain accuracy and conservation to varying degrees.

There is a large body of research that covers the use of embedded boundary / cut-cell approaches applied to static geometries for incompressible flows \cite{almgren1997cartesian, ye1999accurate}, static geometries for compressible flows \cite{colella2006cartesian}, and moving geometries for compressible flows \cite{min2007second, MEINKE2013135, SCHNEIDERS2013786, natarajan2022moving}. 

The cut-cell approach of Udaykumar \cite{udaykumar2001sharp, udaykumar2002interface, udaykumar2003sharp} uses a standard finite volume approach away from the interface and reshapes the control volumes at the cut cells. This cell merging technique was based on \cite{ye1999accurate}, where small cells were merged into their nearest neighbor creating trapezoidal cells. %Then, the weak form of the governing equations are integrated over the non-rectangular control volumes. 
For moving geometries, when newly uncovered cells appear the values are obtained by interpolation from the known values in the surrounding cells and on the moving boundary. This approach is second-order accurate and has been used to solve a variety of problems in fluid-structure interaction \cite{udaykumar2002interface, udaykumar2003sharp} and phase change \cite{udaykumar2001sharp}.  

In \cite{marella2005sharp}, the authors moved away from the cell merging approach of \cite{udaykumar2001sharp} by using a finite difference discretization near the interface and combining the difference operators with the level set function in a manner similar to \cite{gibou2002second, gibou2003level}. Their finite difference approach yielded second-order accuracy and was easily extended to three dimensions. This finite difference approach is not conservative, but the idea did seem to get picked up in a number of ensuing works. In \cite{ng2009efficient}, the authors incorporate a similar idea and directly discretize the projection operator used to solve the incompressible Navier-Stokes equations. They use a high order extrapolation method to extend their computational stencil across the interface when computing the intermediate velocity based on the procedure in \cite{aslam2004partial} and the modification in \cite{min2007second}. Later, this projection method was extended to work with a quadtree/octree adaptive mesh refinement (AMR) strategy in \cite{guittet2015stable} then again to a highly parallel version in \cite{egan2021direct}. 

An alternative to cell merging, that has been extensively used for compressible flows, is that of dynamic stabilization. In flux redistribution schemes, a conservative but potentially unstable update for the solution is first computed using the difference of area-weighted fluxes.
%the standard flux balances for cut cells that are locally computed without regard to the small cell problem. 
The solution in cut cells is then stabilized/regularized by adjusting the update to the solution to be stable, and 
%a only a fraction of the flux and 
redistributing the remainder of the update to neighboring cells in a conservative manner. Flux redistribution schemes are easily extended to three dimensions due to their simplicity, but are often only first-order accurate at the boundaries. Examples of flux redistribution schemes applied to compressible equations can be found in \cite{pember1995adaptive, colella2006cartesian, hu2006conservative, klein2009well, SCHNEIDERS2013786}. 

%\deepak{more of a question than a comment. is the first order accuracy only on the cut cells? i think in general FRD algorithm in hydro is still 2nd order, no?} \ann{FRD definitely drops order at the boundaries -- it is more diffusive than SRD.   The statement is that the solution is 2nd order away from the boundary, 1st order at the boundary, so somewhere in between in a global L1 or L2 norm.  See my rewording in paragraph above} \deepak{rewording looks good}

More recently, the state redistribution algorithm of Berger and Giuliani \cite{berger2021state} has shown great promise as a replacement for flux redistribution (FRD). State redistribution, like flux redistribution, is a post-processing step. One difference between state and flux redistribution schemes is that the post-processing step in SRD is based on the conserved state variables themselves, rather than the updates to the variables. The state redistribution algorithm, like flux redistribution, is fully conservative and has successfully been used to stabilize both finite volume \cite{berger2021state} and discontinuous Galerkin \cite{giuliani2022dgsrd} solvers.
%It can be generalized to high order accuracy. \ann{not sure we want to say higher order here since we're not going to go there}  
In \cite{giuliani2022weighted}, the authors extend the original state redistribution algorithm to three dimensions and include a number of algorithmic improvements. In particular, the authors introduce a weighted version of the algorithm that is easily generalized and demonstrate that state redistribution can stabilize both advective and diffusive contributions to the solution update. We adopt the weighted state redistribution (WSRD) approach for the work presented herein.

It is important to note that cell merging techniques can be viewed as a post-processing step similar to the state redistribution algorithm. The primary difference between these two methods, however, is that the state redistribution supports overlapping neighborhoods whereas cell merging does not. This allows the state redistribution algorithm to be easily extended to three dimensions. Cell merging method, on the other hand, can be difficult to implement in three dimensions as there are many different, possibly incompatible, ways to create non-overlapping cells. For a recent example of cell merging techniques for moving geometries applied to two- and three-dimensions, see \cite{saye2017implicit} and \cite{saye2017part2}. 
}
Numerical solution of partial differential equations in complex geometric domains has a wide range of important scientific and engineering applications.  
Embedded boundary methods, also known as ``cut cell" methods, provide an attractive approach for these types of problems since mesh generation is robust and automatic. The basic idea of the embedded boundary / cut-cell approach is to treat the domain boundary, or solid object within the domain, as a surface intersecting a regular Cartesian grid.
The process of mesh generation for the embedded boundary approach then reduces to intersections of the surface with regular finite volume cells.

The difficulty that arises with this process is that cells cut by the interface may have arbitrarily small volumes.  This is referred to as the ``small cell problem" and is one of the major challenges when using the embedded boundary / cut-cell approaches. 
A number of different approaches have been proposed to address the small cell problem.   A novel dimensionally split approach was proposed by Gokhale {\it et al.} \cite{gokhaleNikosKlein:2018}.  $H$-box methods \cite{mjb-hel-rjl:hbox} have nice theoretical properties but have not been extended to three dimensions.
Cell merging, in which small cells are merged with larger cells, is an intuitive approach that has been successfully used to address the small cell problem, and has seen significant recent progress. See, for example, \cite{MURALIDHARAN2016,saye2017implicit,saye2017part2,GULIZZI2022}.

 An alternative approach to cell merging that has been successfully used for simulation of 
 compressible flow \cite{pember1995adaptive,colella2006cartesian,hu2006conservative, klein2009well,graves2013cartesian,SCHNEIDERS2013786} and incompressible flow 
 \cite{almgren1997cartesian,trebotich:2015}
is called flux redistribution (FRD).  In flux redistribution schemes, a conservative but potentially unstable update for the solution is first computed using the difference of area-weighted fluxes divided by the cell volume.
%the standard flux balances for cut cells that are locally computed without regard to the small cell problem. 
The solution in cut cells is then stabilized/regularized by adjusting the update to the solution to be stable, and 
%a only a fraction of the flux and 
redistributing the remainder of the update to neighboring cells in a conservative manner.
This approach is attractive because it is simple to implement as a postprocessing step. However, the downside is that there is a loss of accuracy at the cut cells and the scheme is not linearity preserving.

% Recently, a new approach called state redistribution (SRD) was proposed for finite volume methods on two-dimensional cut cell grids \cite{BG}.  This is a minimally invasive stabilization technique that is linearity preserving, conservative, and straightforward to implement in two dimensions for hyperbolic conservation laws.  Inspired by flux redistribution, SRD postprocesses the numerical solution by accurately redistributing the solution states in a way that maintains conservation.

Recently, the state redistribution algorithm (SRD) first proposed by Berger and Giuliani \cite{berger2021state} has shown great promise as an alternative to flux redistribution. The state redistribution scheme, like flux redistribution, starts by computing a conservative but potentially unstable update to the solution.
%\MarginPar{We should be more clear about what we mean by post-processing, because when we define the WSRD algorithm we talk about a pre-processing and post-processing step}
However, SRD updates the conserved variables and then redistributes the updated solution, whereas FRD  redistributes the update before updating the conserved variables. State redistribution, like flux redistribution, is fully conservative and has successfully been used to stabilize both finite volume \cite{berger2021state} and discontinuous Galerkin \cite{giuliani2022dgsrd} solvers.
Unlike the FRD approach, SRD is linearity preserving, which leads to improved accuracy at the boundary.
%It can be generalized to high order accuracy. \ann{not sure we want to say higher order here since we're not going to go there}  
Giuliani {\it et al.} \cite{giuliani2022weighted} extend the original state redistribution algorithm to three dimensions and introduce a number of algorithmic improvements. In particular, the authors introduce a weighted version of the algorithm (WSRD) that is easily generalized and demonstrate that state redistribution can stabilize both advective and diffusive contributions to the solution update. The authors also demonstrate that the methodology can be extended to incompressible flows.
Berger and Giuliani \cite{berger2023new} introduce an algorithm for selecting weights that makes WSRD monotone, total variation diminishing and GKS stable in most situations.
%\MarginPar{add ref to new Marsha paper}

The basic concept of embedded boundary approaches, namely, that they are based on an underlying regular finite volume grid, makes them a natural candidate for block-structured adaptive mesh refinement (AMR) of the type initially proposed by Berger and collaborators \cite{berger1989local,bell1994three}.  Several of the algorithms referenced above use this type of AMR.  In AMR algorithms for hyperbolic conservation laws, the data at each level are advanced independently of the other levels, other than using Dirichlet data provided from coarser grids as Dirichlet ``boundary conditions" where appropriate.   
 When the data on grids at different levels reach the same time they must be synchronized to construct a stable, conservative composite solution. When no cut cells are present, this synchronization has two parts. The first is re-defining the coarse solution as the average of the fine solution where possible; this changes solution values only on coarse cells covered by fine cells.  The second part is ``refluxing'', whereby the fluxes used to update the coarse solution are effectively over-written by the time- and space-averaged fine fluxes, thus modifying the coarse data adjacent to, but not covered by, fine grids. Redistribution algorithms for embedded boundary methods require additional steps in the synchronization process when coarse/fine interfaces are close to the embedded boundary.  Pember {\it et al.} \cite{pember1995adaptive} developed a so-called {\it re-redistribution} algorithm that addresses this additional synchronization for flux redistribution.
% \MarginPar{read above paragraph for changes}
%In AMR algorithms for hyperbolic conservation laws, the overall time-stepping approach can be viewed as a recursive operation, in which for each full time step, the data on coarse grids are advanced in time, then the coarse data is interpolated in space and time to supply boundary conditions to advance the data on fine grids, possibly with subcycling in time, and so on recursively until the finest grid level is reached.  When the data on grids at different levels reach the same time they must be synchronized to construct a stable, conservative composite solution from data at different levels that has been updated separately. Embedded boundary methods introduce additional complexity into the synchronization process when coarse /fine interfaces are close to the embedded boundary.  Pember {\it et al.} \cite{pember1995adaptive} developed a so-called {\it re-redistribution} algorithm that address this additional synchronization for flux redistribution.  

The goal of this paper is to extend the weighted state redistribution algorithm to enable AMR without the restriction that the entire embedded boundary be resolved at the same level.  The focus here is on the AMR synchronization process when a coarse/fine boundary is near an embedded boundary.  (If the entire embedded boundary is always at the same, typically finest, level, then no special synchronization is required and the current algorithm as presented in \cite{giuliani2022weighted} is sufficient.)
We will present the algorithm for a system of hyperbolic conservation laws.  The implementation of the re-redistribution algorithm is done in the finite-volume, block-structured adaptive mesh refinement framework of AMReX \cite{AMREX:IJHPCA}.  

The paper is organized as follows. In Section~\ref{sec:math-prelim}, we review the embedded boundary discretization of conservation laws and the basics of AMR.
%, with FRD used as an exemplar of the basic ideas.  
In Section~\ref{sec:SRD}, we review the  weighted state redistribution (WSRD) algorithm and describe how to do re-redistribution for WSRD. In Section~\ref{sec:num}, we validate our numerical scheme using canonical two- and three-dimensional example problems and present an application of the methodology to a more realistic flow problem.  Finally, we summarize our conclusions in Section~\ref{sec:conclusion}.  %\ref{sec:FRD} contains a detailed description of the flux redistribution scheme.
%\MarginPar{I think we should say here what is in the appendices}
The appendices contain 
the details of the flux redistribution algorithm (Appendix A), the re-redistribution algorithm for flux redistribution (Appendix B) and the necessary modifications to the Godunov advance in the presence of cut cells (Appendix C).

%% file: 02_algorithms/prelim_intro.tex
%In this section, we present a simple conservation law equation in order to explain the extension of the WSRD scheme to incorporate re-redistribution.  In a later section we will describe how this is extended to more interesting systems such as compressible and incompressible, variable density flows.
In this section we review the basic concepts of embedded boundary discretizations, discuss the redistribution approach to addressing the small cell problem, and introduce adaptive mesh refinement (AMR).  For the discussion here, we consider a system of hyperbolic conservation laws in three dimensions, 
\begin{eqnarray} \label{eq:conservation_law}
    \frac{\partial U}{\partial t} + \nabla \cdot F(U) = 0
\end{eqnarray}
where $U$ is the state variable and $F(U) = (F^x(U),F^y(U), F^z(U) \;)$ is the flux.

%% file: 02_algorithms/notation.tex
Embedded boundary (EB) methods are used to discretize partial differential equations (PDEs) in complex domains.  The basic idea is to represent an irregular boundary (between ``fluid" and ``body", where we wish to solve the PDEs in the fluid region only) by intersecting an explicitly specified boundary with a uniform finite-volume Cartesian grid with mesh spacings $\Delta x$, $\Delta y$ and $\Delta z$.  This representation 
introduces irregularly shaped cells only adjacent to the boundary. 
Following standard notation, we define each grid cell $(i,j,k)$ to be either 
\textit{body}, \textit{cut}, or \textit{regular}. We define the geometric
volume fraction, $\Lambda$, of each cell to be the fraction of that rectangular cell volume that is inside the fluid region: body cells have $\Lambda = 0,$
regular cells have $\Lambda = 1$, and for cut cells $0 < \Lambda < 1.$ 
Area fractions $a$ are stored on each cell face, again with values in $[0,1],$
representing the fraction of the face not covered by the body. 
%Finally, the location of the cell centroid (which for regular cells is identical to the cell center), and the locations of the face centroids are stored; these are scaled by the cell widths and so have values in $[-0.5,0.5]$ for each coordinate direction.
%\MarginPar{not sure we are using cell centroids. i think we are treating values as living at cell centers here}
%Additionally, there is connectivity information between neighboring cells.
In the AMReX implementation used here, the EB information needed to represent the geometry, such as face area fractions and cell volume fractions, is precomputed and stored in a distributed database at the beginning of the calculation.
For additional details on the embedded boundary implementation, we refer to the 
AMReX documentation \url{https://amrex-codes.github.io/amrex/docs_html/}. 

%In the presence of embedded boundaries, we define area fractions and volume fractions as the fraction of the faces and volumes, respectively, that are in the fluid region rather than the body region.  These fractions have values in $[0,1].$  We call a cell with volume fraction 1 a {\it regular} cell; a cell with volume fraction 0 is denoted a {\it covered} cell, and a cell with volume fraction between 0 and 1 is a {\it cut} cell.  I

Going forward, we will use the notation $A$ to represent the face area of a regular cell multiplied by the area fraction of that particular face, for example $A^x_\iph = a_\iph \Delta y \Delta z$, and $\vol$ to represent the volume of a regular cell multiplied by the volume fraction of that particular cell, i.e. $\vol_\ivec = \Lambda_\ivec \Delta x \Delta y \Delta z$.

%% file: 02_algorithms/finite_volume_intro.tex
For the discussion here, we will consider a single-step integration scheme typical of higher-order Godunov-type discretizations.  Generalizations to a method of lines approach are discussed at the end of Sec. \ref{sec:SRD}.  We let $\U_\ivec$ and $F^x_\iph$, $F^y_\ijph,$ and $F^z_\ijkph$ represent the discretized solution $U$, and the normal flux $F$ on $x$-faces,  $y$-faces and  $z$-faces, respectively.
We define a standard finite volume scheme in the form,
\commentout{
\begin{eqnarray} \label{eq:basic_fv}
    \U_\ivec^{n+1} = \U_\ivec^n &- \Delta t& 
    \frac{\left (A^x_{\iph} F^x_{\iph} - A^x_{\imh} F^x_{\imh} \right )}{\vol_\ivec} \\ &- \Delta t&
    \frac{\left (A^y_{\ijph} F^y_{\ijph} - A^y_{\ijmh} F^y_{\ijmh} \right )}{\vol_\ivec} \nonumber\\ &- \Delta t&  
    \frac{\left (A^z_{\ijkph} F^z_{\ijkph} - A^z_{\ijkmh} F^z_{\ijkmh} \right )}{\vol_\ivec} \nonumber\\
     &- \Delta t& \frac{A^f_\ivec F^f_\ivec}{\vol_\ivec} \nonumber.
\end{eqnarray}
}
\begin{eqnarray*}
U_\ivec^{n+1} &=& U_\ivec^n + \Delta t \; \delta U^c_\ivec 
\end{eqnarray*}
where
\begin{eqnarray} 
\label{eq:cons_up}
\delta U^c_\ivec = 
&-& 
    \frac{\left (A^x_{\iph} F^x_{\iph} - A^x_{\imh} F^x_{\imh} \right )}{\vol_\ivec} \\ 
&-& 
    \frac{\left (A^y_{\ijph} F^y_{\ijph} - A^y_{\ijmh} F^y_{\ijmh} \right )}{\vol_\ivec} \nonumber \\
&-&   
    \frac{\left (A^z_{\ijkph} F^z_{\ijkph} - A^z_{\ijkmh} F^z_{\ijkmh} \right )}{\vol_\ivec} \nonumber \\
&-&  
    \frac{A^f_\ivec F^f_\ivec}{\vol_\ivec} \; \nonumber.
\end{eqnarray} 
The update is conservative by construction, and the time accuracy of the scheme is determined by the details of the flux construction.
Here, $A^f_\ivec$ represents the area of the boundary intersected with cell ($\ivec$) and $F^f_\ivec$ represents the flux in the normal direction at the boundary.
We note that $A^f_\ivec$ and the unit normal to the surface $\mathbf{n}^f_\ivec$  (pointing from body to fluid) can be approximated using 
\begin{align*}
A^f_\ivec \mathbf{n}^f_\ivec &= (A^x_{\iph} - A^x_{\imh}) \; \mathbf{i} + \; 
(A^y_{\ijph} -  A^y_{\ijmh})\; \mathbf{j} + \; 
(A^z_{\ijkph} - A^z_{\ijkmh})\; \mathbf{k}   \; .
\end{align*}
This formula is exact for planar surfaces and gives second-order averages for smooth surfaces.

In the absence of cut cells, where all cells have the same volume $V = \Delta x \Delta y \Delta z$, the CFL condition for an explicit update constrains $\Delta t$ based on the ratio of cell volume $V$ to face areas $A$ as well as the maximum wave speed in the problem, i.e. the time step depends linearly on the mesh spacing.   
The ``small cell problem" discussed in the introduction refers to the fact that 
when an embedded boundary representation of geometry is used, the cut cells can have arbitrarily small volumes, making the update in Eq.~(\ref{eq:cons_up}) unstable if the time step calculation does not take into account the reduction in cell volume \cite{almgren1997cartesian}.  Reducing the time step to account for arbitrarily small cell volumes would be too computationally expensive, thus as
%For explicit time-stepping methods in the context of a standard finite volume discretization, this results in a prohibitively restrictive requirement for the time step 
discussed above there are a number of approaches to addressing the small cell problem, all of which allow the computation of $\Delta t$ based only on the full-cell CFL constraint.  Recalling that the FRD scheme modifies the update itself, we can write the single-level single-timestep advance using FRD as 
\begin{eqnarray} \label{eq:basic_fv_frd}
U_\ivec^{n+1} &=& U_\ivec^n + \Delta t \; R^{FRD} ( \delta U^c_\ivec )
\end{eqnarray}
where $R^{FRD}$ represents the flux redistribution operator.  Details of the flux redistribution algorithm are discussed in \ref{sec:FRD}. Analogously, since WSRD operates on the full solution, we can write 
\begin{eqnarray} \label{eq:basic_fv_srd}
U_\ivec^{n+1} &=& R^{WSRD} (U_\ivec^n + \Delta t \;  \delta U^c_\ivec )
\end{eqnarray}
where $R^{WSRD}$ represents the state redistribution operator, which will be described in detail in Section~\ref{sec:SRD}.

%One type of approach is based on flux redistribution (see Pember {\it et al.} \cite{}).  The idea of flux redistribution is to compute a stable but non-conservative  update for small cells and to redistribute the remainder of the update to neighboring cells so that the overall algorithm is conservative and stable for $\Delta t$ based on the full-cell CFL constraint.
%Operationally flux redistribution modifies the conservative update $\delta U^c_\ivec$ to obtain a new update $\delta U^{EB}_\ivec$.

%and is one of the major challenges when using the embedded boundary / cut-cell approach.  
%In this paper we will discuss two different redistribution schemes for addressing this problem.

%We define the single-level time-stepping algorithm in two steps (for convenience later in expressing the redistribution algorithm): 
\commentout{
As a prelude to discussing how to address the small cell problem, we compute fluxes on each face and define a conservative update 
\begin{align} 
\label{eq:cons_up}
\delta U_\ivec = 
&- 
    \frac{\left (A^x_{\iph} F^x_{\iph} - A^x_{\imh} F^x_{\imh} \right )}{\vol_\ivec} \\ &- 
    \frac{\left (A^y_{\ijph} F^y_{\ijph} - A^y_{\ijmh} F^y_{\ijmh} \right )}{\vol_\ivec} \nonumber
    \\ &-  
    \frac{\left (A^z_{\ijkph} F^z_{\ijkph} - A^z_{\ijkmh} F^z_{\ijkmh} \right )}{\vol_\ivec} \nonumber \\
     &- \frac{A^f_\ivec F^f_\ivec}{\vol_\ivec} \nonumber 
\;\;\;.
\end{align} 

If there were not small cells then the solution update could be expressed as
\begin{eqnarray*} 
U_\ivec^{n+1} = U_\ivec^n + \Delta t \; \delta U_\ivec \;\;\;
\end{eqnarray*} 
}

%% file: 02_algorithms/multilevel_algorithm.tex
We consider a block-structured adaptive mesh refinement algorithm (AMR) in which regions at a given level requiring additional resolution are covered with grid patches of finer resolution.  (For our purposes, a ``level" refers to all the grids with the same resolution, and higher levels denote finer resolution.  We call two levels ``adjacent" if they are not separated by any intermediate levels.) With this type of refinement approach, one can think of the time stepping algorithm as a recursive procedure, starting with the coarsest level.  Data on coarse grids are advanced first, then data from the  coarse grids are interpolated to provide Dirichlet boundary conditions for the grids at the next finer level (if there is one). We refer to boundary cells for fine grids that are filled from coarse data as ``ghost" cells.  Data on the finer grids are then advanced, possibly with smaller time steps, until they reach the same time at the coarse grid. Advancing the fine grids with smaller time steps than the coarse grid is referred to as subcycling.  If there are still finer levels, this process repeats recursively.  Once the data on grids at adjacent levels reach the same time, the data on the two levels are synchronized.  
%At each level, once the coarser grids have supplied boundary conditions, grids at a given level are advanced independent of the other levels.  

There are a number of different facets to an AMR algorithm including refinement criteria, data distribution and load balancing, that we will not consider here. Here we will focus on how the data at adjacent levels are synchronized to create a composite solution.  
%As we consider issues related to embedded boundaries, it is important to differentiate a cell inside the solid body from a coarse cell covered by fine grid in an AMR algorithm;  for the latter we will use the term {\it covered}.
We define the refinement ratio $r_\ell$ to be the ratio of the 
mesh spacing at level $\ell$ to that at level $(\ell+1)$, and let $n_\ell$ be the number of level $\ell+1$ time steps corresponding to a single level $\ell$ time step.
To advance level $\ell,$ $0 \leq \ell \leq \ell_{\mathrm max}$ the following steps are taken.
%Below we summarize the AMR time stepping algorithm

\vspace{.1in}
\noindent Starting with level $\ell = 0,$ 
\begin{itemize}
\item{Fill boundary conditions for $U$ from level $\ell-1$ if level $\ell > 0$,
and from the physical domain boundary conditions. }
\item Define $U^{\ell,n+1}$ using Eq.~(\ref{eq:basic_fv_frd}) or Eq.~(\ref{eq:basic_fv_srd}) above.
%Compute the conservative update, $\delta U,$ using the algorithm discussed above and advance the solution, $U$, at level $\ell$.
%\[
%U^{\ell,n+1}_\ivec = U^{\ell,n}_\ivec + \Delta t^\ell \; \delta U^{\ell}_\ivec 
%\]
%Here $\delta U^{\ell}_\ivec$ is given by 
%uses the conservative update given by Eq. (\ref{eq:cons_up}) for cells sufficiently far from the embedded and uses $\delta U^{EB}_\ivec$ near the embedded boundary.

\commentout{
\begin{eqnarray*} 
\delta U^{\ell,\nph}_\ivec &=& -
\frac{\left (A^\ell_{\iph} F^\ell_{\iph} - A^\ell_{\imh} F^\ell_{\imh} \right )}{V_\ivec} \\ && -
\frac{\left (A^\ell_{\ijph} F^\ell_{\ijph} - A^\ell_{\ijmh} F^\ell_{\ijmh} \right )}{V_\ivec} \\ && -  
\frac{\left (A^\ell_{\ijkph} F^\ell_{\ijkph} - A^\ell_{\ijkmh} F^\ell_{\ijkmh} \right )}{V_\ivec} \\ 
U^{\ell,n+1}_\ivec &=& U^{\ell,n}_\ivec + \Delta t^\ell \; \delta U^{\ell,\nph}_\ivec 
\end{eqnarray*} 
}
\item Then, if $\ell < \ell_{max}$
\begin{itemize}
\item{Advance $U$ at level ($\ell+1$) for $n_\ell$ time steps
with time step $\Delta t^{\ell+1} = \frac{1}{n_\ell} \Delta t^{\ell}.$}
\item Synchronize the data between levels $\ell$ and $\ell+1$
%\begin{itemize}
%\item Volume average $U$ at level $\ell+1$ onto level $\ell$ grids.
%\item Correct $U$ in all level $\ell$ cells adjacent to but not covered by the union of
%level $\ell+1$ grids through an explicit refluxing operation as described below.
\end{itemize}
\end{itemize}

%% file: 02_algorithms/synchronization.tex
\begin{figure}
\centering
  \includegraphics[width=0.5\linewidth, clip=true, trim=0cm 1cm 0cm 0cm]{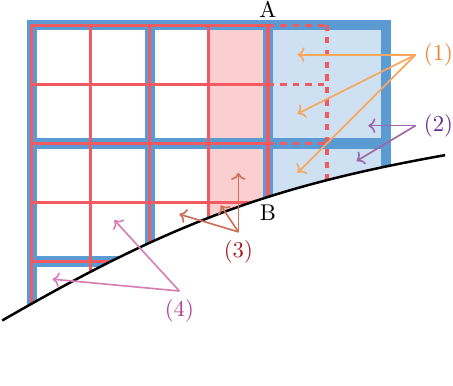}
  \caption{Schematic of the different type of cells and levels close to an EB (black line). Blue grid lines: coarse level, red grid lines: fine level. Orange (1) fine ghost cells; purple (2) uncovered coarse cells; burgundy (3) valid fine cells; and magenta (4) covered coarse cells. The coarse-fine interface lies between A and B. Note that we have not labeled all of the cells to the left of $\overline{AB}$; all of the coarse cells left of $\overline{AB}$ are covered and all of the fine cells left of $\overline{AB}$ are valid. Here, cells with a face that is part of $\overline{AB}$ are considered to be at the coarse-fine boundary. Shaded blue cells indicate the coarse cells at the coarse-fine interface; shaded red cells denote the fine cells at the coarse-fine interface.}
  \label{fig:amr-eb}
\end{figure}

As mentioned in the introduction, there are several steps in the synchronization procedure.  First, we want to replace the coarse solution by the average of the fine solution where possible.   Here we define a ``covered" cell at level $\ell$ as one which is overlaid by $r_{\ell}^d$ level $(\ell+1)$ grid cells, where $d$ is the spatial dimension of the problem. In addition we define a coarse cell to be ``at the coarse-fine boundary" if it shares a face with fine cells but is not covered, and a fine cell to be ``at the coarse-fine boundary" if it shares a face with a coarse cell that is not covered. These definitions are illustrated in Figure~\ref{fig:amr-eb}.  For EB methods we replace the solution in covered coarse grid cells with the volume-fraction-weighted average of the fine grid solution.  We refer to this step as the {\it average down} step.
%\MarginPar{define this and refluxing here so we can refer to them later}

%\MarginPar{add to caption once cells "at" the c-f boundary are labeled}

The other steps in the synchronization address issues that arise from advancing the levels independently.  Although the operations are defined procedurally, we note that these operations, in combination with the level-by-level update, implicitly define a composite discretization. (See \cite{bell2011adaptive} for a more general discussion of this perspective.)
%\MarginPar{missing ref}
One of these issues is that at the coarse-fine boundary, the flux computed on the coarse grid does not match the average of the fine grid fluxes.  (We note that the average of the fine grid area fractions does in fact equal the coarse grid area fraction, but the values of the flux itself need not match.)
Correcting for this mismatch is referred to as {\it refluxing}.
We compute the difference between the coarse grid flux and the time and space average of the fine grid fluxes on every face along the coarse-fine boundary. This difference
is used to correct the solution in the coarse grid cells at the coarse-fine boundary, i.e. coarse cells adjacent to the faces where the fluxes are corrected but not covered by finer grid (Figure~\ref{fig:amr-eb}).
 
For example, if coarse cell $(I,J,K)$ is located at the coarse-fine boundary and cell $(I-1,J,K)$ is covered, then the refluxing step updates $U_{I,J,K}^\ell$ by adding 
\begin{align}
%U^\ell_\IVEC := U^\ell_\IVEC &-
%\frac{  \Delta t^{\ell}  }{\vol_\IVEC}  A^{x,\ell}_{I-\half,J,K} F^{x,\ell}_{I-\half,J,K}
%+ \frac{1}{\vol_\IVEC} \sum_{n=1,n_\ell} \sum_{j,k} 
%\Delta t^{\ell+1}  
%A^{x,\ell+1}_{i+\half,j,k} F^{x,\ell+1}_{i+\half,j,k} \\
%&= U^\ell_\IVEC  + 
\frac{1}{\vol_\IVEC} \delta F^{x,\ell}_{I-\half,J,K} \nonumber
\end{align}
where
%\MarginPar{Added a minus right after the equals sign}
\begin{equation}
\label{eq:reflux}
\delta F^{x,\ell}_{I-\half,J,K} \equiv   - \Delta t^\ell A^{x,\ell}_{I-\half,J,K} F^{x,\ell}_{I-\half,J,K}
+ \sum_{n=1,n_\ell} \sum_{j,k} 
\Delta t^{\ell+1}  
A^{x,\ell+1}_{i+\half,j,k} F^{x,\ell+1}_{i+\half,j,k,n}
\end{equation}
represents the extensive flux correction associated with edge $({I-\half,J,K})$.
Here $i$ is the index for fine cells at the coarse/fine boundary and the summation in $j$ and $k$ ranges over the fine grid cells at the coarse/fine boundary that cover cell $(I-1,J,K)$ and $n$ refers to substeps of the fine grid corresponding to a coarse grid time step.
Note that Eq. (\ref{eq:reflux}) incorporates the orientation of the face; if the coarse cell at the boundary is instead on the left side of the fine grid the signs of the flux terms are reversed.

For non-EB problems, and when cut cells are sufficiently far from coarse/fine boundaries, 
refluxing is the only synchronization that is needed to make the composite solution conservative.
However, when there are cut cells near the coarse/fine boundary, additional synchronization is required
in order to account for the quantity of each conserved variable that has effectively crossed the coarse-fine boundary
during the redistribution step which occurred independently at the coarse and fine levels.  This additional {\it redistribution synchronization} is analogous to the refluxing procedure in that during the redistribution step at each level we accumulate the mismatch, then once the coarse and fine data have reached the same time, we modify the coarse data at the coarse/fine boundary to correct the mismatch.  

While this process is sufficient to make the composite solution conservative, there is an additional step required to maintain stability. 
Specifically, since the coarse cells whose data is modified by refluxing and the redistribution synchronization step can be arbitrarily small,  we need to redistribute the corrections from these procedures.
Details of these two steps, which collectively are referred to as re-redistribution, were first discussed in Pember {\it et al.} \cite{pember1995adaptive} in the context of flux redistribution; 
see \ref{sec:FRRD} for a summary of this approach.
%See \ref{sec:FRRD} for a detailed discussion of the additional synchronization steps for flux redistribution.

%% file: 02_algorithms/SRD_intro.tex
%State redistribution is a modern approach to stabilize cut-cells that arise in finite volume discretization using embedded boundary representations. If we examine the evolution equation, we see that the stability of this method requires the time step, $\Delta t$, to be proportional to the cell size, $\Delta x_i$. This creates a prohibitively restrictive requirement on the time step as the size of cut cells can be arbitrarily small. 

State redistribution \cite{berger2021state,giuliani2022weighted} has recently emerged as an alternative to flux redistribution for addressing the small cell problem associated with embedded boundary methods.
%State redistribution addresses the problem by creating a temporary merging neighborhood where small cut cells are merged with their neighbors until the neighborhood volume is greater than a specified threshold. 
State redistribution addresses the problem by creating logical merging neighborhoods where small cut cells are merged with their neighbors until the neighborhood volume is greater than a specified threshold. Unlike traditional cell-merging algorithms, a single cell can be in multiple neighborhoods.
As with FRD, the first step in advancing the solution in time is to create a 
conservative but possibly unstable update for each of the state variables using a finite volume formulation as in Eq. (\ref{eq:cons_up}).  The solution is then provisionally updated using the possibly unstable update, and this new-time provisional solution is modified by the state redistribution procedure. This procedure itself has two steps: first, the values of the provisionally updated state variables are used to define values associated with each neighborhood.  In the second step, the neighborhood values are used to define new solution values on the original cells.  We note that solution values in both cut and regular cells are used to define the neighborhood solutions, and that, as with FRD, the SRD procedure modifies the final solution in both cut and regular cells.  Unlike FRD however, the use of linear interpolation in this second step of SRD preserves linearity of the full solution where possible.

For a wide range of applications, it would be desirable to couple state redistribution with adaptive mesh refinement.  Current (W)SRD implementations impose a restriction that cells on or near the embedded boundary must all be at the same level of refinement; a coarse/fine boundary cannot intersect the embedded boundary.  Here we will introduce a ``re-redistribution" algorithm for WSRD that removes this restriction.  First, we review how the neighborhoods are defined and how the solution is advanced.

\subsection{Defining the neighborhoods}\label{sec:nbhds}

To create merging neighborhoods we first identify 
``small'' cells, i.e. cells with a volume fraction less than a specified target threshold, $V_{target}$.  These cells are then logically merged with neighboring cells until the neighborhood volume, i.e. the sum of the volume of all cells in the 
neighborhood reaches the target value. The neighborhoods in WSRD do not form a tessellation of the computational domain as is done in traditional cell merging.  Instead they are logical constructs where cells can potentially be part of several neighborhoods.  Here we use
\begin{align*}
    V_{target} = \Delta x \Delta y \Delta z / 2.
\end{align*}
This choice for the target threshold is informed by \cite{berger2021state}, but in practice the appropriate threshold may depend on factors such as the details of the numerical scheme and boundary conditions. 

There are several ways to choose which neighbor a small cell merges with. Using the simplest option, {\it normal merging}, we first merge small cells with their (cut or regular) neighbor in the direction of the boundary normal. If the addition of the neighbor's volume is not sufficient for the target threshold, additional cells are added to the merging neighborhood. In our implementation, our second choice neighbor is based on the next largest component of the normal. This creates an ``L-shaped'' neighborhood and we automatically merge the cell in the same plane as the two other neighbors to turn this neighborhood into a $2 \times 2$ neighborhood. In 3D if the threshold target is still not met, we add the additional cells necessary to create a $2 \times 2 \times 2$ neighborhood. An alternative process, called {\it central merging}, merges a small cell with all of its (cut or regular) neighbors in the $3 \times 3 \times 3$ box centered on the small cell. 

% we want to discuss the differences between the merging processes, explain why the smaller box is better, 

The main motivation for WSRD (as opposed to the initial formulation of SRD) is that a large cut cell that is just below the threshold volume requires less stabilization than a very small cell.  Observing that every cell is always in its own neighborhood, we define $M_{i,j,k}$ to be the set of cells contained in the neighborhood associated with cell $(i,j,k)$, and $M^-_{i,j,k} = M_{i,j,k} - \{(i,j,k)\}$, i.e. $M^-_{i,j,k}$ is the set of cells contained in the neighborhood of cell $(i,j,k)$ with the exception of cell $(i,j,k)$ itself. We define $W_{i,j,k}$ to be the set of
indices $(r,s,t)$ such that cell $(i,j,k)$ is in the neighborhood associated with $(r,s,t)$, i.e. so that cell $(i,j,k)$ is in $M_{r,s,t}.$ We define $W_{i,j,k}^- = W_{i,j,k} - \{(i,j,k)\}$.  Finally, we define $N_{i,j,k}$ to be the number of neighborhoods that cell $(i,j,k)$ is in, i.e. the size of $W_{i,j,k}$.

For the weighted version of SRD, at every regular or cut cell $(i,j,k)$ we define 
\begin{equation}\label{eqn:beta}
{\beta}_{i,j,k} =\begin{cases}
			(V_{target} - V_{i,j,k}) \;\; /  \displaystyle\sum_{(r,s,t) \in M^-_{i,j,k} } \, V_{r,s,t} \, & \text{if $V_{i,j,k} < V_{target}$}\\
            0 & \text{if $V_{i,j,k} \geq V_{target}$}
		 \end{cases}
\end{equation}
and then 
\begin{equation}
\label{eqn:alpha}
{\alpha}_{i,j,k} =  1 \, - \,  \frac{1}{N_{i,j,k}} \, \sum_{(r,s,t) \in W^-_{i,j,k} } \, \beta_{r,s,t} \, .
\end{equation}
%The idea of $\beta$ is to control the extent to which redistribution changes the solution. The closer $V_{i,j,k}$ is to $V_{target},$ the less the redistribution operation changes  the solution in $(i,j,k)$. 
Equations \eqref{eqn:beta} and \eqref{eqn:alpha} imply  that  $0 \le \alpha_{i,j,k}, \beta_{i,j,k} \le 1$.   We observe that $\beta_{i,j,k}$ depends on the volume fractions of the cells in the neighborhood of $(i,j,k),$ while $\alpha_{i,j,k}$ depends on the values of $\beta$ in cells whose neighborhoods contain $(i,j,k).$
%We note that there are a number of other strategies for defining the SRD weights including variants where the weights depend on the solution.  

As $V_{i,j,k}$
becomes larger, the weighted algorithm increases the dependence of the final solution in cell $(i,j,k)$ on the provisional solution in $(i,j,k)$, unlike the original SRD algorithm. 
This weighted SRD algorithm may have somewhat different stability properties than the original SRD algorithm, especially for one-dimensional test cases that are not representative of Cartesian cut cell meshes in higher dimensions.
 
%The formulae weight the contribution of cell $(i,j,k)$ itself to the quantities denoted 
%\MarginPar{HUH?}
%by $(i,j,k)$ by $\alpha$, and weight the contributions of cells $(r,s,t) \neq (i,j,k)$ by $\beta$.  
We define the weighted volume, $\widehat V_{i,j,k}$, of the neighborhood associated with $({i,j,k})$ by
 \begin{equation}
 \label{eqn:nbhd_vol}
{\widehat V}_{i,j,k} =  \alpha_{i,j,k} V_{i,j,k} +  \beta_{i,j,k} \sum_{(r,s,t) \in M^-_{i,j,k}} \,  \frac{V_{r,s,t}}{N_{r,s,t}}.
\end{equation}
and the centroid of the neighborhood associated with $(i,j,k)$ by 
\begin{align}
\begin{aligned}
({\widehat x}_{i,j,k},{\widehat y}_{i,j,k},{\widehat z}_{i,j,k}) =  \frac{1}{\widehat V_{i,j,k}} \biggl(
&\alpha_{i,j,k} V_{i,j,k} (x_{i,j,k},y_{i,j,k},z_{i,j,k}) \, + \\
& \; \beta_{i,j,k}  \sum_{(r,s,t) \in M^-_{i,j,k} }  \frac{V_{r,s,t}}{N_{r,s,t}}(x_{r,s,t},y_{r,s,t},z_{r,s,t}) \biggr).
\end{aligned}
\end{align}
where $(x_{i,j,k}, y_{i,j,k}, z_{i,j,k})$ is the original centroid of cell $(i,j,k)$.

In an efficient implementation, all of the steps in this subsection can be done when the EB information is defined at the start of a simulation or, in the case of AMR, after any step in which changes in the grid hierarchy would require it. 
%\MarginPar{Moved this from above}

\subsection{Advancing the solution}\label{sec:adv}

As noted in Sec. \ref{sec:math}, we start by defining a conservative update, $\delta U_{i,j,k}^c$ using 
Eq. (\ref{eq:cons_up}).  Here for simplicity we initially focus on a single step integration scheme, in particulate a piecewise-linear variant of Godunov's method of Collela \cite{colella1990multidimensional} (2D) and Saltzman \cite{SALTZMAN1994} (3D) that has been modified in and near cut cells.  The details of how we modify the algorithm are discussed in \ref{sec:EBGdnv}.  (Berger and Giuliani \cite{berger2021state} use a similar method in two dimensions that uses a somewhat different approach to dealing with cut cells.)
\commentout{
To advance the solution, we first compute a provisional, but possibly unstable cell update
 \begin{equation}
     \widehat{U}_{i,j,k}=U^n_{\ivec}-\frac{\Delta t}{V_{\ivec}}\displaystyle\sum_{\ell\in\text{faces}}\textbf{F}^*_{\ell}\cdot \textbf{n}_{\ell}A_{\ell}
 \end{equation}
 where $A_{\ell}$ is a face area, $\textbf{F}^*_{\ell}$ is the numerical flux, and $\textbf{n}_{\ell}$ is the outward pointing normal.
 }
 
 For WSRD, once we have computed $\delta U_{i,j,k}^c$ we define a provisional update
 \begin{equation}
     \widehat{U}_{i,j,k}=U^n_{\ivec}+ \Delta t \; \delta U_{i,j,k}^c
 \end{equation}
 
The next step is to define the  weighted solution average for each neighborhood as
 \begin{equation}
 \label{eqn:qhat_alpha}
\widehat{Q}_{i,j,k} =  \frac{1}{{\widehat V}_{i,j,k}} \, \left( \alpha_{i,j,k} V_{i,j,k} \widehat{U}_{i,j,k} \, + \beta_{i,j,k} \,
\sum_{(r,s,t) \in M^-_{i,j,k}}   
\frac{V_{r,s,t}}{N_{r,s,t}}  \widehat{U}_{r,s,t} \right).
 \end{equation}
 We observe that if $(i,j,k)$ is a regular cell, or a cut cell with $V_{\ivec} > V_{target},$ the neighborhood associated with $(i,j,k)$ contains only itself, and $\widehat{Q}_{i,j,k} = \widehat{U}_{i,j,k}.$

To achieve second order accuracy where possible, in both the original and weighted SRD approaches, we define a gradient, $(\widehat{\sigma}_x, \widehat{\sigma}_x, \widehat{\sigma}_x)$ in all neighborhoods.  To do so, we define a $3\times3\times3$ ($3\times3$ in 2D) box centered on $(i,j,k)$ for all $(i,j,k)$ whose neighborhoods contain at least two cells; otherwise we simply set the gradient to zero. (Note that this box has no relation to the size or shape of the merging neighborhoods themselves.)
We then use the least squares equation to fit the values of $\widehat{Q}$ in all cut or regular cells in this box, treating the value of each $\widehat{Q}_{r,s,t}$ as defined at the centroid of that neighborhood,
 $(\widehat{x}_{r,s,t},\widehat{y}_{r,s,t},\widehat{z}_{r,s,t})$.  
If the distance between centroids in any direction is below a given threshold (half the mesh spacing in that direction), then the gradient stencil is increased to have width $5$ rather than $3$ in that direction.  If the increased stencils either don't contain enough points or
don't have sufficient separation of coordinates, we set the gradient to 0.
In order to not generate new maxima or minima of $\widehat{Q},$ we also limit the components of the gradient with a Barth-Jesperson-style limiter. 
This allows us to define a linear function for neighborhood $(i,j,k)$ of the form
 \begin{equation}\label{eq:qrecon}
\widehat{q}_{i,j,k}(x,y,z) = \widehat{Q}_{i, j,k} + \widehat{\sigma}_{x,i,j,k}(x - \widehat{x}_{i,j,k}) + \widehat{\sigma}_{y,i,j,k}(y - \widehat{y}_{i,j,k}) + \widehat{\sigma}_{z,i,j,k}(z - \widehat{z}_{i,j,k}) .
\end{equation}
We then compute the final solution,
 \begin{equation} \label{eq:final_update_alpha}
        U^{n+1}_{i,j,k} =   \alpha_{i,j,k} \; \widehat{q}_{i,j,k}(x_{i,j,k},y_{i,j,k},z_{i,j,k}) + 
        \frac{1}{N_{i,j,k}}\sum_{(r,s,t)  \in W^-_{i,j,k}} \beta_{r,s,t} \; 
        \widehat{q}_{r,s,t}(x_{i,j,k},y_{i,j,k},z_{i,j,k})  .
\end{equation}
We observe that if a small cut cell $(i,j,k)$ is in no neighborhoods but its own, then $\alpha_{i,j,k} = 1$ and ${U}^{n+1}_{i,j,k} = \widehat{q}_{i,j,k}(x_{i,j,k},y_{i,j,k},z_{i,j,k})$ regardless of how many other cells are in its neighborhood. 
\\

%% file: 02_algorithms/srd_matrix.tex
\subsection{Matrix form}

We can also express the WSRD algorithm in matrix form, which will be useful in defining the re-redistribution algorithm.  We first define the matrix $A$ such that
\[
\widehat{V} = A V,
\]
where $V$ is the vector of cell volumes and $\widehat{V}$ is the vector of neighborhood volumes. This is the matrix form of Equation~(\ref{eqn:nbhd_vol}).
Thus, a nonzero entry $a_{p,q}$ of $A$ corresponds to the weighted contribution of cell-$p$ to neighborhood-$q$.  (Tying back to our previous notation, cell-$p$ might be $(i,j,k$) and cell-$q$ might be cell $(r,s,t)$ in 3D.)
We note that the matrix $A$ can be stored in a sparse format. Specifically, for each neighborhood we only need the values corresponding to the nonzero entries and an associated index. 
%\MarginPar{not a full sentence}

\commentout{
To make it concrete, for the case considered in Figure \ref{fig:3x3}, the matrix $A$ for the weighted algorithms is given by 
\MarginPar{This matrix corresponds to an image that we have not included}
\[
A = 
\begin{bmatrix}
1 & 0 & 0 & 0 & 0 & 0 & 0 & 0 & 0 \\
0 & {\alpha_2} & 0 & 0 & 0 & 0 & 0 & 0 & 0 \\
0 & \frac{\beta_3}{2} & 1 & 0 & 0 & 0 & 0 & 0 & 0 \\
0 & 0 & 0 & {\alpha_4} & 0 & 0 & 0 & 0 & 0 \\
0 & 0 & 0 & 0 &  {\alpha_5} & 0 & 0 & 0 & 0 \\
0 & 0 & 0 & 0 &  \frac{\beta_6}{4} & {\alpha_6} & 0 & 0 & 0 \\
0 & 0 & 0 & \frac{\beta_7}{2} & 0 & 0 & 1 & 0 & 0 \\
0 & 0 & 0 & 0 &  \frac{\beta_8}{4} & 0 & 0 & {\alpha_8} & 0 \\
0 & 0 & 0 & 0 &  \frac{\beta_9}{4} & \frac{\beta_9}{2} & 0 & \frac{\beta_9}{2} & 1 
\end{bmatrix}.
\]
In $A$ the entries $a_{ij}$ are the inverse of cell $j$'s overlap count, if cell $i$ is in cell $j$'s neighborhood, multiplied by the corresponding previously defined $\alpha$ and $\beta$.
}
The definition of the $\alpha$'s and the $\beta$'s guarantees that the nonzero elements of $A$ are positive.
An additional property of $A$ is that the columns sum to 1, i.e., 
\[
e^T A = e^T,
\]
where $e$ is a vector of all 1's. This property ensures that the WSRD algorithm is conservative.
Although not pursued here, this suggests that any collection of non-negative weights that reach the target neighborhood volumes with
$e^T A = e^T$ can be used to define the SRD algorithm.  This opens the possibility of more sophisticated approaches, including ones that might depend on the local solution.

We can now write Equation~(\ref{eqn:qhat_alpha}) as
\[
\widehat{Q} = \mathrm{Diag}(\widehat{V})^{-1} A \;\mathrm{Diag}(V) \widehat{U}.
\]
Furthermore, if the gradients are zeroed in Eq. (\ref{eq:qrecon}) then we have simply
\[
U^{n+1} = A^T \widehat{Q}.
\]
\commentout{Conservation then follows from
\begin{align*}
V^T U^{n+1} &= V^T A^T \; \mathrm{Diag}(\widehat{V})^{-1} A \; \mathrm{Diag}(V) \widehat {U}
= 
 \widehat{V}^T \; \mathrm{Diag}(\widehat{V})^{-1} A \; \mathrm{Diag}(V) \widehat {U}, \\
 &= e^T A \; \mathrm{Diag}(V) \widehat{U} = V^T \widehat{U} .
\end{align*}}
Including the gradients in Eq. (\ref{eq:qrecon}) redistributes mass to preserve linearity but does not alter the conservation properties of the method. Continuing, 
%Conservation is maintained because in Eq. (\ref{eqn:final_update_alpha}) 
%each $\widehat{q}_{r,s,t}$ can be written as
%\[
%\widehat{q}_{r,s,t} = \widehat{Q}_{r,s,t} + \delta \widehat{q}_{r,s,t}
%(x_{i,j,k},y_{i,j,k}, z_{i,j,k})
%\]
%where $ \delta \widehat{q}_{r,s,t} (x_{i,j,k},y_{i,j,k}, z_{i,j,k}) $ %represents a
%correction to $\widehat{Q}_{r,s,t}$ due to the reconstructed slopes.
%By construction {\marginpar \MJB{is this relevant? Only for linear? Matrix proof is good}}
%\[
% \sum_{i,j,k}  \alpha_{i,j,k} \delta %\widehat{q}_{i,j,k}(x_{i,j,k},y_{i,j,k},z_{i,j,k})
% +  \frac{1}{N_{i,j,k}}\sum_{(r,s,t)  \in W^-_{i,j,k}} \beta_{r,s,t} \; 
%       \delta \widehat{q}_{r,s,t}(x_{i,j,k},y_{i,j,k},z_{i,j,k}) = 0 \;\;\;.
%\]
we can write \eqref{eq:final_update_alpha}  as
\begin{align*}
U^{n+1} = A^T \widehat{Q} + \left[\mathrm{Diag}(x)A^T-A^T \mathrm{Diag}(\widehat{x}) \right ]\widehat{\sigma}_x
 &+\left[ \mathrm{Diag}(y)A^T-A^T \mathrm{Diag}(\widehat{y}) \right ]\widehat{\sigma}_y \\
 &+ \left[ \mathrm{Diag}(z)A^T-A^T \mathrm{Diag}(\widehat{z}) \right ]\widehat{\sigma}_z,
\end{align*}
where $(x,y,z)$ are the centroids of the original cells, $\widehat{x},\widehat{y},\widehat{z}$
are vectors of the components of the weighted centroids of the neighborhoods and $\widehat{\sigma}_x,\widehat{\sigma}_y, \widehat{\sigma}_x$ are vectors of the components of the gradients from the neighborhood gradient reconstruction discussed above.
%\MarginPar{are the sigmas matrices, vectors or numbers?   Should they have hats as in (12)?}

%As an aside, we note that
%\begin{align*}
%V^T \left[\mathrm{Diag}(x)A^T-A^T \mathrm{Diag}(\widehat{x}) \right ]\sigma_x  &=
%\left [ x^T \mathrm{Diag}(V)A^T- V^T A^T \mathrm{Diag}(\widehat{x}) \right ]\sigma_x %\\
%&= \left [ \widehat{x}^T \mathrm{Diag}(\widehat{V}) - \widehat{V}^T \mathrm{Diag}%(\widehat{x}) \right ] \sigma_x = 0 ,
%\end{align*}
%and analogously for $\sigma_y$ and $\sigma_z.$

\subsection{Coupling to AMR}\label{sec:coupling_to_amr}
\vspace{.1in}

As discussed earlier, when cut cells are sufficiently near to the coarse/fine boundary,
%in a traditional block-structured AMR algorithm, the data on the coarse level, $\ell$, are advanced independently of the data on the fine level and then the data on the fine level, $\ell+1$, are advanced, possibly with subcycling in time. When the coarse and fine solutions reach the same time, they need to be synchronized to form a composite solution. When no cut cells are near a coarse/fine boundary, 
%the synchronization between levels is simply the average down and refluxing operations discussed in}} %Section~\ref{sec:refluxing}.
%However, in the vicinity of the intersection of the fluid-body interface with the coarse/fine boundary,
the synchronization operations for AMR must include {\it refluxing the redistribution} as well as {\it redistributing the reflux}
%also account for the effect of redistribution on the solution and the refluxing algorithm defined in Section~\ref{sec:refluxing} must be modified to account for the geometry.
in order to define a conservative and stable composite solution.  Borrowing some of the notation from \cite{pember1995adaptive}, which described re-redistribution for algorithms using flux redistribution, here we present the re-redistribution algorithm for algorithms using weighted state redistribution.

We first briefly introduce some additional notation that will be helpful:  fine {\it valid} cells are cells inside level $(\ell+1)$ grids, and fine {\it ghost} cells are fine cells that are not inside level $(\ell+1)$ grids but are adjacent to fine valid cells.   Fine ghost cells overlay coarse uncovered cells at the coarse/fine boundary.

We recall that in order to perform the standard refluxing operation at the coarse/fine boundary we accumulate the difference of coarse and fine fluxes in face-centered
$\delta F$'s as shown in Eq.~\ref{eq:reflux}.  In order to reflux the redistribution, we analogously compute cell-centered $\delta R$'s that track the extensive amount of each variable that effectively crosses the coarse/fine interface in the redistribution step at each level.  Because the data associated with a fine ghost cell or a coarse covered cell is not part of the composite solution, contributions to or from those cells are not correctly accounted for in the composite solution. When redistributing at the coarse level, $\ell$, we include in $\delta R$ the extensive contributions that went 
%There are four quantities that need to be captured in the re-redistribution step, both for flux and state redistribution. 
%Specifically, we need to capture the extensive contribution
\begin{enumerate}[label=(\Alph*)]
    \item from coarse uncovered cells to coarse covered cells %, denoted $\delta R^{\ell,u \rightarrow c}$
    \item from coarse covered cells to coarse uncovered cells %, denoted $\delta R^{\ell,c \rightarrow u}$ 
\end{enumerate}
When redistributing at the fine level, $\ell+1,$ we include in $\delta R$  the extensive contributions that went  
\begin{enumerate}[label=(\alph*)]    
    \item from fine valid cells to fine ghost cells %, denoted $\delta R^{\ell+1,v \rightarrow g}$
    \item from fine ghost cells to fine valid cells %, denoted $\delta R^{\ell+1,g \rightarrow v}$
\end{enumerate}

%All of these contributions correspond to transfer of a quantity across the coarse/fine boundary that is, in some sense, lost when we construct the composite solution.  As an example, if density is the conserved quantity being updated, then the quantities captured in $\delta R$ have units of mass. 

\subsubsection{Computing $\delta R$}
\vspace{.1in}

Here we describe how to compute $\delta R$ for the WSRD algorithm. To do so, we first create a matrix, $\mathcal{R},$ that captures the contributions of $\widehat{U}$ to $U^{n+1}$ in nearby cells.

If we set the gradient terms in Eqs. (\ref{eq:qrecon})-(\ref{eq:final_update_alpha}) to 
zero, 
from the matrix form of WSRD we can write the final update (in terms of $\widehat{U}$) as
\[
\Diag(V) U^{n+1} = \Diag(V)  A^T \mathrm{Diag}(\widehat{V})^{-1} A \;\mathrm{Diag}(V) \widehat{U}
\]
From this relation we can see that the matrix
\[
\mathcal{R} \equiv \Diag(V)  A^T \mathrm{Diag}(\widehat{V})^{-1} A \;\mathrm{Diag}(V) \Diag(\widehat{U})
\]
specifies how values of $\widehat{U}$ affect the final update.  In particular, the row of $\mathcal{R}$ corresponding to $(I,J,K)$ specifies, in extensive form, the contribution of each of the $\widehat{U}$'s to $U^{n+1}_{I,J,K}$.  We note that $\mathcal{R}$ is a sparse matrix with the same sparsity pattern as $A^T A$, making the nonzeros easy to identify using the sparse representation of $A$.

Let us consider the contributions from the redistribution operation at each level. 
We can define $\mathrm{R}_{I,J,K}$ and $\mathrm{C}_{I,J,K}$ to be the row and column of  in $\mathcal{R}$, respectively, corresponding to cell $(I,J,K)$.
At the coarse level, we associate with uncovered cell $(I,J,K)$ the sum of its contributions to covered cells minus the sum of contributions of covered cells to $(I,J,K)$; i.e., for uncovered cell $(I,J,K)$ we define
%\[
%\delta R_{I,J,K}^\ell = \sum_{c\epsilon G(u)} \delta R_{I,J,K}^{\ell,u \rightarrow c} -
%\sum_{c\epsilon S(u)} \delta R_{I,J,K}^{\ell,c \rightarrow u}
%\]
%where $G(u)$ is the set of cells that u gives information to $S(u)$ is the set of cells that give information to $u$.
\begin{equation}
\delta R_{I,J,K}^\ell = \sum_{(II,JJ,KK)\epsilon \mathcal{C}^\ell} \mathcal{R}^\ell_{\mathrm{R}_{II,JJ,KK},\mathrm{C}_{I,J,K} }
 - \; \mathcal{R}^\ell_{\mathrm{R}_{I,J,K},\mathrm{C}_{II,JJ,KK} }
%\sum_{II,JJ,KK \eps \mathcal{U \delta R_{I,J,K}^{\ell,c \rightarrow u}
\label{eq:dr_srd_coarse}
\end{equation}
where $\mathcal{C}^\ell$ is the set of covered cells at level $\ell$ and $\mathcal{R}^\ell$ is matrix $\mathcal{R}$ at level $\ell$.
Here the first and second terms on the right hand side correspond to contributions of type (a) and type (b), respectively.
%where $G(u)$ is the set of cells that u gives information to $S(u)$ is the set of cells that give information to $u$.
Similarly, at the fine level, we associate with ghost cell $(i,j,k)$ the sum of contributions to cell $(i,j,k)$ from valid cells minus the sum of contributions of cell $(i,j,k)$ to valid cells; i.e.,
%\[
%\delta R_{i,j,k}^{\ell+1} = \sum_{v\epsilon S(g)} \delta R_{i,j,k}^{\ell,v \rightarrow g} -
%\sum_{v\epsilon G(g)} \delta R_{i,j,k}^{\ell,g \rightarrow v} \;\;.
%\]
\begin{equation}
\delta R_{i,j,k}^{\ell+1} = \sum_{(ii,jj,kk)\epsilon \mathcal{V}^{\ell+1}} \mathcal{R}^{\ell+1}_{\mathrm{R}_{i,j,k},\mathrm{C}_{ii,jj,kk} }
 - \; \mathcal{R}^{\ell+1}_{\mathrm{R}_{ii,jj,kk},\mathrm{C}_{i,j,k} }
%\sum_{II,JJ,KK \eps \mathcal{U \delta R_{I,J,K}^{\ell,c \rightarrow u}
\label{eq:dr_srd_fine}
\end{equation}
where $\mathcal{V}^{\ell+1}$ is the set of valid cells at level $\ell+1$ and $\mathcal{R}^{\ell+1}$ is the matrix $\mathcal{R}$ at level $\ell+1$.
In this case, the first and second terms on the right hand side correspond to contributions of type (c) and type (d), respectively.
We note that in both Eq. (\ref{eq:dr_srd_coarse}) and Eq. (\ref{eq:dr_srd_fine}), $\mathcal{R}$ is sparse with a known sparsity pattern, making the summations easy to compute.

We now need to describe how to treat the gradient terms in 
Eqs. (\ref{eq:qrecon})-(\ref{eq:final_update_alpha}). One possible approach would be to capture explicitly how the gradient depends on the
$\widehat{U}$'s; however, limiting makes this a nonlinear process and the potentially large stencil for the gradient dramatically enlarges the stencil for the redistribution.
Alternatively, we will view the gradients in Eq. (\ref{eq:final_update_alpha})
as fixed and then let
\[
\delta \widehat{q}^{I,J,K}_{R,S,T} = 
\widehat{\sigma}_{x,R,S,T}(x_{I,J,K} - \widehat{x}_{R,S,T}) + 
\widehat{\sigma}_{y,R,S,T}(y_{I,J,K} - \widehat{y}_{R,S,T}) + 
\widehat{\sigma}_{z,R,S,T}(z_{I,J,K} - \widehat{z}_{R,S,T})
\]
be the contribution of neighborhood $(R,S,T)$ to the final update of cell $(I,J,K)$ as in Eq. (\ref{eq:final_update_alpha}).
The issue here is that we cannot explicitly link the correction back to the $\widehat{U}$'s.
One could potentially associate this update directly with the neighborhood.  
%This would effectively treat $\delta q^{I,J,K}_{R,S,T}$ as a correction to $\delta R^\ell_{I,J,K}$ from cell $R,S,T$.  
This would correspond to treating the entire $\delta q$ correction for neighborhood $(R,S,T)$ as originating from cell $(R,S,T)$. Although such an approach would be conservative, it could potentially generate large corrections for cells with very small volumes, which would introduce an instability.
Instead we define
\[
\left (\delta q^{I,J,K}_{R,S,T} \right )_{II,JJ,KK} =
\frac{1}{\widehat{V}_{R,S,T}} A_{\mathrm{R}_{R,S,T},\mathrm{C}_{II,JJ,KK}} \; V_{II,JJ,KK}
\; \delta \widehat{q}^{I,J,K}_{R,S,T} \quad ,
\]
which distributes $ \delta \widehat{q}^{I,J,K}_{R,S,T}$ over the cells $(II,JJ,KK)$ that are in neighborhodd $(R,S,T)$.
If $(I,J,K)$ is uncovered and $(II,JJ,KK)$ is covered then we add
\begin{equation}
V_{II,JJ,KK} A_{\mathrm{R}_{R,S,T},\mathrm{C}_{II,JJ,KK}}
\left (\delta q^{II,JJ,KK}_{R,S,T} \right )_{I,J,K}
-V_{I,J,K} A_{\mathrm{R}_{R,S,T},\mathrm{C}_{I,J,K}}
\left (\delta q^{I,J,K}_{R,S,T} \right )_{II,JJ,KK}
\label{eq:dq_cor}
\end{equation}
to
$\delta R_{I,J,K}^\ell$.  Note that in Eq. (\ref{eq:dq_cor}) the indices of $A$ correspond to multiplication by $A^T$.
Analogous corrections are also added to $\delta R_{i,j,k}^{\ell+1}$.

Similar to FRD, we associate the WSRD redistribution correction with the uncovered coarse cells at the coarse/fine boundary. 
Specifically, for $(I,J,K)$ at the coarse/fine boundary, we define 
\begin{equation*}
\delta \mathbf{R}^\ell_\IVEC = \delta R^{\ell}_\IVEC 
+ \sum_{i,j,k} \delta R^{\ell+1}_\ivec
%+ \sum_{i,j,k \, | \, C^{\ell+1}_\ivec \, \subset \, C^\ell_{II,JJ,KK}} \delta R^{\ell+1}_\ivec
%\label(eq:delR}
\end{equation*}
where the sum is over fine grid ghost cells $(i,j,k)$  that overlay cell $(\IVEC)$.  
The quantity $\delta \mathbf{R}^\ell_\IVEC $ represents the difference between redistribution on the coarse and fine levels in extensive form.

\subsubsection{Updating the solution}\label{sec:final_update}
\vspace{.1in}

Once we have computed $\delta \mathbf{R}^\ell_\IVEC,$
if $\vfrac_\IVEC = 1$ then we can simply update
\[
U^\ell_\IVEC := U^\ell_\IVEC + \frac{1}{\vol_\IVEC} \,\delta \mathbf{R}^\ell_\IVEC
\]
to complete the synchronization so that the resulting scheme is conservative.
However, if $\vfrac_\IVEC < 1$ adding the entire increment $\delta R^\ell_\IVEC$ to $U^\ell_\IVEC$ could introduce an instability.  
Instead the update needs to be redistributed as well.  This redistribution can be done using either the methodology used in FRD or using a WSRD approach.  Here we describe the FRD approach for this ``redistribution of reflux" step.  Specifically, we update 
\[
U^\ell_\IVEC := U^\ell_\IVEC + \frac{\vfrac_\IVEC}{\vol_\IVEC} \,\delta \mathrm{R}^\ell_\IVEC
\]
and redistribute the remainder, $(1-\vfrac_\IVEC)  \,\delta \mathbf{R}^\ell_\IVEC$ to neighboring cells (neighboring in the FRD sense, see Appendix A; not a WSRD neighborhood). In particular we define
\[
\vol_{nbh} = \sum_{II,JJ,KK \, \in \, \NFRD(\IVEC) } \vol_{II,JJ,KK}
\]
where $\NFRD(I,J,K)$ is the set of regular and cut cells $(II,JJ,KK)$ such that none of the indices differs from $(II,JJ,KK)$ by more than one and where $(II,JJ,KK)$ can be reached from $(I,J,K)$ by a monotone path of regular and cut cells in index 
space; i.e., a path in which none of the indices are both incremented and decremented.

We then update cells $(II,JJ,KK)$ in $\NFRD(I,J,K)$ by
\[
U^\ell_{II,JJ,KK} := U^\ell_{II,JJ,KK} + \frac{(1-\vfrac_\IVEC)  \,\delta \mathbf{R}^\ell_\IVEC}{\vol_{nbh}}
\]
This redistribution process can update covered coarse cells.  In that case the correction to the covered coarse cell must be interpolated to the fine cells that cover the coarse grid cell. In particular, if
$(II,JJ,KK)$ is covered, we use piecewise constant interpolation to update
\[
U^m_{i,j,k} := U^m_{ii,jj,kk} + \frac{(1-\vfrac_\ivec)  \,\delta \mathbf{R}^\ell_\ivec}{\vol_{nbh}} \qquad \forall  \, \ivec \,| \,  C^m_\ivec \subset C^\ell_{II,JJ,KK} \; \mathrm{for} \;  m>\ell 
\]
at all cells $(i,j,k)$ at level $(\ell+1)$ that overlay cell $(II,JJ,KK)$ at level $\ell.$
%where $\mathbf{C}^\ell$ is the set of covered coarse cells at level $\ell$.
%This update corresponds to piecewise constant interpolation of $\frac{(1-\vfrac_\IVEC)  \,\delta R^\ell_\IVEC}{\vol_{nbh}}$ to the finer grid.

For cut cells ($\vfrac_\ivec < 1$), adding the entire flux correction from Eq. (\ref{eq:reflux}) can also introduce an instability.  Consequently, reflux increments also need to be redistributed. Here we have again used the FRD style redistribution for the reflux increments.  We note that, operationally, this can be done by adding the reflux corrections to $\delta \mathrm{R}$ before redistributing it.

\vspace{0.1in}
\noindent \textit{Method-of-lines temporal integration}
\vspace{0.1in}

The methodology described above can be easily extended to work using a method-of-lines discretization.  Berger and Giuliani \cite{berger2021state} show how to use state redistribution in the context of a simple predictor / corrector temporal discretization. Here we will derive essentially the same scheme but from a different perspective.  The alternate perspective provides useful insight into the requirements for the synchronization.

The WSRD process can be viewed as defining an operator, $R^{WSRD}$ such that
\begin{equation}\label{eq:op}
U^{n+1} = R^{WSRD} ( U^n + \Delta t \; \delta U ^c )
\end{equation}
In the absence of slopes (Eq. (\ref{eq:qrecon})) 
\[
R^{WSRD} = A^T \Diag(\widehat{V})^{-1} A \Diag(V)
\]
Starting with Eq. (\ref{eq:op}) we define an effective $\delta U^{n,WSRD}$ using
\begin{align}\label{eq:molform}
U^{n+1} &= R^{WSRD} ( U^n + \Delta t \; \delta U^c ) \nonumber \\
&=U^n + \Delta t \; \frac {R^{WSRD} ( U^n + \Delta t \:\delta U^c ) - U^n}{\Delta t} \nonumber \\
&= U^n + \Delta t \; \delta U^{n,WSRD}
\end{align}
We note that, when the reconstructed slopes are not considered,
$R^{WSRD}$ is linear and we can rewrite
\[
\delta U^{n,WSRD} = R^{WSRD}(\delta U^c) + \frac{1}{\Delta t} (R^{WSRD}-I)U^n
\]
The properties of $A$ show that the second term on the right hand side is an (anisotropic) numerical diffusion term with magnitude that scales with 
$\frac{1}{\Delta t}$.

We can then define the predictor / corrector scheme as
\begin{align}\label{eq:pcor}
U^* &= U^n +  \Delta t \; \delta U^{n,WSRD} \nonumber \\
U^{n+1} &= \frac{U^n + U^* + \Delta t \; \delta U^{*,WSRD}}{2} \nonumber \\
&= U^n + \Delta t \; \frac{\delta U^{n,WSRD}+\delta U^{*,WSRD}}{2}
\end{align}
This form shows that in forming $\delta R_{I,J,K}^\ell$ and $\delta R_{i,j,k}^{\ell+1}$ using Eqs. (\ref{eq:dr_srd_coarse}) and (\ref{eq:dr_srd_fine}), respectively, the contributions need to be multiplied by $\half$ to reflect their weightings in Eq. (\ref{eq:pcor}).

This basic approach can be extended to higher-order temporal integration strategies with suitable weights in defining the $\delta R$ that reflect the underlying temporal quadrature rule.  (We note that for higher-order temporal discretizations in conjunction with AMR and subcyling, special care must be taken in defining boundary conditions for fine grid from coarse data to avoid order reduction.  See, for example, \cite{EmmettETAL2019}).

%% file: 04_results/results_2d.tex
\subsection{Conservation properties}\label{sec:num_res_conservation}
\begin{figure}
     \centering
     \begin{subfigure}[b]{0.48\textwidth}
         \centering
         \includegraphics[width=\textwidth,page=1]{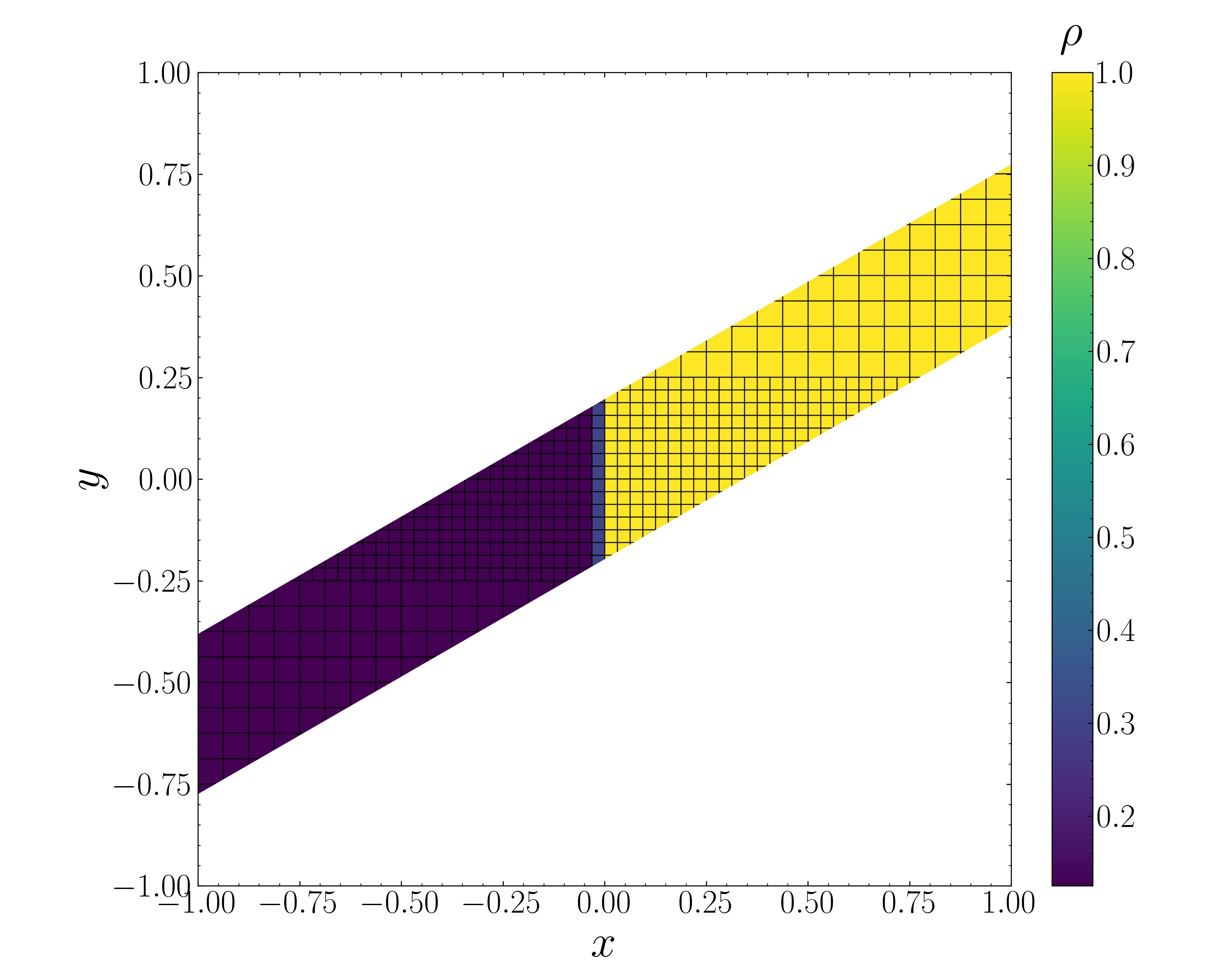}
         \caption{$t=0.0$.}
         %\label{fig:y equals x}
     \end{subfigure}
     \hfill%
     \begin{subfigure}[b]{0.48\textwidth}
         \centering
         \includegraphics[width=\textwidth,page=2]{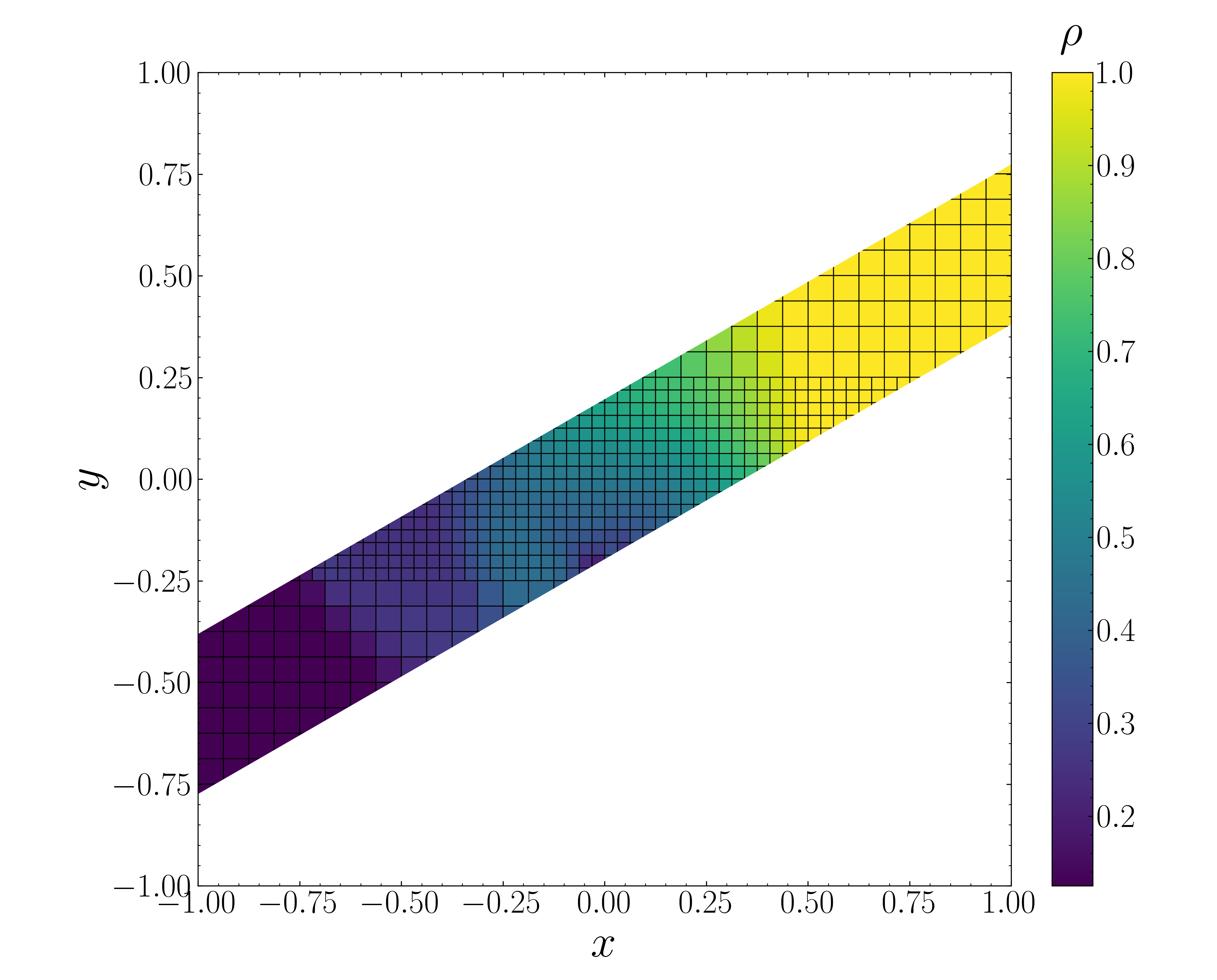}
         \caption{$t=0.4$.}
         %\label{fig:y equals x}
     \end{subfigure}
     \caption{Density profiles for the two-dimensional simulation of flow in a rotated channel with fixed mesh refinement.}
     \label{fig:cons_soln}
\end{figure}

We first demonstrate that finite volume embedded boundary methods  without the re-redistribution step violate conservation when the coarse-fine interface intersects the embedded boundary, and that the re-redistribution step restores conservation.
For this demonstration we choose very simple cases.  In the three-dimensional case, the fluid region is the interior of a cylinder whose axis passes through $(0,0,0)$ in a [-2:2,-2:2,-2:2] domain and is 30$^\circ$ off the x-axis in the y direction.  In the two-dimensional case we take the intersection of that cylinder with the $z=0$ plane, forming a fluid region bounded above and below by lines slanted at 30$^\circ$ off the x-axis in the y-direction.  The radius of the cylinder is $0.172$. The initial data in both cases is zero velocity with a discontinuous one-dimensional profile for density and pressure, namely $\rho = 0.125$ and $p = 0.1$ for $x \leq 0$ and $\rho = 1.0$ and $p = 1.0$ for $x > 0.$ The region of the domain
define by $-2 \leq x \leq 2$, $-0.125 \leq y \leq 0.125$ and, in 3D, $-2 \leq z \leq 2$
%[-2:2, -.125:.125, -2:2] 
is refined by a factor of two relative to the base level. Figure~\ref{fig:cons_soln} illustrates the solution at times $t=0$ and $t=0.4$, the bounding geometry and region of grid refinement.
Here we have intentionally presented results from a low resolution simulation to emphasize the role of the cut cells. When adequately resolved, at early times the initial discontinuity evolves into a rarefaction wave expanding to the right, a contact discontinuity and a shock wave propagating to the left.  The shock reflects off the upper boundary and generates an expansion wave at the lower boundary. Later in the evolution, the different waves begin to interact with each other, producing a fairly complex flow even in this relatively simple geometry.

\begin{figure}
     \centering
     \begin{subfigure}[b]{0.48\textwidth}
         \centering
        \includegraphics[width=\textwidth,page=1]{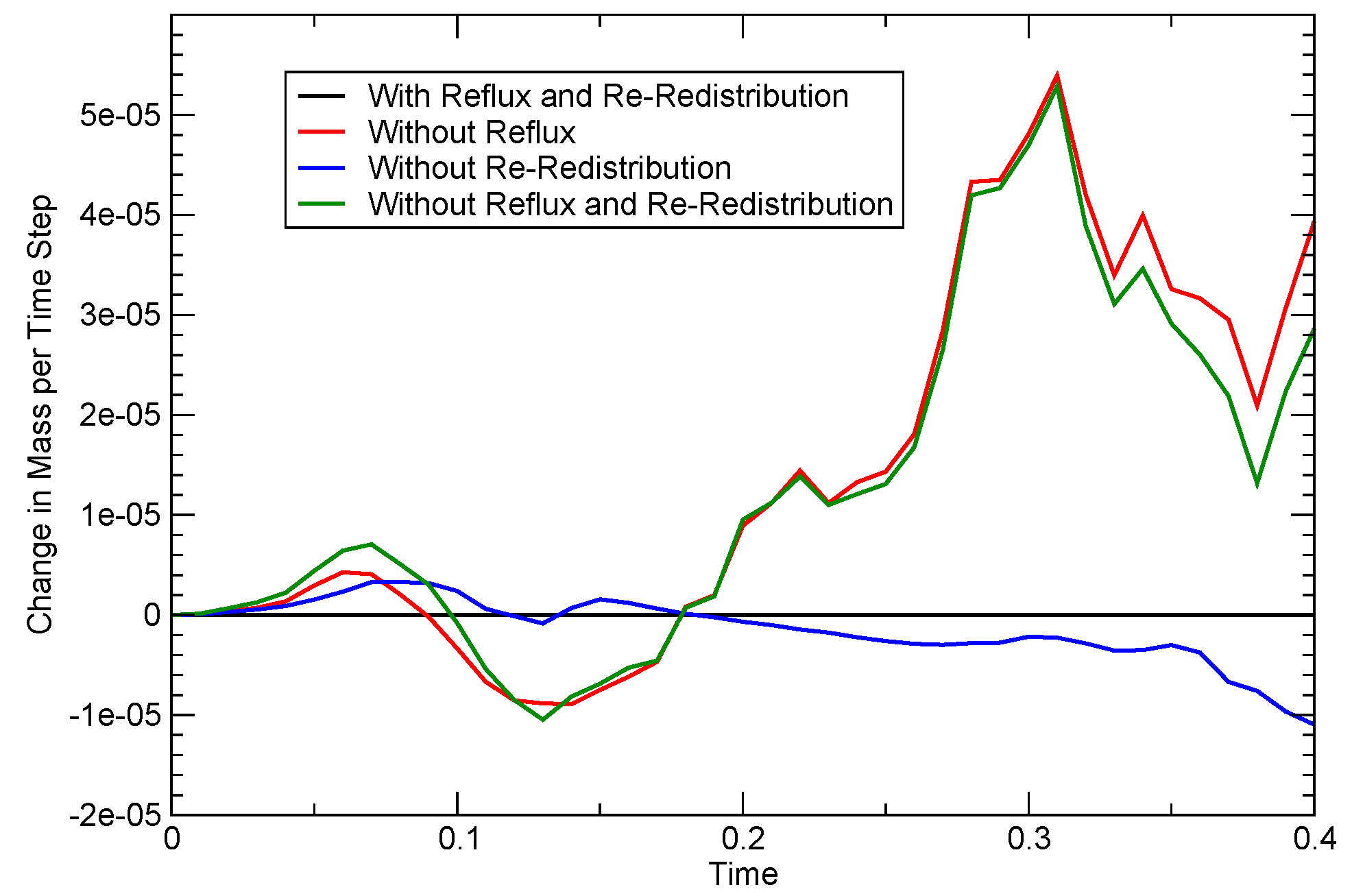}
         \caption{Two-dimensional simulation with Godunov methodology and state redistribution.}
         %\label{fig:y equals x}
     \end{subfigure}
     \hfill%
     \begin{subfigure}[b]{0.48\textwidth}
         \centering
       \includegraphics[width=\textwidth,page=2]{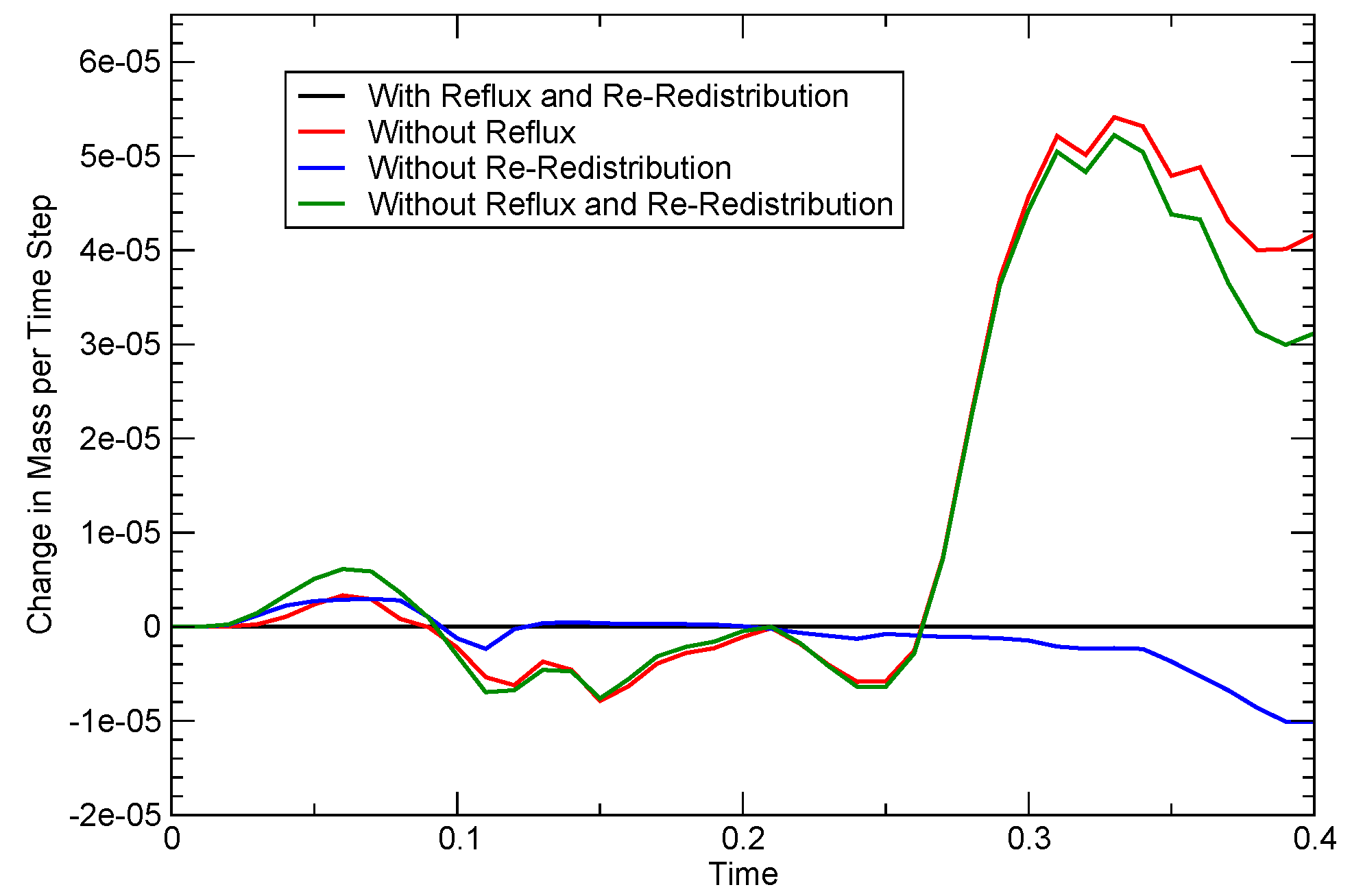}
         \caption{Two-dimensional simulation with method-of-lines methodology and state redistribution.}
         %\label{fig:y equals x}
     \end{subfigure}
     
\vspace{.2in}
     \begin{subfigure}[b]{0.48\textwidth}
         \centering
        \includegraphics[width=\textwidth,page=3]{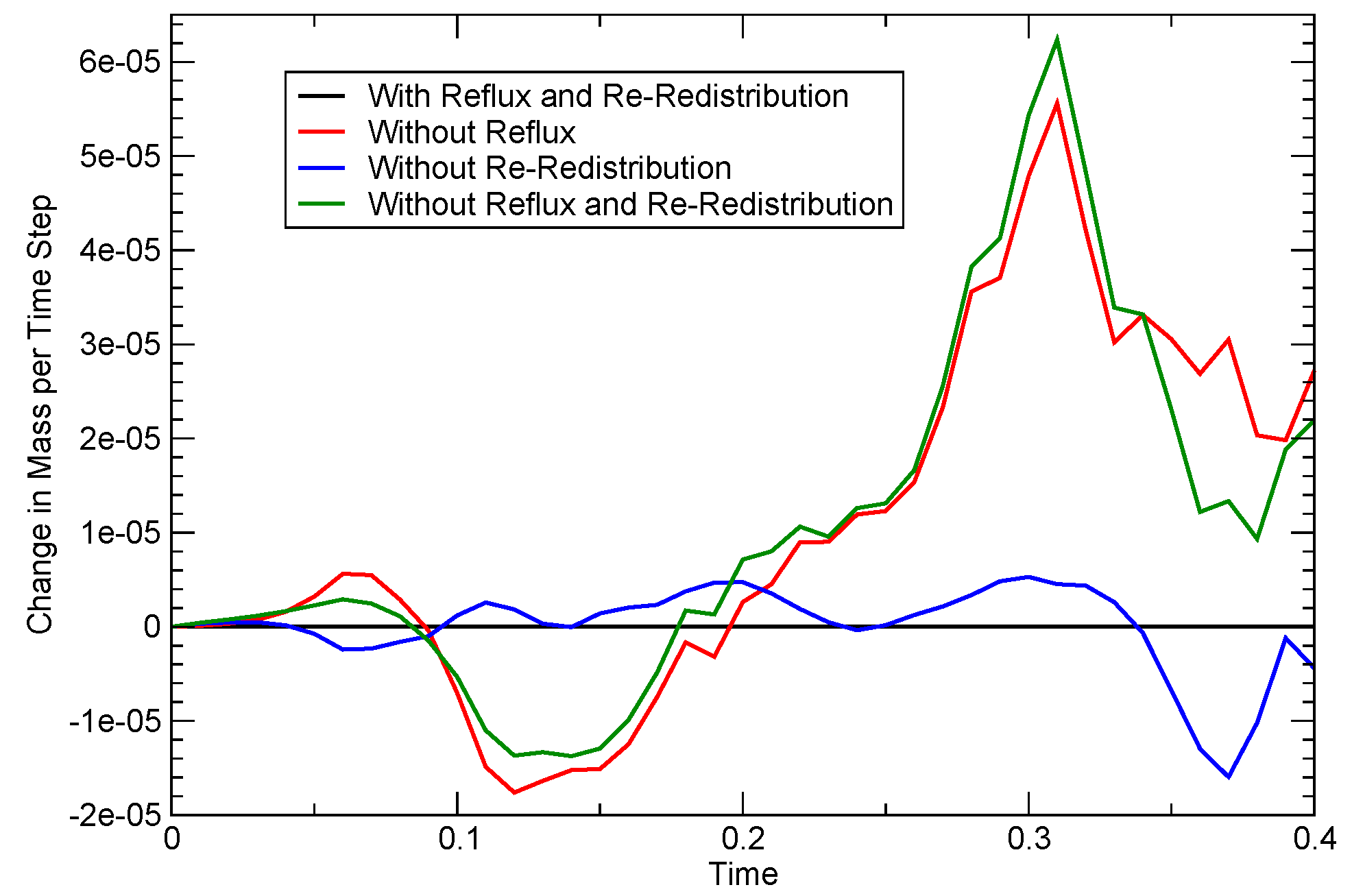}
         \caption{Two-dimensional simulation with Godunov methodology and flux redistribution.}
         %\label{fig:y equals x}
     \end{subfigure}
     \hfill%
     \begin{subfigure}[b]{0.48\textwidth}
         \centering
       \includegraphics[width=\textwidth,page=4]{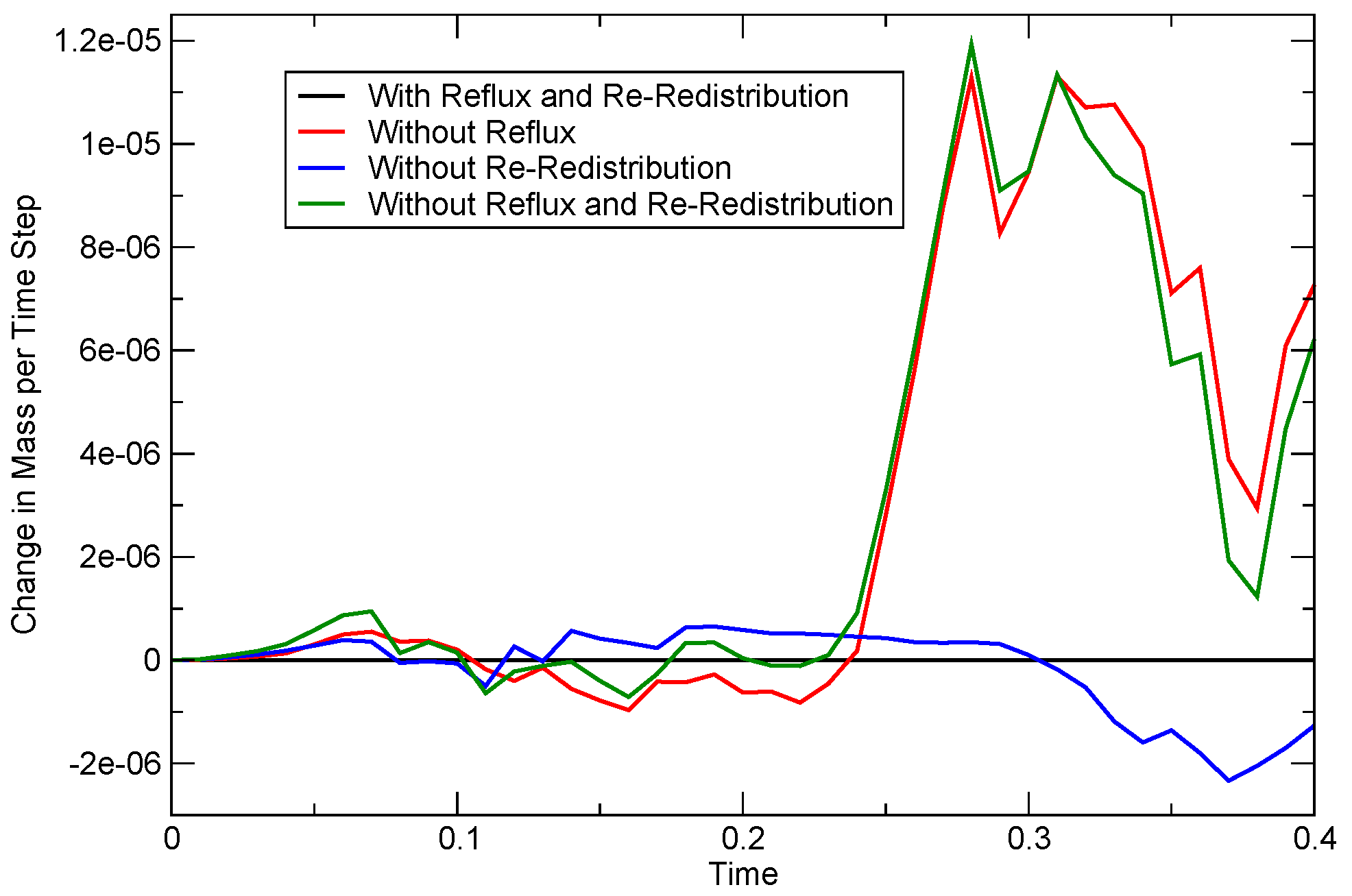}
         \caption{Three-dimensional simulation with Godunov methodology and state redistribution.}
         %\label{fig:y equals x}
     \end{subfigure}
     \caption{Change in mass per timestep as function of time}
     \label{fig:cons_mass}
\end{figure}

Figure~\ref{fig:cons_mass} shows the change in mass per time step in the entire domain.  It is hard to scale this plot in a meaningful way, so here we compare the rate of mass loss that would occur without refluxing on the same plot as the rate of mass loss due to not re-redistributing.   In Figure~\ref{fig:cons_mass}(a) we show two-dimensional results using the Godunov methodology and state redistribution.  In Figure~\ref{fig:cons_mass}(b) we show two-dimensional results using the method-of-lines methodology and state redistribution in order to demonstrate the lack of sensitivity of the results to the specific advection scheme, and that state re-redistribution behaves correctly for a multi-step integrator.  In Figure~\ref{fig:cons_mass}(c) we show two-dimensional results using the Godunov methodology and flux redistribution to compare the behavior of flux re-redistribution and state re-redistribution.  And finally in Figure~\ref{fig:cons_mass}(d) we show the three-dimensional results using the Godunov methodology and state redistribution.

We note that the flux of mass per time step across the coarse/fine interface over the time shown here ranges from roughly $10^{-7}$ to $10^{-3}$ in two dimensions and $10^{-8}$ to $4.\times 10^{-4}$ in three dimensions.  Thus we can see that the synchronization accounts for at most a few percent of the total mass crossing the coarse/fine interface.  It should also not be surprising that the effect of refluxing is greater than that of re-redistribution; the refluxing corrections occur over the entire coarse/interface while the re-redistribution only affects the solution at the part of the coarse/fine interface that is near the embedded boundary.  We also note, by contrasting the solutions with and without the synchronization steps, that the difference in solutions with and without the synchronization steps is also on the order of a few percent in max norm.

We mention a few caveats here.  This problem is relatively simple and run for a relatively short time.  For more complicated problems run for longer, it is quite possible for the methodology to fail in the absence of either refluxing or the redistribution synchronization.   In addition, recall that the final step of the synchronization methodology is the redistribution of the correction after it is added to the coarse solution.  It is straightforward to demonstrate that if the cylinder is shifted so that there are uncovered coarse cells with sufficiently small volume fractions (e.g., less than $10^{-4}$ in this case) at the coarse/fine boundary, the methodology fails in the first time step with negative energy and pressure, due to too large an update being added to that cell in the synchronization step.

Finally, we note that while we have shown the behavior of the algorithm in terms of mass conservation, the re-redistribution methodology acts on all conserved variables, thus maintaining the appropriate conservation properties of the entire solution.

%% file: 04_results/pele-sod.tex
\subsection{Sod shock tube validation}\label{sec:sod-validation}

Using the same two-dimensional geometry described in Section~\ref{sec:num_res_conservation}, we perform a validation test case using the Sod shock tube problem.  For this problem  the initial discontinuity is normal to the centerline of the channel, i.e., the left and right states are defined with respect to $x'>0$ where $x' = x \cos(\theta) + y \sin(\theta)$ and $\theta$ is the angle of the channel. This case is similar to validation case considered by Gulizzi {\it et al.} ~\cite{GULIZZI2022}. The AMR refinement is based on a dynamic tagging criteria: additional levels are added where the local change in density between adjacent cells is greater than 0.05. Simulations were performed with one level of AMR and three base grid resolutions: $\Delta x_0 = 0.0125,  0.00625,$ and $0.003125$. Coarse-fine interfaces intersect the EB in this problem, and these intersections are dynamic, i.e., they move depending on the flow features. These simulations used the method-of-lines methodology and the proposed state re-redistribution.

\begin{figure}
     \centering
     \begin{subfigure}[b]{0.48\textwidth}
         \centering
         \includegraphics[width=\textwidth,page=1]{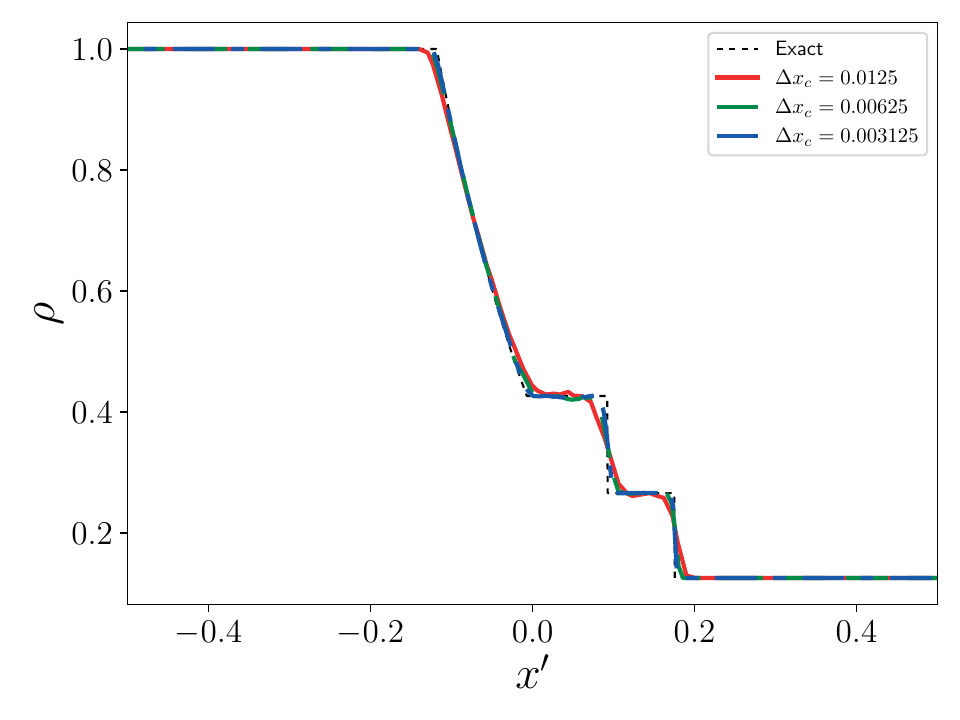}
         \caption{Density}
         %\label{fig:y equals x}
     \end{subfigure}
     \hfill%
     \begin{subfigure}[b]{0.48\textwidth}
         \centering
         \includegraphics[width=\textwidth,page=2]{04_results/pele-sod/profiles.pdf}
         \caption{Pressure}
         %\label{fig:y equals x}
     \end{subfigure}\\
     \begin{subfigure}[b]{0.48\textwidth}
         \centering
         \includegraphics[width=\textwidth,page=3]{04_results/pele-sod/profiles.pdf}
         \caption{Velocity}
         %\label{fig:y equals x}
     \end{subfigure}
     \hfill%
     \begin{subfigure}[b]{0.48\textwidth}
         \centering
         \includegraphics[width=\textwidth,page=4]{04_results/pele-sod/profiles.pdf}
         \caption{Temperature}
         %\label{fig:y equals x}
     \end{subfigure}
     \caption{Profiles for the Sod shock tube problem in a rotated channel at $t=0.1$.}
     \label{fig:pele-sod}
\end{figure}

The simulation results are presented in Figure~\ref{fig:pele-sod}. The data is taken through a line in the center of the cylinder. There is good agreement across all flow variables and mesh resolutions with the exact solution. As the base grid is refined, the numerical solution converges to the exact solution. As expected from previous results using the state re-redistribution, mass is conserved throughout the simulation.

\subsection{Flow past a cylinder}\label{sec:cyl-validation}

\begin{figure}
     \centering
     \begin{subfigure}[t]{0.32\textwidth}
         \centering
         \includegraphics[width=\textwidth,page=1]{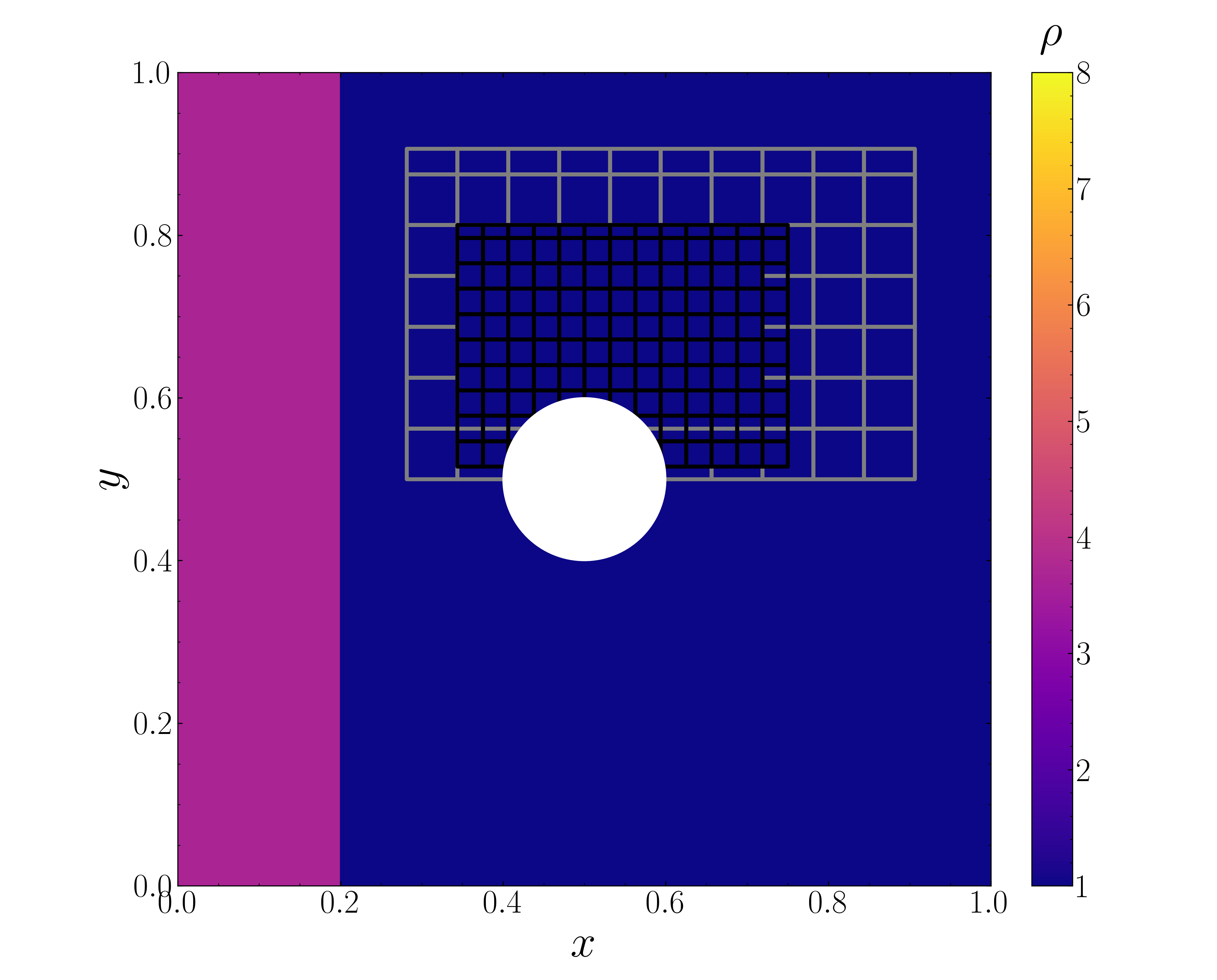}
         \caption{$t=0.0$, showing full domain and refinement levels.}
         \label{fig:cylinder-amr-grids}
     \end{subfigure}
     \hfill%
     \begin{subfigure}[t]{0.32\textwidth}
         \centering
         \includegraphics[width=\textwidth,page=2]{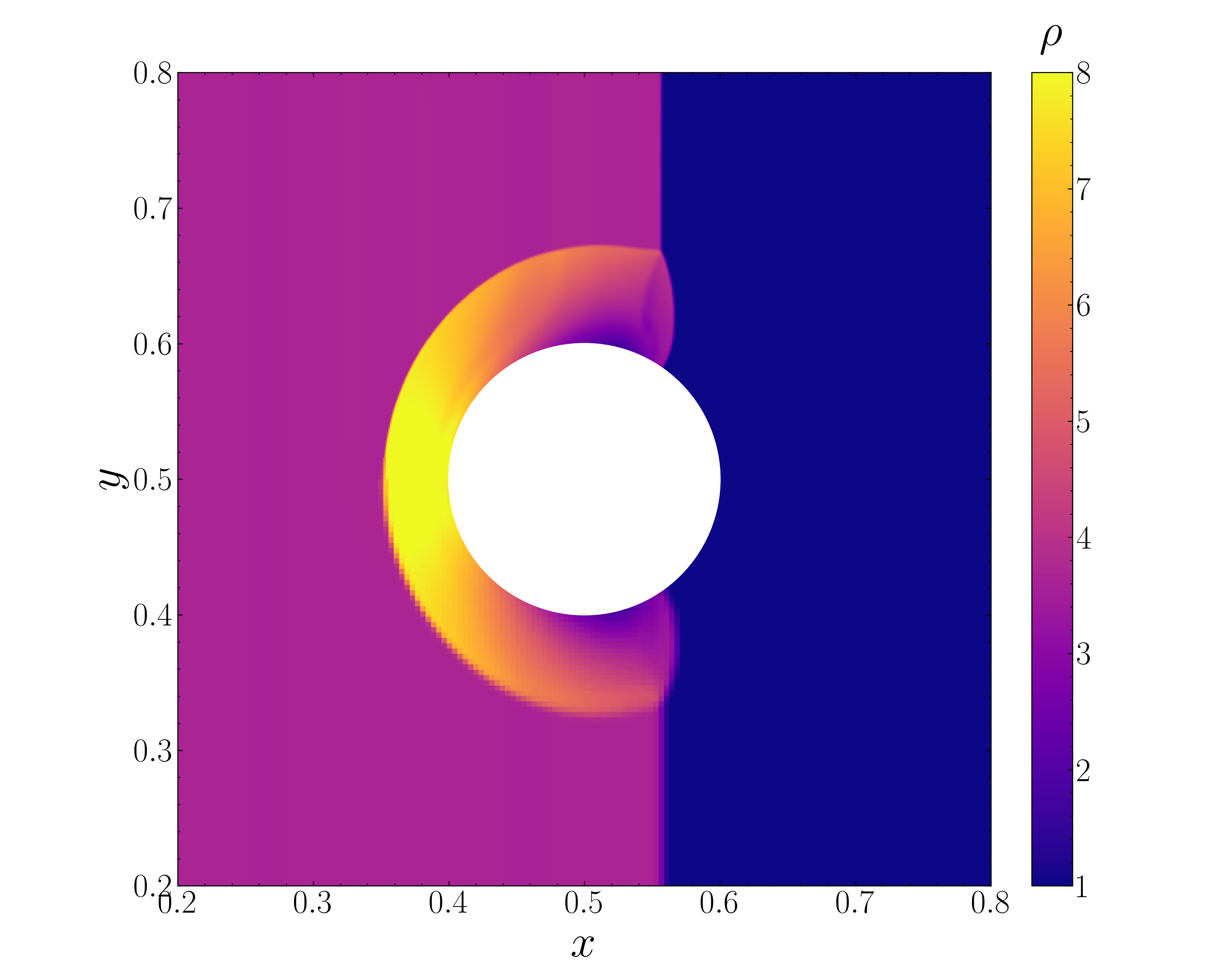}
         \caption{$t=0.068$, close-up view.}
         %\label{fig:y equals x}
     \end{subfigure}
     \hfill%
     \begin{subfigure}[t]{0.32\textwidth}
         \centering
         \includegraphics[width=\textwidth,page=3]{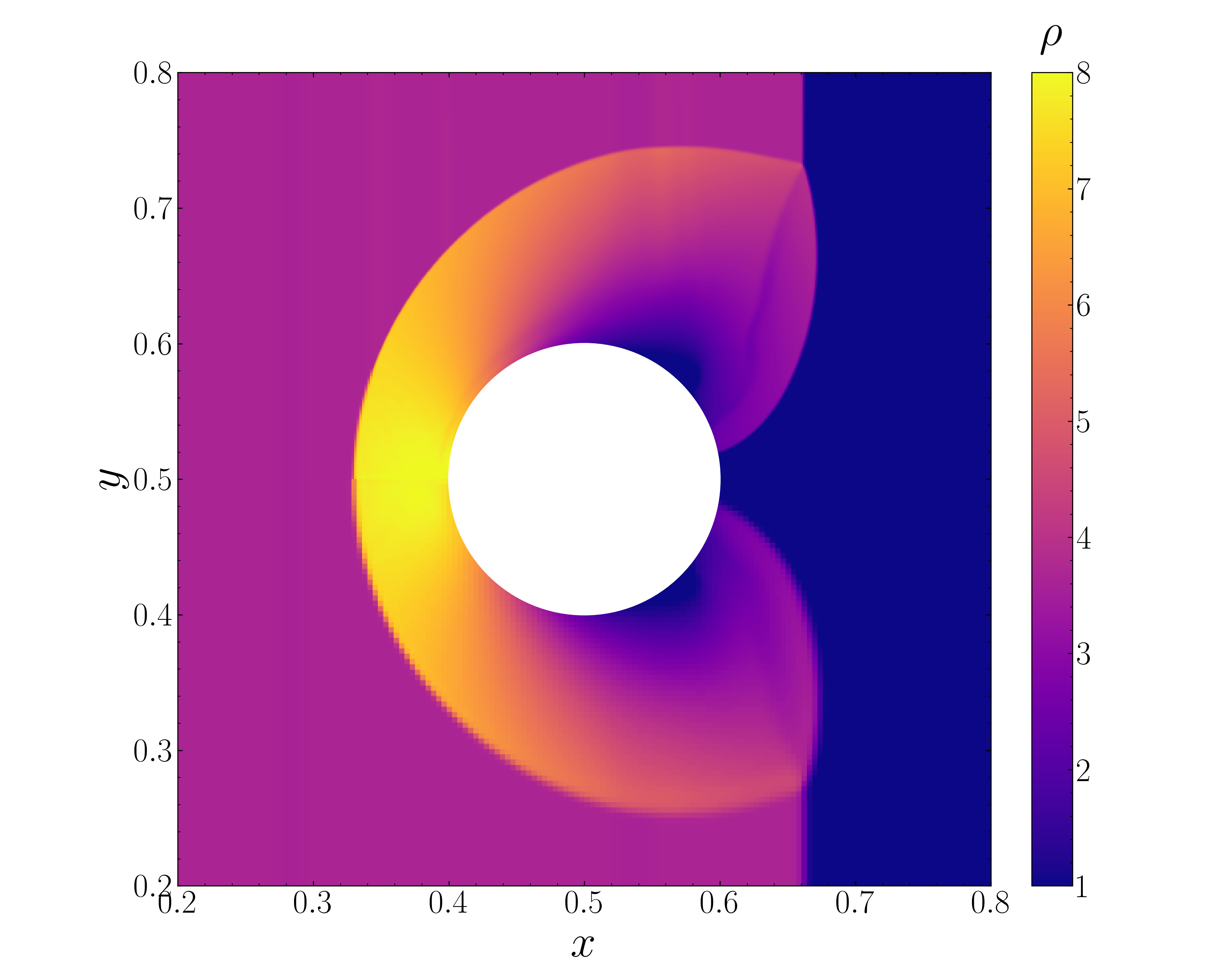}
         \caption{$t=0.088$, close-up view.}
         %\label{fig:y equals x}
     \end{subfigure}\\
     \begin{subfigure}[t]{0.32\textwidth}
         \centering
         \includegraphics[width=\textwidth,page=1]{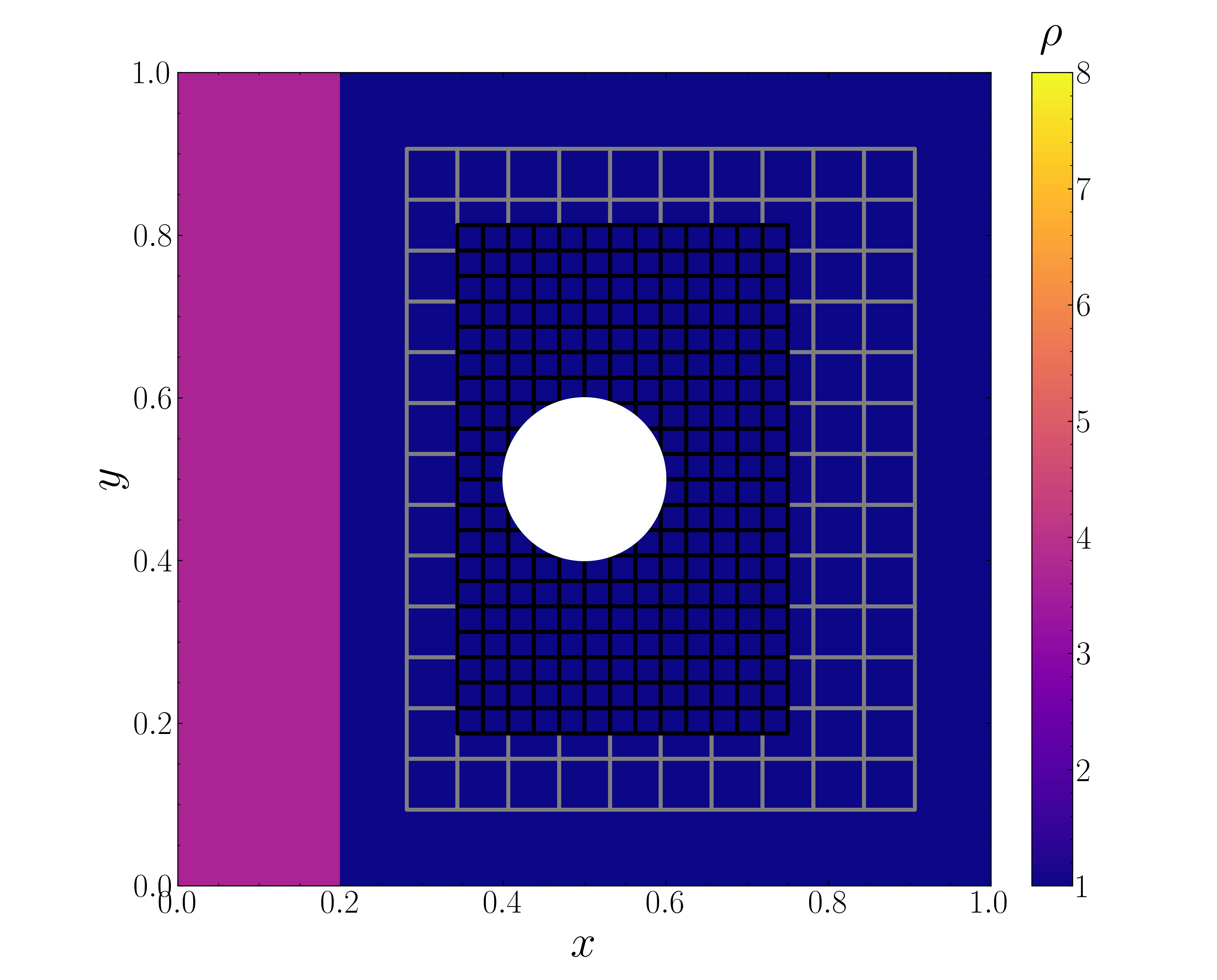}
         \caption{$t=0.0$, showing full domain and refinement levels.}
         \label{fig:cylinder-full-grids}
     \end{subfigure}
     \hfill%
     \begin{subfigure}[t]{0.32\textwidth}
         \centering
         \includegraphics[width=\textwidth,page=2]{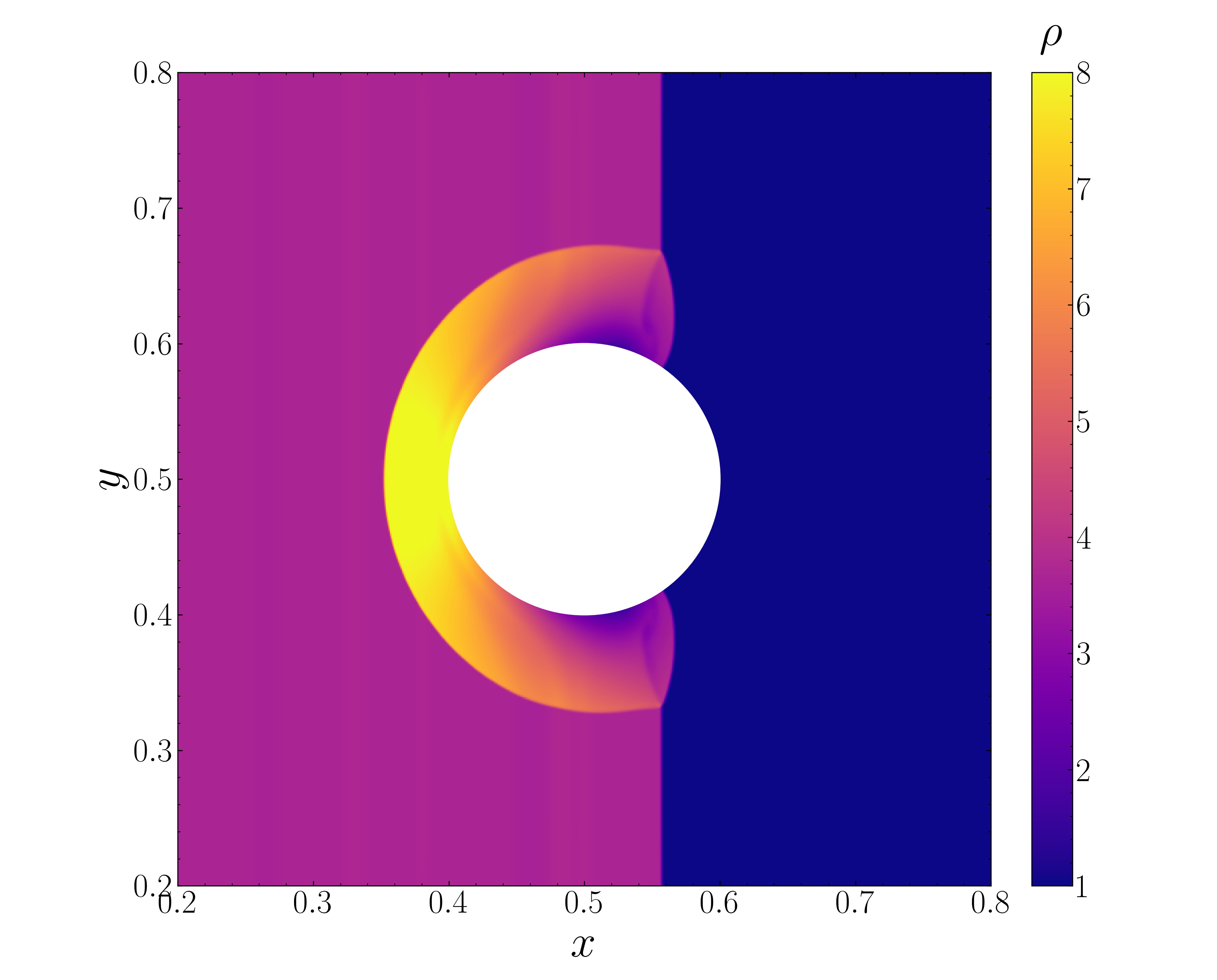}
         \caption{$t=0.068$, close-up view.}
         %\label{fig:y equals x}
     \end{subfigure}
     \hfill%
     \begin{subfigure}[t]{0.32\textwidth}
         \centering
         \includegraphics[width=\textwidth,page=3]{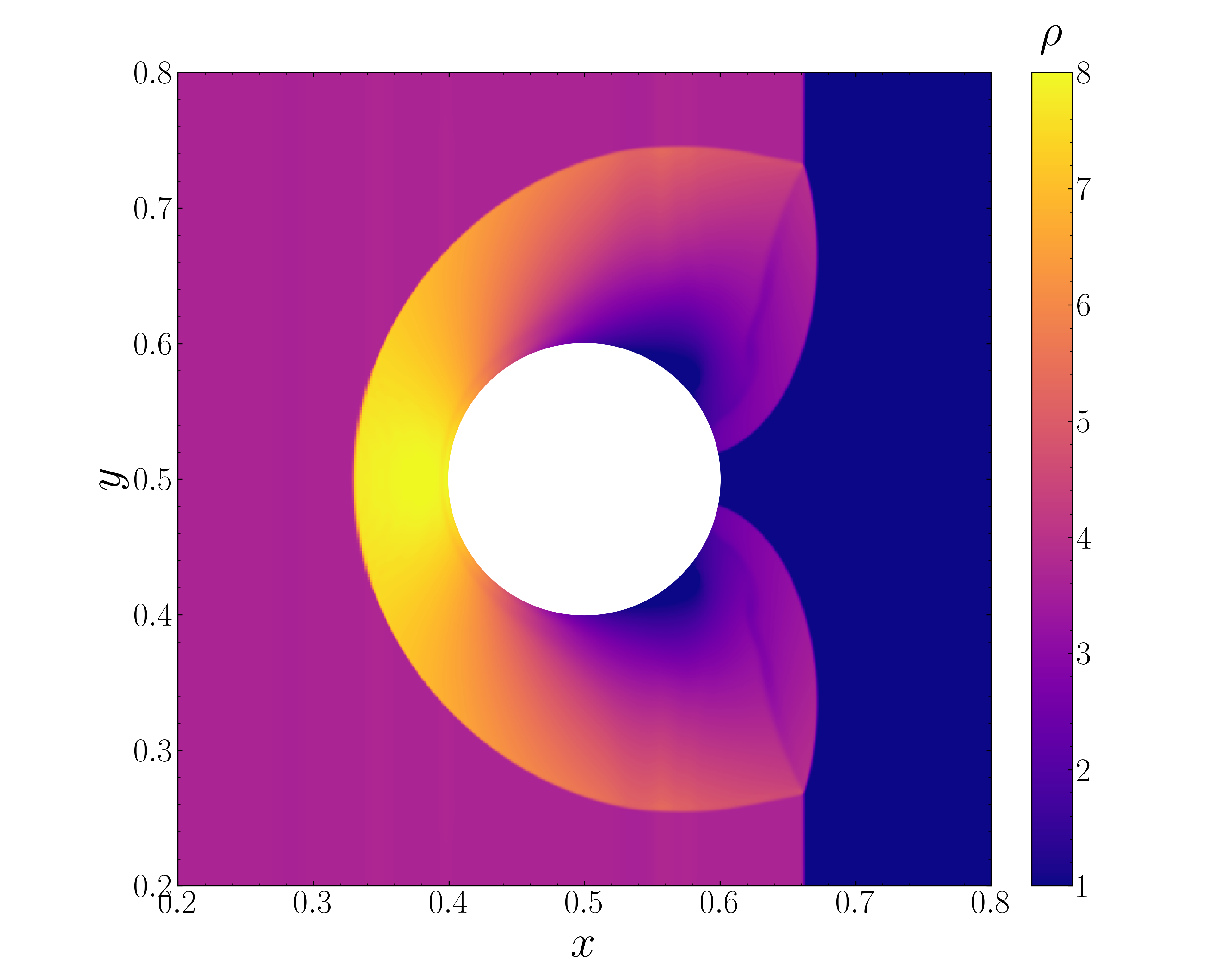}
         \caption{$t=0.088$, close-up view.}
         %\label{fig:y equals x}
     \end{subfigure}
     \caption{Profiles of density for the case of a shock interacting with a cylinder. Top row: partial refinement around the cylinder; bottom row: full refinement around the cylinder.}
     \label{fig:cylinder-amr}
\end{figure}

Here we consider two-dimensional flow past a cylinder and contrast the results with cases in which the cylinder is fully refined and with experimental data. This case has been extensively used in other works to validate numerical approaches~\cite{GULIZZI2022, ni2022immersed, jiang2022development, yang1987computation,chaudhuri2011use,zoltak1998hybrid,whitham1957new} by comparing to experimental results~\cite{bryson1961diffraction}. For this case, the domain is $[0,1] \times [0,1]$ with the center of the cylinder located at $(0.5, 0.5)$. A planar shock wave is initialized at $x=0.2$ and travels towards the cylinder at a Mach number $M=2.81$. The initial right state, upstream of the shock, is $\rho=1.0$, $u=v=0$, and $p=\nicefrac{1}{\gamma-1}$, where $\gamma=1.4$. The left state, downstream of the shock, is given by the shock jump relations~\cite{toro2013riemann}. Outflow conditions are specified at the $x$-boundaries and walls are specified at the $y$-boundaries. The base grid for the simulation is $256 \times 256$ cells. Two levels of static refinement are used. We intentionally refine only the top half of the cylinder at the first level, and slightly less at the finest, i.e., second, level. The locations of the refinement regions are shown in Figure~\ref{fig:cylinder-amr}(a). This problem differs from the earlier cases in that there are two levels of refinement, and the flow is around an object rather than within a constrained geometry. Comparisons are made to a case where the refinement region encapsulates the cylinder entirely, instead of bisecting it, as shown in Figure~\ref{fig:cylinder-amr}(d).
%\MarginPar{Is there a copyright issue with the experimental image?}

\begin{figure}
     \centering
     \begin{subfigure}[t]{0.48\textwidth}
         \centering
         \includegraphics[width=\textwidth,page=1]{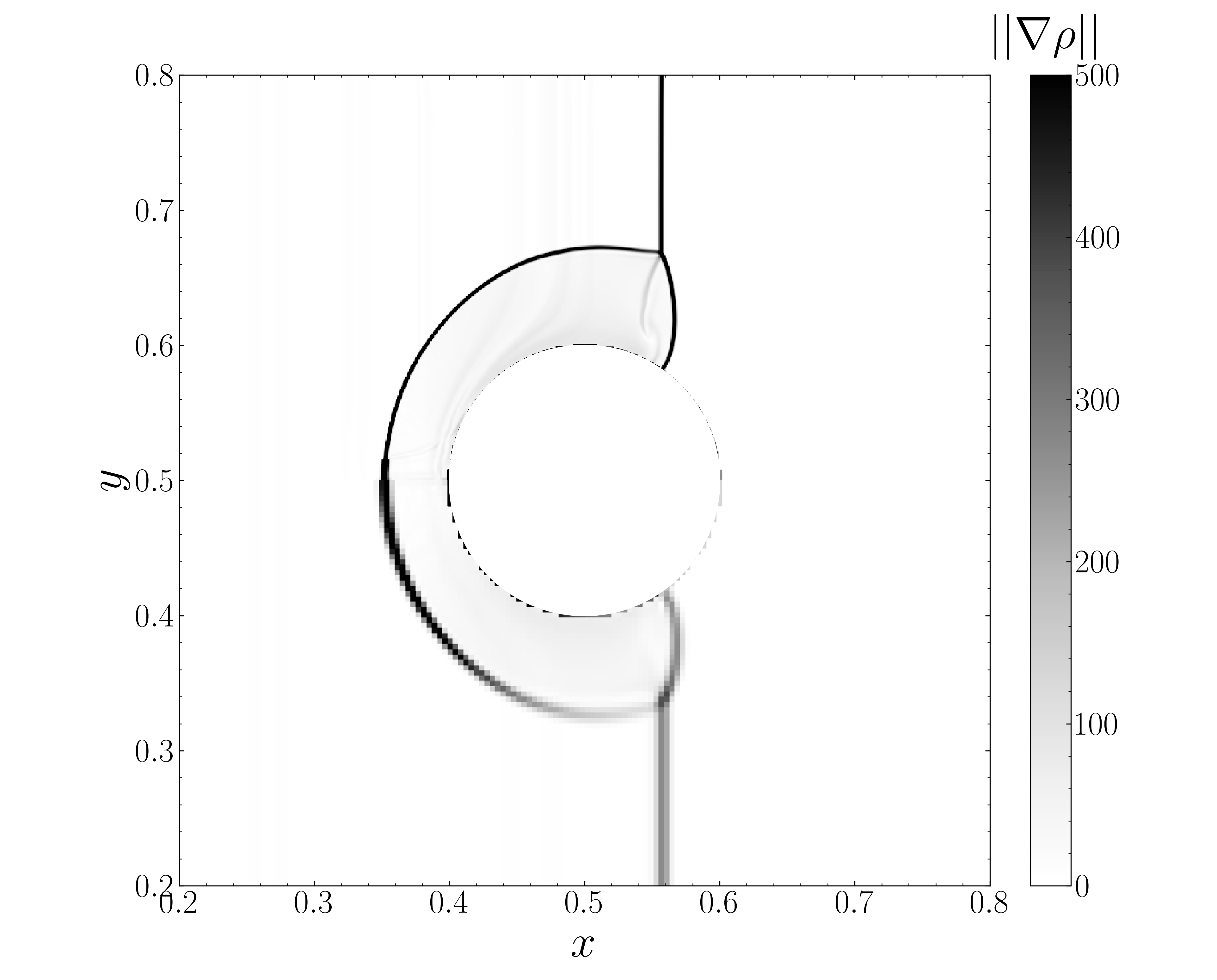}
         %\caption{$t=0.068$.}
         %\label{fig:cylinder-amr-grids}
     \end{subfigure}
     \hfill%
     \begin{subfigure}[t]{0.48\textwidth}
         \centering
         \includegraphics[width=\textwidth,page=2]{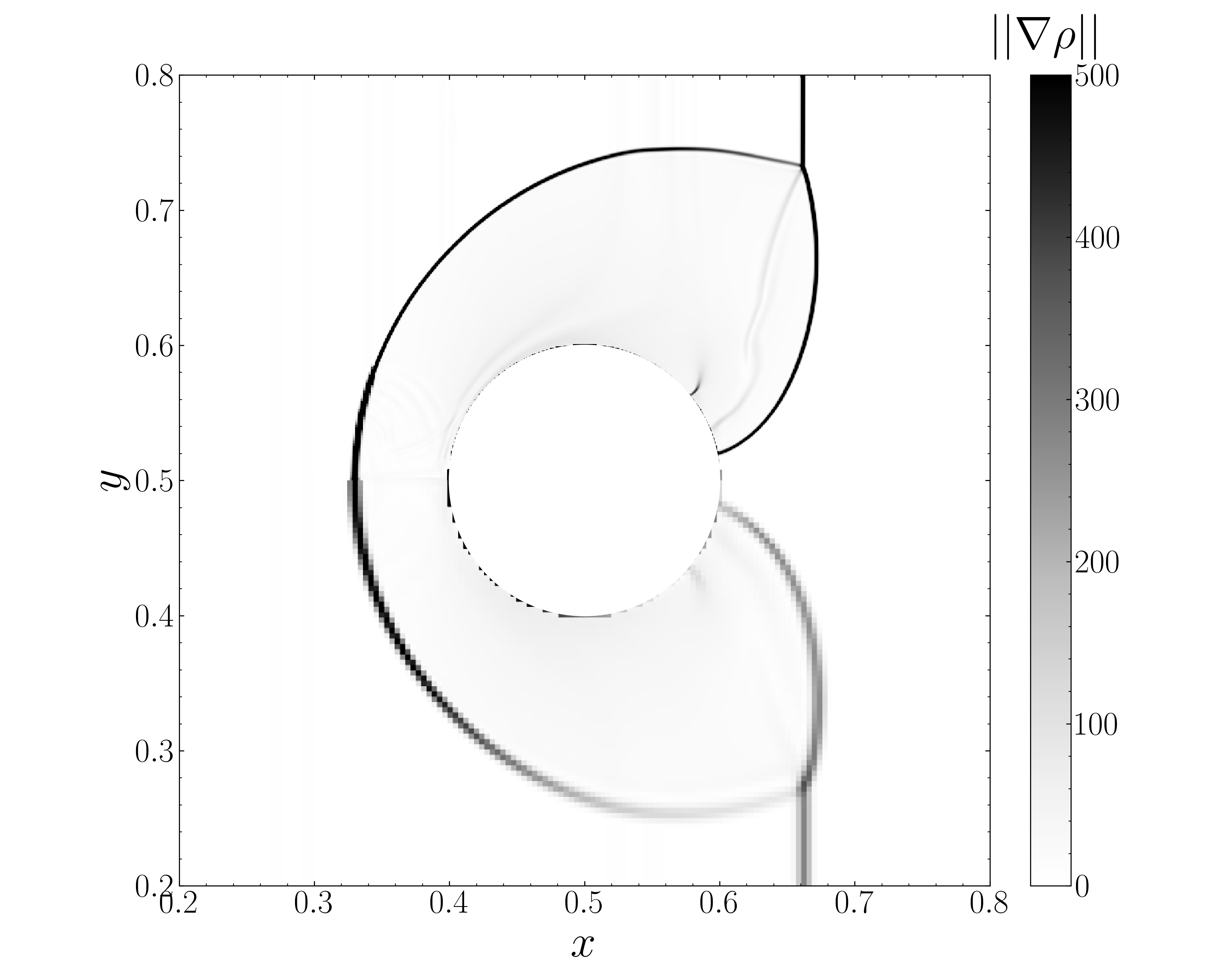}
         %\caption{$t=0.088$.}
         %\label{fig:y equals x}
     \end{subfigure}\\
     \begin{subfigure}[t]{0.48\textwidth}
         \centering
         \includegraphics[width=\textwidth,page=1]{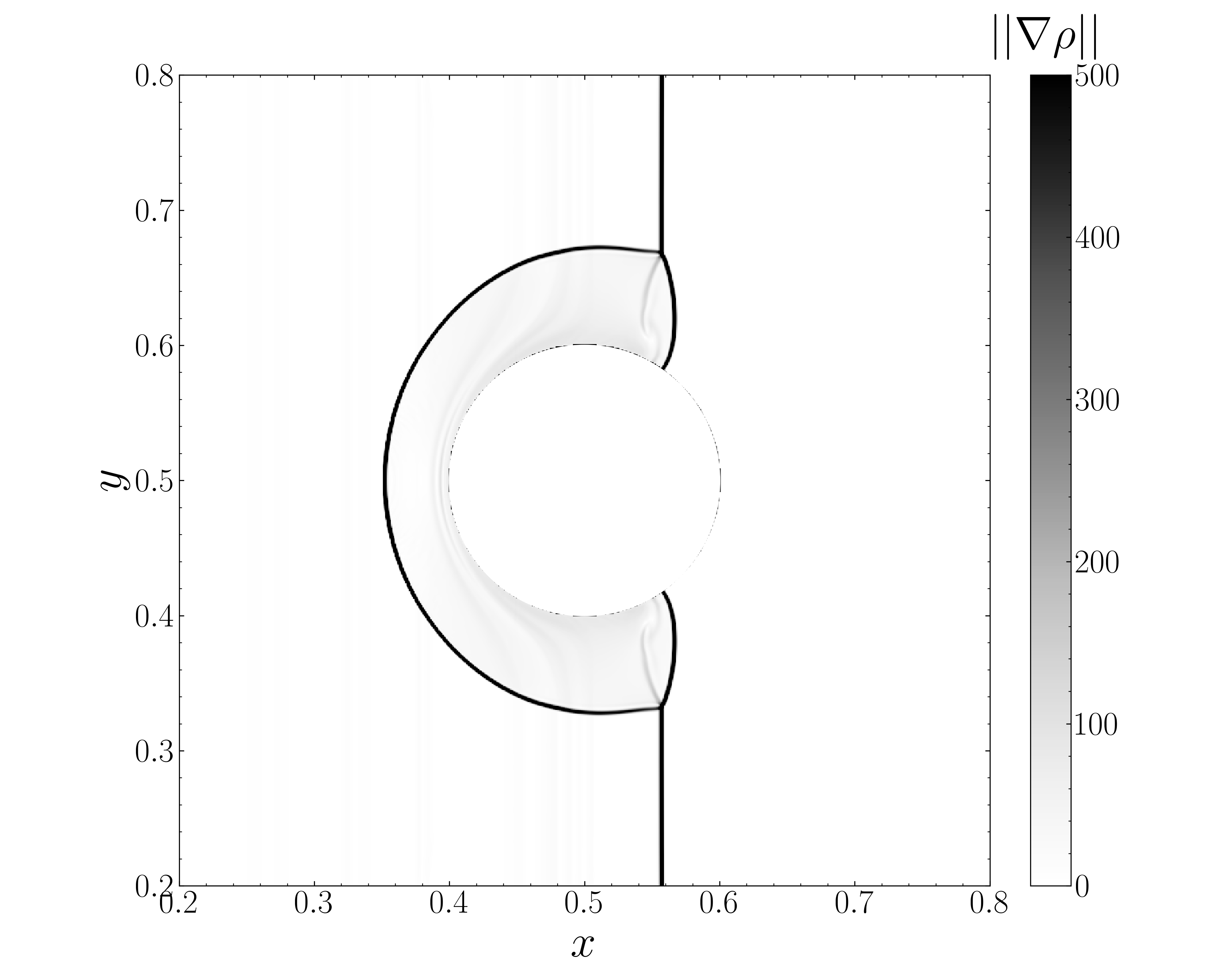}
         \caption{$t=0.068$.}
         %\label{fig:cylinder-full-grids}
     \end{subfigure}
     \hfill%
     \begin{subfigure}[t]{0.48\textwidth}
         \centering
         \includegraphics[width=\textwidth,page=2]{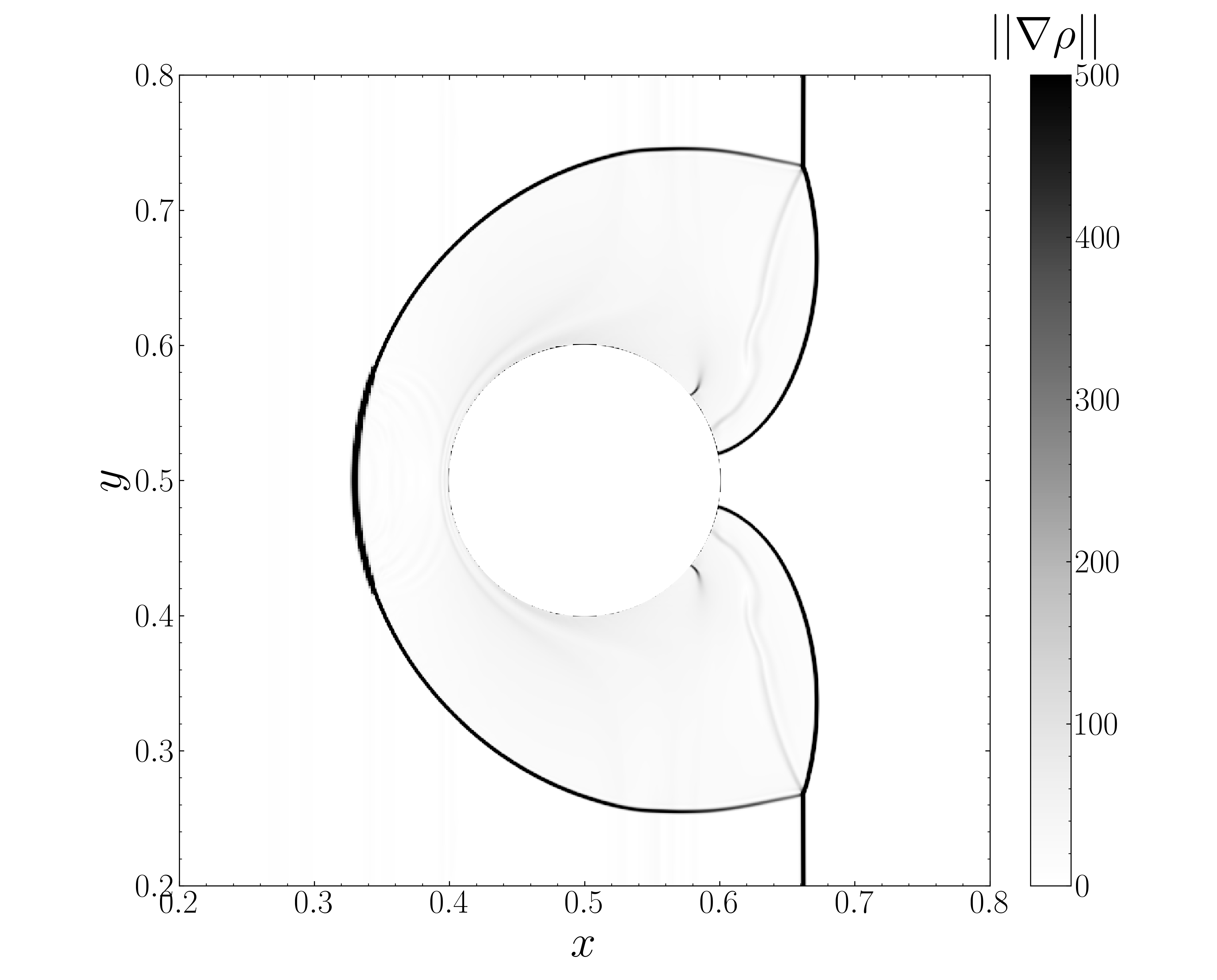}
         \caption{$t=0.088$.}
         %\label{fig:y equals x}
     \end{subfigure}
     \caption{Profiles of density gradient magnitude for the case of a shock interacting with a cylinder. Top row: partial refinement around the cylinder; bottom row: full refinement around the cylinder.}
     \label{fig:cylinder-amr-sch}
\end{figure}

\begin{figure}
\centering
  \includegraphics[width=0.5\linewidth]{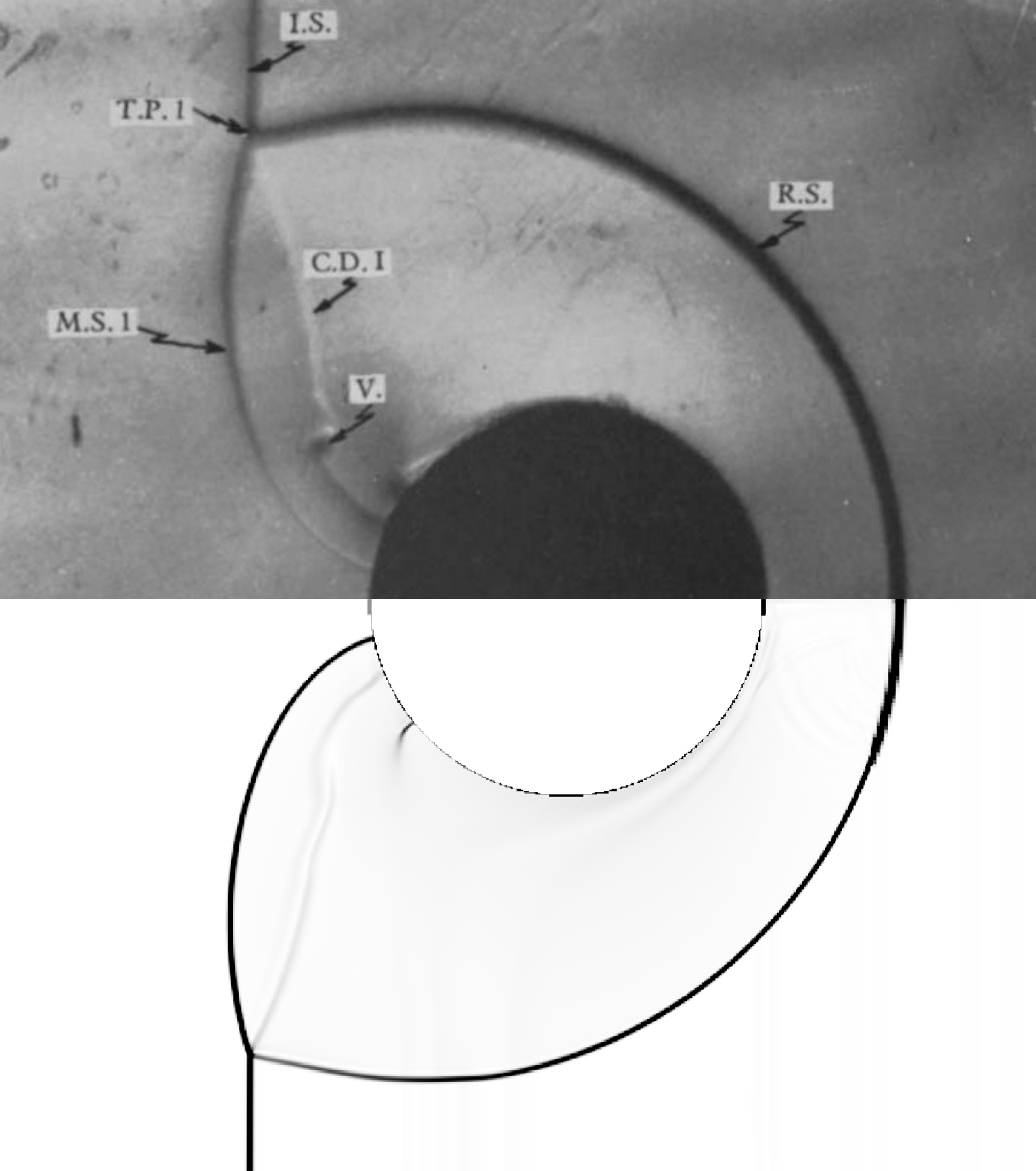}
  \caption{Comparison between the experimental Schlieren photograph~\cite{bryson1961diffraction} and the numerical solution at $t=0.088$ (case where only half the domain is refined, the section of the domain where the refinement is performed is shown here on the bottom of this figure). Note: the figure from the reference is mirrored compared to the case description.}
  \label{fig:cylinder-comp}
\end{figure}

Simulation results are shown in Figures~\ref{fig:cylinder-amr} and \ref{fig:cylinder-amr-sch}. For the case where the top half of the cylinder if refined, the solution in the refined region is similar to the case with full refinement around the cylinder. The simulation results also retain all the features, e.g., the reflected shock (R.S.), the primary Mach stem (M.S.1), the primary contact discontinuity (C.D.1), and the primary triple point (T.P.1), observed in the experimental results, Figure~\ref{fig:cylinder-comp}.

%% file: 04_results/pelec-challenge.tex
\subsection{Flow in a compression ignition engine}\label{sec:pelec-ic}
For this case, we use the new state re-redistribution scheme for three-dimensional simulation of a production-level case, a compression ignition engine. These simulations were performed with PeleC.%~\cite{PeleC_IJHPCA,PeleViz,Sitaraman2021}, a compressible reacting flow solver, part of the Pele suite of solvers. 

\subsubsection{Case description}

The case chosen for this demonstration is a large scale simulation of a simplified piston-bowl geometry, at realistic engine conditions, with a dual pulse fuel injection. In contrast with previous cases described in this work, this case involves solving the compressible Navier-Stokes equations, including viscous stresses, thermal conduction, and species diffusion~\cite{PeleC_IJHPCA}. We note that, similarly to what was shown for state redistribution in \cite{giuliani2022weighted}, the inclusion of the diffusive fluxes does not change the re-redistribution algorithm. 

\begin{figure}
\centering
  \includegraphics[width=0.5\linewidth]{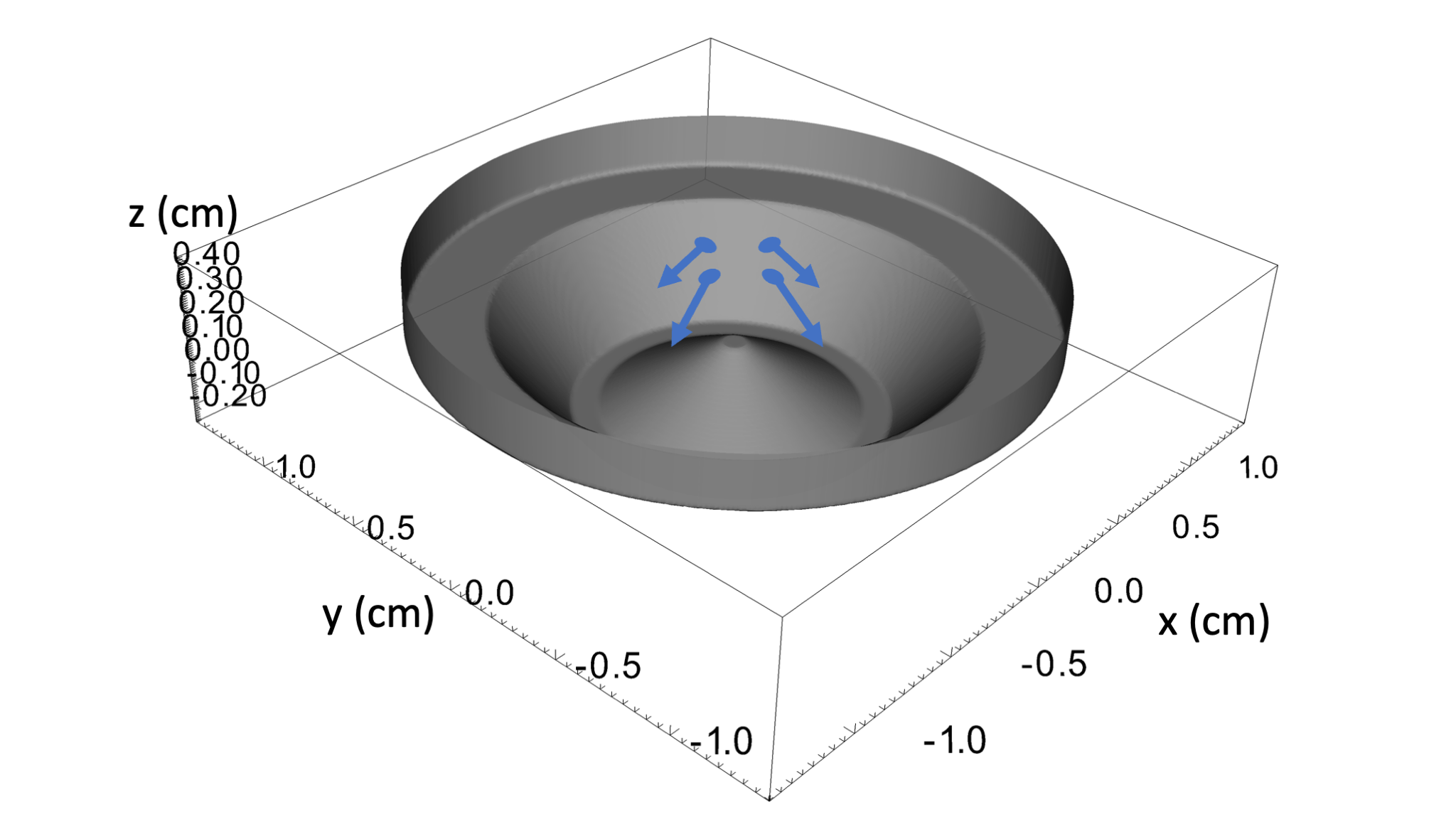}
  \caption{Computational domain, EB geometry, and jet locations for the compression ignition engine case.}
  \label{fig:pelec-geometry}
\end{figure}

This case is ideally suited to demonstrating the advantages of adaptive refinement crossing the the EB as it contains four different jets interacting with a geometry, where, for most of the simulation run time, it is not necessary to have full EB refinement. This case forms the basis of the performance benchmark used by PeleC on high performance computing systems~\cite{PeleC_IJHPCA}. The geometry contains a 24mm diameter cylinder and a piston head, leading to a 2.29mm top volume height and a 3.92mm cavity extending into the piston head. The geometry and jet positions are shown in Figure~\ref{fig:pelec-geometry}. This small volume, representative of a piston near top-dead-center, forms the computational domain. For these demonstration simulations, a relatively coarse Cartesian base grid is chosen: 128 $\times$ 128 $\times$ 32 cells. 

The initial conditions in the chamber are a mixture of methane (CH4) and vitiated air at an equivalence ratio $\phi = 0.5$. The initial pressure and temperature in the domain are 6079500.0 Pa and 900K, respectively. Turbulent velocity fluctuations are imposed on the initial velocity field to represent a low-swirl flow. 

A diesel surrogate, $n$-dodecane (NC12H26), is directly injected through four jets, angled at a $45^\circ$ angle, at the $xy$-plane at the top of the domain, 10 mm downstream of the actual diesel nozzle. A 35 transported species mechanism is used for the finite rate chemistry. The diameter of the jets is 680 $\mu$m. They are assumed to be in a fully gaseous phase (fully evaporated liquid jets). The mixture composition of the gaseous jets are 45\% of $n$-dodecane and 55\% of the initial chamber mixture. They are injected at $470$K and a normal velocity of 28 m/s. A precursor turbulent pipe flow is used to superimpose turbulent fluctuations on the jet velocities. The injection duration is 0.5 ms and the jet injection velocity is ramped down with a hyperbolic tangent function.

In previous version of this test case, as reported elsewhere~\cite{PeleC_IJHPCA}, we had the following options that avoided having coarse-fine interfaces cross the EB: (1) refine the EB entirely to the finest level present in the domain, or (2) force the EB to be refined only up to a specified level, typically less than the finest level present in the domain and do not allow for refined sections at the EB. Option 1 leads to a computationally unfeasible number of cells, even on the largest supercomputers. For example, if one is resolving the ignition kernels and therefore requiring six levels of AMR on a base grid of 512 $\times$ 512 $\times$ 128 cells, even without including the required cells on the EB, there would 68 billion cells in the interior of the domain. Including the cells on the EB would lead to $\mathcal{O}(10)$ times more cells. Option 2 is the one currently pursued for production simulations of this case but leads to underresolved physics at the EB, particularly when the jets start interacting with the geometry. A third option is enabled through this work: refining the EB using different resolutions to respond to the physical phenomena at the EB.
 
\subsubsection{Results}
We assess the three refinement options using this case by comparing flow features and performance metrics. For the cases presented here, in order to maintain a reasonable computational cost while still demonstrating the capabilities of the proposed scheme, the simulations use one level of adaptive mesh refinement and chemical reactions are not evaluated (though all the relevant species are transported). Simulations were conducted to $t=0.9$ms. They were performed on 4 nodes (36 Intel Xeon-Gold Skylake-6154 processor cores per node) on NREL's Eagle supercomputer~\cite{eagle-ref}.

\begin{figure}
     \centering
     \begin{subfigure}[b]{0.32\textwidth}
         \centering
         \includegraphics[width=\textwidth,page=1, clip=true, trim=0.4cm 4cm 0.4cm 0.5cm]{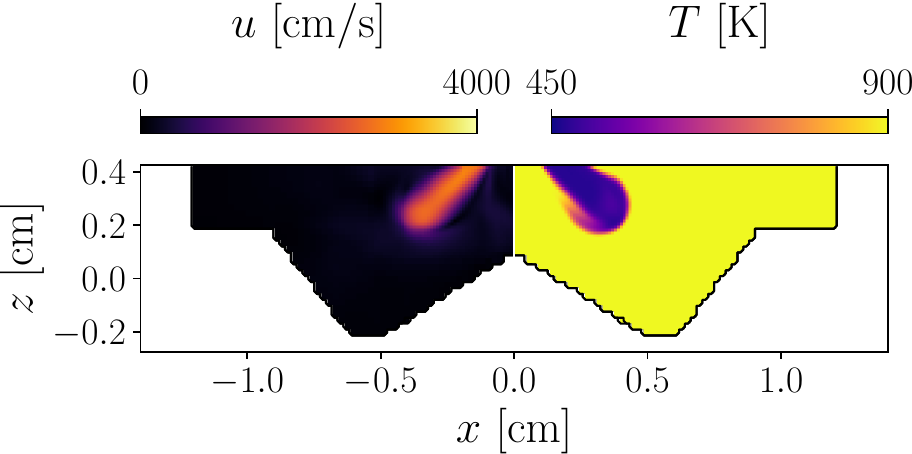}
         %\caption{$t=0.2$ms}
         %\label{fig:y equals x}
     \end{subfigure}
     \hfill%
     \begin{subfigure}[b]{0.32\textwidth}
         \centering
         \includegraphics[width=\textwidth,page=1, clip=true, trim=0.4cm 4cm 0.4cm 3cm]{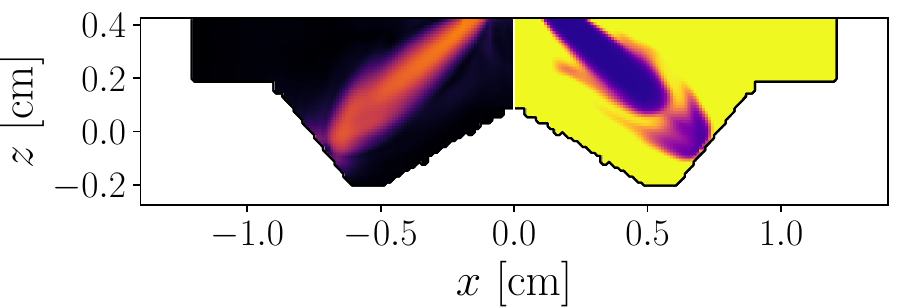}
         %\caption{$t=0.45$ms}
         %\label{fig:three sin x}
     \end{subfigure}
     \hfill%
     \begin{subfigure}[b]{0.32\textwidth}
         \centering
         \includegraphics[width=\textwidth,page=1, clip=true, trim=0.4cm 4cm 0.4cm 3cm]{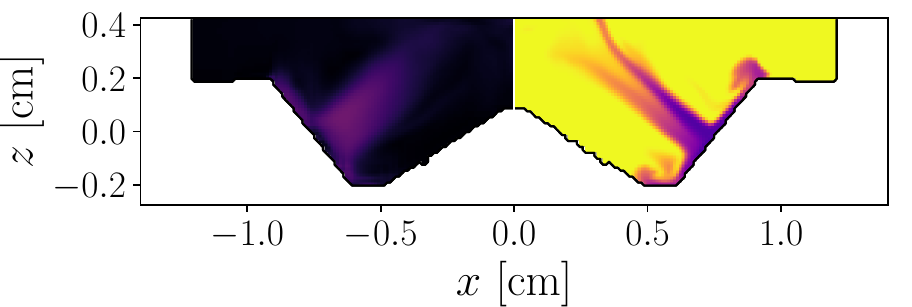}
         %\caption{$t=0.9$ms}
         %\label{fig:five over x}
     \end{subfigure}\\
     \begin{subfigure}[b]{0.32\textwidth}
         \centering
         \includegraphics[width=\textwidth,page=1, clip=true, trim=0.4cm 4cm 0.4cm 3cm]{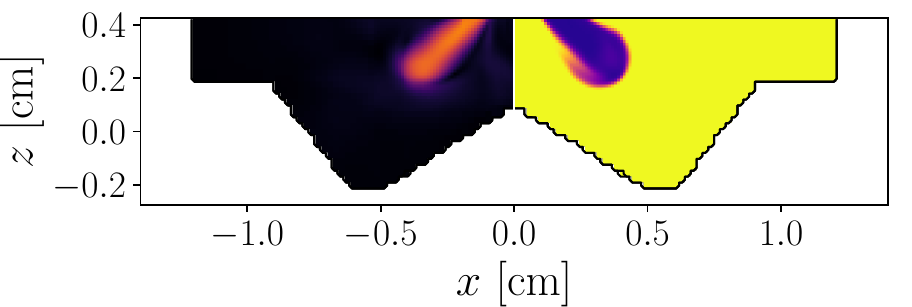}
         %\caption{$t=0.2$ms}
         %\label{fig:y equals x}
     \end{subfigure}
     \hfill%
     \begin{subfigure}[b]{0.32\textwidth}
         \centering
         \includegraphics[width=\textwidth,page=1, clip=true, trim=0.4cm 4cm 0.4cm 3cm]{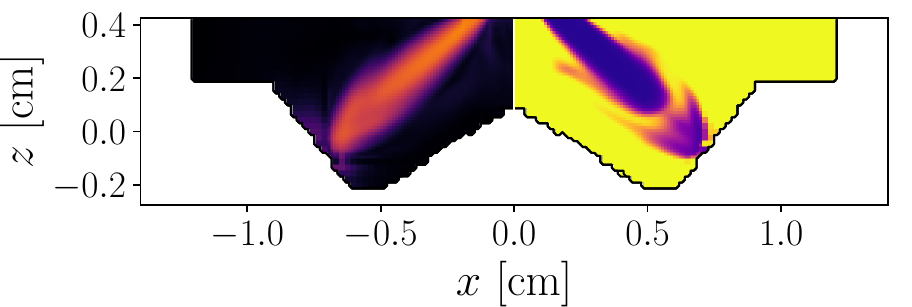}
         %\caption{$t=0.45$ms}
         %\label{fig:three sin x}
     \end{subfigure}
     \hfill%
     \begin{subfigure}[b]{0.32\textwidth}
         \centering
         \includegraphics[width=\textwidth,page=1, clip=true, trim=0.4cm 4cm 0.4cm 3cm]{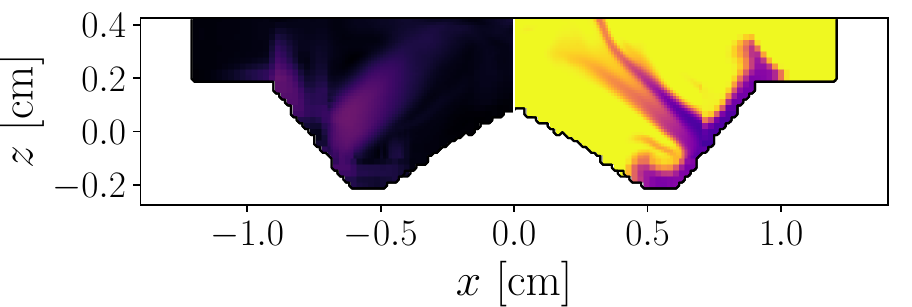}
         %\caption{$t=0.9$ms}
         %\label{fig:five over x}
     \end{subfigure}\\
     \begin{subfigure}[b]{0.32\textwidth}
         \centering
         \includegraphics[width=\textwidth,page=1, clip=true, trim=0.4cm 4cm 0.4cm 3cm]{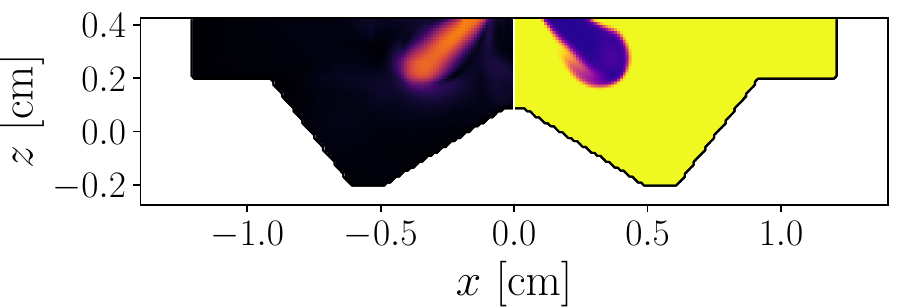}
         \caption{$t=0.2$ms}
         %\label{fig:y equals x}
     \end{subfigure}
     \hfill%
     \begin{subfigure}[b]{0.32\textwidth}
         \centering
         \includegraphics[width=\textwidth,page=1, clip=true, trim=0.4cm 4cm 0.4cm 3cm]{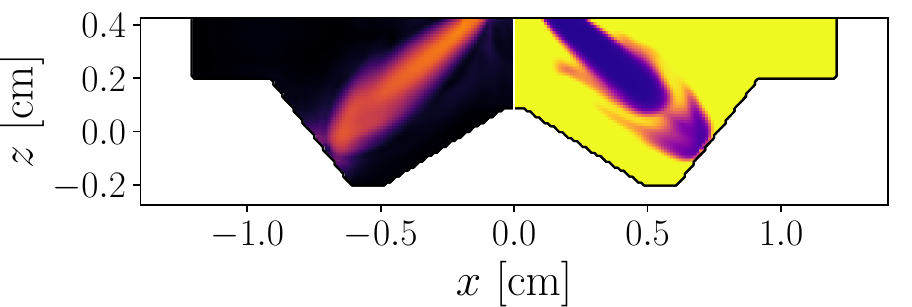}
         \caption{$t=0.45$ms}
         %\label{fig:three sin x}
     \end{subfigure}
     \hfill%
     \begin{subfigure}[b]{0.32\textwidth}
         \centering
         \includegraphics[width=\textwidth,page=1, clip=true, trim=0.4cm 4cm 0.4cm 3cm]{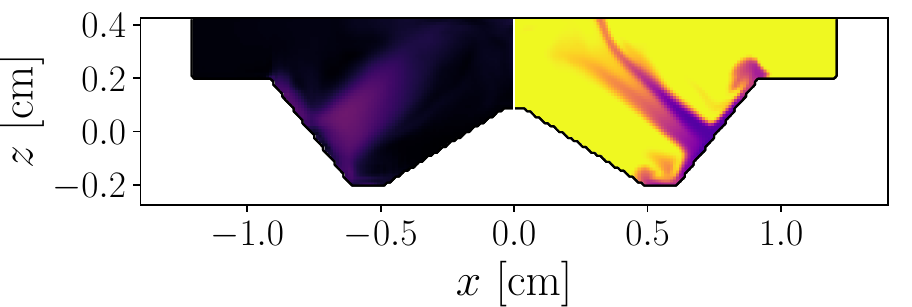}
         \caption{$t=0.9$ms}
         %\label{fig:five over x}
     \end{subfigure}
     \caption{Pseudocolor vizualization of the velocity magnitude and temperature fields. Top row: using state re-redistribution, middle row: refining the EB to one level coarser than the finest level, bottom row: refining the EB to the maximum number of AMR levels.}
     \label{fig:piston-bowl-flow}
\end{figure}

\begin{figure}
     \centering
     \begin{subfigure}[b]{0.32\textwidth}
         \centering
         \includegraphics[width=\textwidth,page=1, clip=true, trim=0.4cm 4cm 0.4cm 0.5cm]{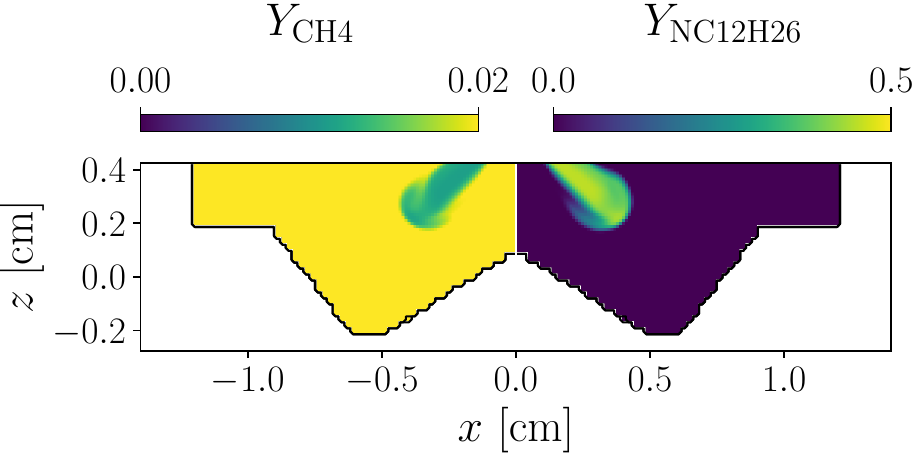}
         %\caption{$t=0.2$ms}
         %\label{fig:y equals x}
     \end{subfigure}
     \hfill%
     \begin{subfigure}[b]{0.32\textwidth}
         \centering
         \includegraphics[width=\textwidth,page=1, clip=true, trim=0.4cm 4cm 0.4cm 3cm]{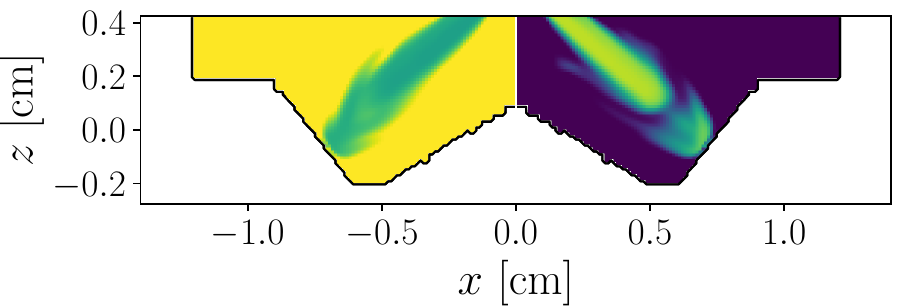}
         %\caption{$t=0.45$ms}
         %\label{fig:three sin x}
     \end{subfigure}
     \hfill%
     \begin{subfigure}[b]{0.32\textwidth}
         \centering
         \includegraphics[width=\textwidth,page=1, clip=true, trim=0.4cm 4cm 0.4cm 3cm]{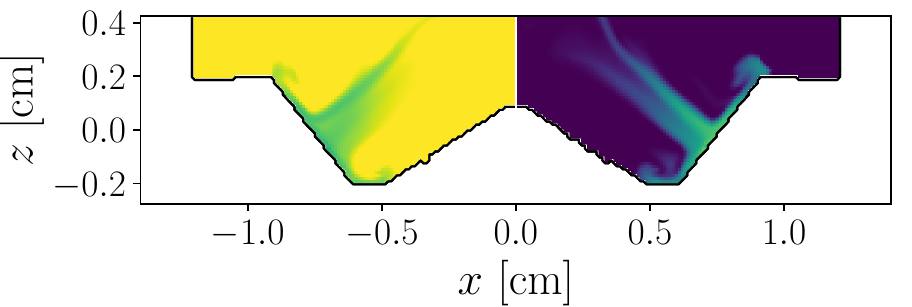}
         %\caption{$t=0.9$ms}
         %\label{fig:five over x}
     \end{subfigure}\\
     \begin{subfigure}[b]{0.32\textwidth}
         \centering
         \includegraphics[width=\textwidth,page=1, clip=true, trim=0.4cm 4cm 0.4cm 3cm]{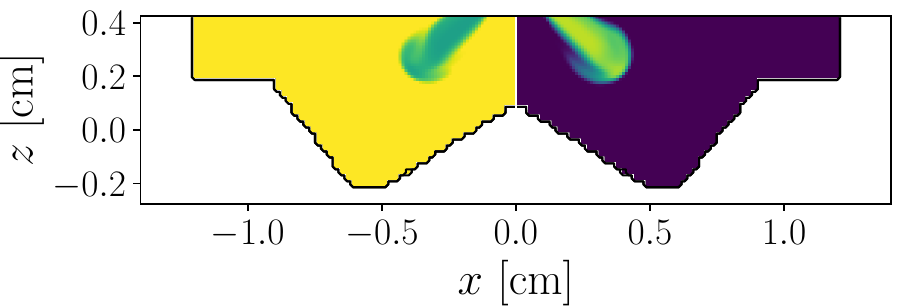}
         %\caption{$t=0.2$ms}
         %\label{fig:y equals x}
     \end{subfigure}
     \hfill%
     \begin{subfigure}[b]{0.32\textwidth}
         \centering
         \includegraphics[width=\textwidth,page=1, clip=true, trim=0.4cm 4cm 0.4cm 3cm]{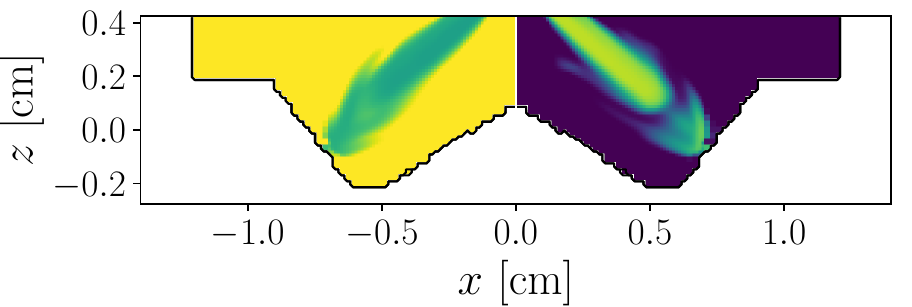}
         %\caption{$t=0.45$ms}
         %\label{fig:three sin x}
     \end{subfigure}
     \hfill%
     \begin{subfigure}[b]{0.32\textwidth}
         \centering
         \includegraphics[width=\textwidth,page=1, clip=true, trim=0.4cm 4cm 0.4cm 3cm]{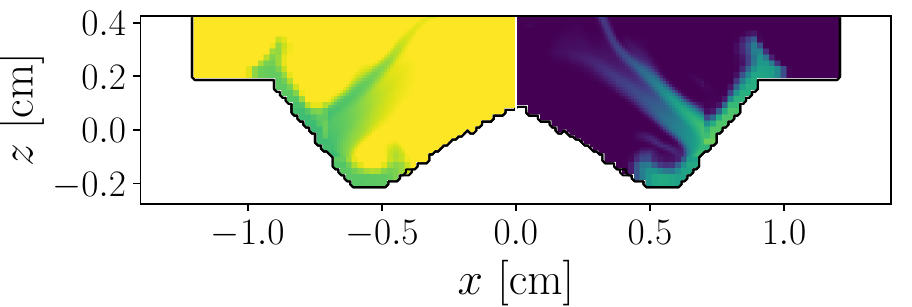}
         %\caption{$t=0.9$ms}
         %\label{fig:five over x}
     \end{subfigure}\\
     \begin{subfigure}[b]{0.32\textwidth}
         \centering
         \includegraphics[width=\textwidth,page=1, clip=true, trim=0.4cm 4cm 0.4cm 3cm]{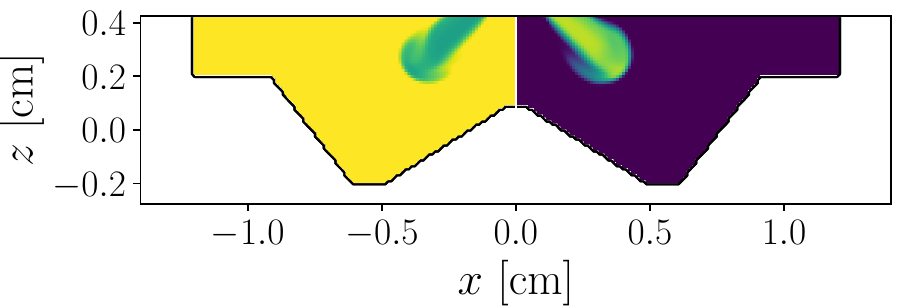}
         \caption{$t=0.2$ms}
         %\label{fig:y equals x}
     \end{subfigure}
     \hfill%
     \begin{subfigure}[b]{0.32\textwidth}
         \centering
         \includegraphics[width=\textwidth,page=1, clip=true, trim=0.4cm 4cm 0.4cm 3cm]{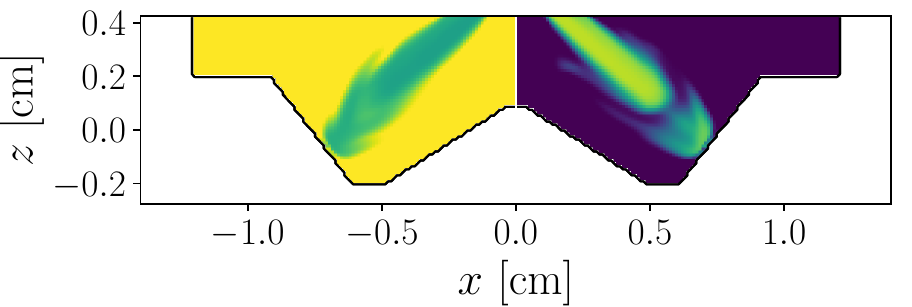}
         \caption{$t=0.45$ms}
         %\label{fig:three sin x}
     \end{subfigure}
     \hfill%
     \begin{subfigure}[b]{0.32\textwidth}
         \centering
         \includegraphics[width=\textwidth,page=1, clip=true, trim=0.4cm 4cm 0.4cm 3cm]{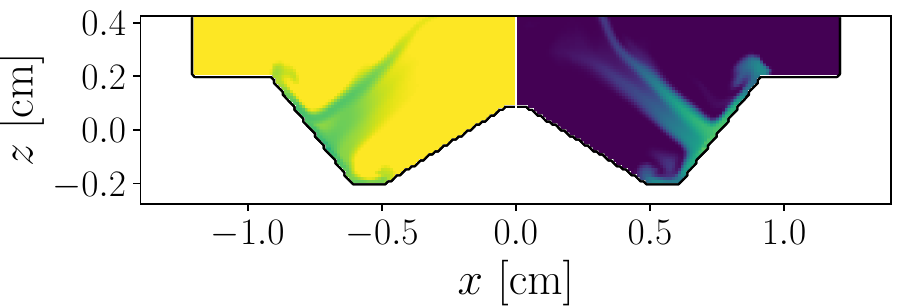}
         \caption{$t=0.9$ms}
         %\label{fig:five over x}
     \end{subfigure}
     \caption{Pseudocolor vizualization of the CH4 and NC12H26 ($n$-dodecane) mass fraction fields. Top row: using state re-redistribution, middle row: refining the EB to one level coarser than the finest level, bottom row: refining the EB to the maximum number of AMR levels.}
     \label{fig:piston-bowl-flow-ms}
\end{figure}

Flow field comparisons for the three options for the magnitude of velocity and temperature are shown in Figure~\ref{fig:piston-bowl-flow}  and  for the mass fractions of methane and $n$-dodecane in Figure~\ref{fig:piston-bowl-flow-ms}. From these flow field visualizations, it is clear that the option that forces the EB to be refined at a level coarser than the finest available level (middle row in the figures) shows some notable differences compared to the the other two cases (full refinement at the EB and adaptive refinement at the EB). For this option, as the jet interacts with the wall, starting around $t=0.45$ms, the jet features appear more diffuse. The jet spreads further along the EB (over the inflection point and into the bottom of the piston bowl). The case using adaptive refinement at the EB (top row) presents flow feature identical to the fully refined case (bottom row). Table~\ref{tab:pelec-ic} shows the performance metrics of these cases at two different times ($t=0.45$ms and $t=0.9$ms): the number of cells on the finest level, the percentage of the domain that is covered by the finest levels, and the time per time step. At $t=0.45$, just as the jets reach the walls, the time per time step for the case with adaptive refinement at the EB and the case with a coarse EB refinement are similar and the percent of the domain that is covered by the finest level is approximately the same. The case with full refinement of the EB is approximately four times more computationally expensive, while achieving the same flow field results. At $t=0.9$ms, after the jets have significantly interacted with the walls, the case with adaptive refinement at the EB has increased in computational cost by $50\%$ but still only has about half the cells on the finest level compared to the case with full EB refinement and is a 3.2 times faster computationally.

\begin{table}[]
\centering
\begin{tabular}{lccc}
\toprule
Max level at EB  & $\#$ cells on finest level & $\%$ finest level & time per step (s)\\
\midrule
at $t=0.45$ms: & & & \\
$\quad$1 (adaptive) & 491520 & 11.7 & 10.1\\
$\quad$0 & 327680 & 7.8 & 10.4 \\
$\quad$1 (full) & 1982464 & 47.3 & 44.0\\
\midrule
at $t=0.9$ms: & & & \\
$\quad$1 (adaptive) & 1150976 & 27.4 & 14.4\\
$\quad$0 & 319488 & 7.6 & 10.0 \\
$\quad$1 (full) & 2068480 & 49.0 & 45.3 \\
\bottomrule
\end{tabular}
\caption{Performance comparison at $t=0.45$ms and $t=0.9$ms. Top row: using state re-redistribution, middle row: refining the EB to one level coarser than the finest level, bottom row: refining the EB to the maximum number of AMR levels.}
\label{tab:pelec-ic}
\end{table}

%% file: conclusions.tex
Weighted state redistribution is a recently developed approach for addressing the ``small cell" problem that arises in embedded boundary discretizations.  When coupled with adaptive mesh refinement, additional synchronization steps are needed when the interfaces between different levels of refinement are close to the embedded boundary to ensure that the overall method is both stable and conservative. Consequently, previous implementations of (weighted) state redistribution with AMR have required the entire embedded boundary to be resolved at the same level, typically the finest level in the problem.

In this paper, we have developed a synchronization procedure for WSRD that removes this restriction.  
\commentout{The new approach consists of several basic steps.  First we identify, in extensive form, what the weighted state redistribution algorithm moves across coarse/fine boundaries during the independent integration of the solution at different levels.  That information is then used to correct solution in cells near the coarse/fine boundary, making the overall update conservative.  When correcting the solution in cut cells, an additional redistribution of the synchronization correction is needed to ensure stability.}
The resulting methodology was tested in two and three space dimensions for systems of hyperbolic conservation laws.  We demonstrated that the algorithm is conservative and stable in fairly simple geometries using both a single-step piecewise-linear Godunov method and a predictor/corrector method of lines temporal integration scheme. We also presented validation cases for a Sod shock tube in slanted channel and shock reflection from a sphere.
Finally, we demonstrated the methodology on a more complex model for a compression ignition engine that includes multiple species and diffusive transport.  In this more complex example we illustrate the potential computational savings resulting from not requiring the entire embedded boundary to be resolved at the finest level.

The methodology here can be extended to more general compressible flows.  In particular, following the overall approach outlined in \cite{giuliani2022weighted}, the approach could be extended to include reactions.  This approach can also be extended to incompressible flows and more general low Mach number models. We also note that there are a number of potential variations to the basic WSRD approach that could be considered.

%% file: 02_algorithms/FRD_intro.tex
Flux redistribution is one approach to dealing with the small cell problem.  
%The basic idea is to compute both a conservative but unstable update to the solution, $\delta U^\cons_\ivec$, and a stable but non-conservative update, $\delta U^\ncons_\ivec$.
The basic idea is to use the conservative but unstable update to the solution, $\delta U^\cons_\ivec$, from Eq. (\ref{eq:cons_up}) to compute a stable but non-conservative update, $\delta U^\ncons_\ivec$.
A portion of these two updates is used to update the cell and the remainder is ``redistributed" to neighboring cells so that the resulting algorithm is conservative and stable.
We first define $\delta U^\cons_\ivec$ using Eq. (\ref{eq:cons_up}) for each valid cell in the domain.
\commentout{
\[(\nabla \cdot {F})^c_i = \dfrac{1}{\mathcal{V}_i} \sum_{f=1}^{N_f} ({F}_f\cdot{n}_f) A_f \]

Here $N_f$ is the number of faces of cell $i$, $\vec{n}_f$ and $A_f$
are the unit normal and area of the $f$-th face respectively,
and %\mathcal{V}_i$ is the volume of cell $i$ given by

\[\mathcal{V}_i = (\Delta x \Delta y \Delta z)\cdot \mathcal{K}_i \]

where $\mathcal{K}_i$ is the volume fraction of cell$i$.

Now, a conservative update can be written as

\[\frac{ \rho^{n+1} \phi^{n+1} - \rho^{n} \phi^{n} }{\Delta t} = - \nabla \cdot{F}^c \]
}

Then, as done in Trebotich and Graves \cite{trebotich:2015}, for each cell cut by the EB geometry, we compute the non-conservative update, 
\[
\delta U^\ncons_\ivec = \dfrac{\sum\limits_{ii,jj,kk \in \NFRD(i,j,k) } \vfrac_{ii,jj,kk}
\; \delta U^\cons_{ii,jj,kk} }{\sum\limits_{ii,jj,kk\in \NFRD(i,j,k) } \vfrac_{ii,jj,kk} }\]
where $\NFRD(i,j,k)$ is the set of cells $(ii,jj,kk)$ such that none of the indices differs by more than one from $(i,j,k)$ and where $(ii,jj,kk)$ can be reached from $(i,j,k)$ by a monotone path of valid cells in index space; i.e., a path in which none of the indices are both incremented and decremented. 
%The monotone path criterion prevents including points that go around corners.
(The approach taken here is not unique.  See Pember {\it et al.} \cite{pember1995adaptive} and Colella {\it et al.} \cite{colella2006cartesian} for an alternative approach to computing the nonconservative updates.)

For each cell cut by the EB geometry, we compute the initial convective update $\delta U^{EB,*}_\ivec$ using
\begin{equation}
    \delta U^{EB,*}_\ivec = \vfrac_\ivec \delta U^{\cons}_\ivec +(1-\vfrac_\ivec) \delta U^{\ncons}_\ivec
    \label{eq:init_update}
\end{equation}
Eq. (\ref{eq:init_update}) provides a stable but non-conservative update for cell-$\ivec$.
A direct computation shows that for cut cells, the extensive difference between $\delta U^{EB}$ and the conservative update is given by
\begin{align}
    \delM_\ivec &=  \vol_\ivec (1- \vfrac_\ivec)[\delta U^{\cons}_\ivec -\delta U^{\ncons}_\ivec] \\
    &= \Delta x \Delta y \Delta z \; \vfrac_\ivec (1- \vfrac_\ivec)[\delta U^{\cons}_\ivec -\delta U^{\ncons}_\ivec] \nonumber
\end{align}
For the method to be conservative, we need to redistribute $\delM$ to the neighboring cells.
This is done by modifying Eq. (\ref{eq:init_update}) using
\begin{equation}
    \delta U^{EB}_{ii,jj,kk} =    \delta U^{EB,*}_{ii,jj,kk} + \sum_{i,j,k | ii,jj,kk\in \NFRD(i,j,k) } w^\ivec_{ii,jj,kk}
    \delM_\ivec\ %, \qquad \forall \; ii,jj,kk\in N(i,j,k)
    \label{eq:redist_M}
\end{equation}
where the weights, $w^\ivec_{ii,jj,kk}$, are given  
\begin{equation}
    w^\ivec_{ii,jj,kk} = \dfrac{1}{\sum\limits_{ii,jj,kk\in \NFRD(\ivec)}  \vol_{ii,jj,kk}}
    = \frac{1}{\Delta x \, \Delta y \, \Delta z} \dfrac{1}{\sum\limits_{ii,jj,kk\in \NFRD(\ivec)}  \vfrac_{ii,jj,kk}}
\end{equation}
We note that a number of alternative weightings for the redistribution are possible.  For example, in gas dynamics a density weighted redistribution given by 
\[
    w^\ivec_{ii,jj,kk} = \dfrac{\rho_{ii,jj,kk}}{\sum\limits_{ii,jj,kk\in \NFRD(\ivec)}  \rho_{ii,jj,kk} \vol_{ii,jj,kk}}
\]
is often used.
\commentout{
Note that $\nabla \cdot{F}_i^{EB}$ gives an update for $`\rho \phi$ ; i.e.,

\[\frac{(\rho \phi_i)^{n+1} - (\rho \phi_i)^{n} }{\Delta t} = - \nabla \cdot{F}^{EB}_i \]
}
%Typically, the redistribution neighborhood for each cell is one that can be
%reached via a monotonic path in each coordinate direction of unit length (see,

\commentout{
We note that the redistribution operations discussed above can be written in terms of operations with sparse matrices.  We first define an averaging operator, $\mathbf{A}^1$,
by
\[
\mathbf{A}^1_{R(i,j,k),C(ii,jj,kk)} = 
\begin{cases}
   \frac{ \vfrac_{ii,jj,kk} }{\sum\limits_{ii,jj,kk\in N(i,j,k) } \vfrac_{ii,jj,kk} } & \text{if } 0< \vfrac_{i,j,k} <1 \text{ and } \\
   & ii,jj,kk\in N(i,j,k) \\
   & \\
   0 & \text{otherwise}
\end{cases}
\]
where $R$ and $C$ map grid indices into rows and columns of a matrix.
Eq. (\ref{eq:init_update}) can then be written as
\begin{equation}
\delta U^{EB} = \left[ \mathbf{\vfrac} + (\mathbf{I}- \mathbf{\vfrac}) \mathbf{A}_1  \right ] \delta U^\cons
\label{eq:preup_matrix}
\end{equation}
where $\mathbf{\vfrac}$ is a diagonal matrix of volume fractions and $\delta U^\cons$ is a vector of the conservative updates.
We also note the
\[
\delM = \Delta x \Delta y \Delta z \; \mathbf{\vfrac} (\mathbf{I}- \mathbf{\vfrac})
(\mathbf{I} - \mathbf{A}_1 ) \delta U^\cons
\]

We now define a second averaging operator, $\mathbf{A}^2$ by
\[
\mathbf{A}^2_{R(i,j,k),C(ii,jj,kk)} = 
%\begin{cases}
 \sum_{ii,jj,kk \; | \; \ivec \in N(ii,jj,kk) \text{ and } 0 < \vfrac_{ii,jj,kk}<1}   w^{ii,jj,kk}_{i,j,k} 
\]
Eq. (\ref{eq:redist_M}) can be written
as
\[
\delta U^{EB} := \delta U^{EB} + \mathbf{A}^2 \delM = \delta U^{EB} + \mathbf{A}^2 \Delta x \Delta y \Delta z \; \mathbf{\vfrac} (\mathbf{I}- \mathbf{\vfrac})
(\mathbf{I} - \mathbf{A}_1 ) \delta U^\cons
\]
Using Eq. (\ref{eq:preup_matrix}) then gives
\[
\delta U^{EB} = \left [ \left( \mathbf{\vfrac} + (\mathbf{I}- \mathbf{\vfrac}) \mathbf{A}_1  \right ) + \mathbf{A}^2 \Delta x \Delta y \Delta z \; \mathbf{\vfrac} (\mathbf{I}- \mathbf{\vfrac})
(\mathbf{I} - \mathbf{A}_1 ) \right ] \delta U^\cons
\]
}

%% file: 02_algorithms/flux_re_redistribution.tex
This appendix assumes that the reader has first read the introduction to Section~\ref{sec:coupling_to_amr}.  Here we give the details for how to compute $\delta R$ for flux redistribution algorithms.  The approach described here is taken from Pember {\it et al.} \cite{pember1995adaptive} and the reader is referred to that paper for a more detailed description of the algorithm.

\commentout{
For coarse grid cells at the coarse/fine boundary away from the EB boundary, the synchronization between levels is simply the reflux correction discussed in \ref{sec:refluxing}.
However, in the vicinity of the intersection of the fluid-body interface with a coarse grid/fine grid boundary, the refluxing algorithm defined in \ref{sec:refluxing}
must be modified to account for the geometry and an additional correction is needed to account for redistribution across the coarse/fine boundary. The approach described here is taken from Pember {\it et al.} \cite{pember1995adaptive} and the reader is referred to that paper for a more detailed description of the algorithm.

There are four quantities that need to be captured in the re-redistribution step for  FRD.
Specifically, we need to capture the extensive contribution of
\begin{enumerate}[label=(\alph*)]
    \item coarse uncovered cells to coarse covered cells, denoted $\delta R^{\ell,u \rightarrow c}$
    \item coarse covered cells to coarse uncovered cells, denoted $\delta R^{\ell,c \rightarrow u}$ 
    \item fine valid cells to fine ghost cells, denoted $\delta R^{\ell+1,v \rightarrow g}$
    \item fine ghost cells to fine valid cells, denoted $\delta R^{\ell+1,g \rightarrow v}$
\end{enumerate}
Here, fine {\it valid} cells are cells inside level $(\ell+1)$ grids, and fine {\it ghost} cells are fine cells that are not inside level $(\ell+1)$ grids but are adjacent to fine valid cells.   
Fine valid ghost cells overlay coarse uncovered cells.}

%Correcting for redistribution across the coarse/fine boundary, referred to as re-redistribution, affects coarse cells adjacent to the coarse fine boundary that are (a) uncovered cells that redistribute into coarse cells that are covered by the fine grid and (b) uncovered cells in the neighborhood of  %Re-redistribution requires data from the redistribution on both the coarse and fine grids.
On the coarse grid, for all coarse uncovered cells $(\IVEC)$ at the coarse/fine boundary, we define
\[
\delta R^{\ell,\nto}_\IVEC = \sum_{II,JJ,KK \,\in \, \NFRD(\IVEC)\, \cap \,\mathbf{C}^\ell} \vol_{II,JJ,KK} \,
w^{\ell,\IVEC}_{II,JJ,KK} \,
    \delM^\ell_\IVEC
    \]
which accounts for all contributions of type (a), 
and
\[
\delta R^{\ell,\otn}_\IVEC =
\sum_{(II,JJ,KK) \, | \, \IVEC \,\in \, \NFRD(II,JJ,KK) \; \mathrm{and} \; (II,JJ,KK) \in \mathbf{C}^\ell}
\vol_\IVEC \,
w^{\ell,II,JJ,KK}_\IVEC \,
    \delM^\ell_{II,JJ,KK}
\]
which accounts for all contributions of type (b). 
We note here that any given coarse cell can be have contributions of both type (a) and type (b).
Here $\mathbf{C}^\ell$ is the set of covered coarse cells at level $\ell.$

We also need to capture the effect of redistribution on the fine grid.  
%Again, there are two types of contribution:  (c) ghost cells at the coarse / fine boundary that redistribute into valid fine cells, and (d) ghost cells that are in the neighborhood of a fine cell.
At all fine ghost cells, we define
\[
\delta R^{\ell+1,\gtv}_\ivec = \sum_{ii,jj,kk \,\in \, \NFRD(\ivec)\, \cap \,\mathbf{V}^{\ell+1}} \vol_{ii,jj,kk}
w^{\ell+1,\ivec}_{ii,jj,kk}
    \delM^{\ell+1}_\ivec
\]
which accounts for all contributions of type (c),
and
\[
\delta R^{\ell+1,\vtg}_\ivec = \sum_{ii,jj,kk \, | \, \ivec \,\in \, \NFRD(ii,jj,kk)\;
\mathrm{and} \; ii,jj,kk \,\in \,\mathbf{V}^{\ell+1}} \vol_{\ivec}
w^{\ell+1,ii,jj,kk}_{\ivec}
    \delM^{\ell+1}_{ii,jj,kk}
\]
which accounts for all contributions of type (d).
Here $\mathbf{V}^{\ell+1}$ is the set of valid cells at level $\ell+1$

We then define
\begin{equation*}
\delta \mathbf{R}_\IVEC^\ell
\ = \delta R^{\ell,\nto}_\IVEC - \delta R^{\ell,\otn}_\IVEC
+ \sum_{i,j,k \, | \, C^{\ell+1}_\ivec \, \subset \, \mathbf{C}^\ell_{II,JJ,KK}} (\delta R^{\ell+1,\vtg}_\ivec - \delta R^{\ell+1,\gtv}_\ivec)
%\label(eq:delR}
\end{equation*}
where the sum is over all fine grid ghost cells that overlay cell $(\IVEC)$.
The quantity $\delta R^\ell_\IVEC $ represents the difference between redistribution on the coarse and fine grid in extensive form.

See Section~{\ref{sec:final_update}} for how $\delta R^\ell_\IVEC $ is used to update the coarse and fine solutions.

\commentout{
If $\vfrac_\IVEC = 1$ then we can simply update
\[
U^\ell_\IVEC := U^\ell_\IVEC + \frac{1}{\vol_\IVEC} \,\delta R^\ell_\IVEC
\]
to complete synchronization so that the resulting scheme is conservative.
However, if $\vfrac_\IVEC < 1$ adding the entire increment $\delta R^\ell_\IVEC$ to $U^\ell_\IVEC$ could introduce an instability.  Instead update 
\[
U^\ell_\IVEC := U^\ell_\IVEC + \frac{\vfrac_\IVEC}{\vol_\IVEC} \,\delta R^\ell_\IVEC
\]
and redistribute the remaining ``mass", $(1-\vfrac_\IVEC)  \,\delta R^\ell_\IVEC$ to neighboring cells. In particular we define
\[
\vol_{nbh} = \sum_{II,JJ,KK \, \in \, \NFRD(\IVEC) } \vol_{II,JJ,KK}
\]
and then update cells in the neighborhood by
\[
U^\ell_{II,JJ,KK} := U^\ell_{II,JJ,KK} + \frac{(1-\vfrac_\IVEC)  \,\delta R^\ell_\IVEC}{\vol_{nbh}}
\]
This redistribution process can update covered coarse cells.  In that case the fine cells that cover the coarse grid cell must also be updated. In particular, if
$(II,JJ,KK) \, \in  \mathbf{C}^{\ell}$ then we update
\[
U^m_{i,j,k} := U^m_{ii,jj,kk} + \frac{(1-\vfrac_\ivec)  \,\delta R^\ell_\ivec}{\vol_{nbh}} \qquad \forall  \, \ivec \,| \,  C^m_\ivec \subset C^\ell_{II,JJ,KK} \; \mathrm{for} \;  m>\ell 
\]
This update corresponds to piecewise constant interpolation of $\frac{(1-\vfrac_\IVEC)  \,\delta R^\ell_\IVEC}{\vol_{nbh}}$ to the finer grid.

We note that for cut cells ($\vfrac_\ivec < 1$), adding the entire flux correction from Eq. (\ref{eq:reflux}) can also introduce an instability.  Consequently, reflux increments also need to be redistributed using the above procedure.  Operationally this can be done by adding the reflux corrections to $\delta R$ before redistributing it.}

\commentout{

Away from the EB boundary, the correction to the coarse grid solution due to coarse grid/fine grid interactions consists of adding a flux correction $\delta F$ to the coarse grid cells that border fine grids. This refluxing process is summarized above in \ref{sec:refluxing}.

In the vicinity of the intersection of the fluid-body interface with a coarse grid/fine grid boundary, the refluxing algorithm defined in \ref{sec:refluxing} must be replaced with a composite refluxing/redistribution algorithm that makes corrections to the solution in selected cells in order to maintain conservation. The new algorithm is needed for the following reasons:
\begin{enumerate}
    \item to account for the effect of geometry on the flux increment $\delta F$;
    \item to account for the additional transport of ``mass" due to redistribution across coarse-fine grid boundaries;
    \item to make corrections in a manner that is both stable and conservative.
\end{enumerate}
(Throughout this discussion, the term ``mass" refers to integrated values of the conservative quantities, i.e., mass, momentum, and energy as opposed to density, momentum density, and energy density.)

The composite refluxing/redistribution algorithm is applied only to coarse cells which would normally be affected by refluxing and which also satisfy the following condition:
\begin{enumerate}        
    \item The coarse cell is a mixed cell, or
    \item the coarse cell is a neighbor of a mixed coarse cell underlying the adjacent fine grid.
\end{enumerate}

The refluxing algorithm described in \ref{sec:refluxing} is used for all other coarse cells at coarse-fine grid boundaries. 
The composite algorithm consists of the following steps in an applicable coarse cell $i$,$j$,$k$ on level $L$:
\begin{enumerate}
    \item Compute a flux increment $\delta F^L_{i,j,k}$ to account for the transport of "mass" due to the difference between the level $L$ and the level $L+1$ fluxes. 
    \item Compute a redistribution increment $\delta R^L_{i,j,k}$ to account for the transport of "mass" due to redistribution across the level $L$/level $L+1$ grid boundaries.
    \item Define the coarse-fine grid correction, $\delta M^{L,c-f}_{i,j,k}$ by 
        \begin{equation}
            \delta M^{L,c-f}_{i,j,k} = \delta F^L_{i,j,k} + \delta R^L_{i,j,k}.
        \end{equation}
    \item Add the stable portion of these increments to the level $L$ solution in cell $i,j,k$, $U^{L}_{i,j,k}$, by the following: 
        \begin{equation}
            U^{L}_{i,j,k} = U^{L}_{i,j,k} + \frac{\delta M^{L,c-f}_{i,j,k}}{\Delta x\Delta y \Delta z}.
        \end{equation}
    \item If $\Lambda_{i,j,k}<1$, 
        \begin{enumerate}
            \item Redistribute $(1-\Lambda_{i,j,k})\delta M^{L,c-f}_{i,j,k}$ to the coarse grid neighbors of $i,j,k$ using volume weighted versions of (2.27) and (2.29).
            \item If a coarse grid neighbor $l,m,n$ of $i,j,k$ underlays a level $M$ cell, $M\geq L+1$, distribute the "mass" of the redistribution increment $(1-\Lambda_{i,j,k})\delta M^{L,c-f}_{i,j,k}$ onto the level $M$ cell in a volume-weighted manner.
        \end{enumerate}
    \end{enumerate}
Step 1 is identical to the operation defined in (3.4) and the accompanying text with the modification that the cell areas are weighted by area fractions. For example, equation (3.4), which computes the flux increment in the case that coarse grid cell $i,j,k$ at level $L$ shares its right $x$-cell face with a level $L+1$ grid boundary, is modified to be
\begin{align}
\begin{split}
    \delta F^{L}_{i,j,k} &= \Delta t a^L_{i+1/2,j,k}\Delta y_L \Delta z_L F^L_{i+1/2,j,k}\\
    &-\sum_{J=1,r_L}\sum_{m=r_Lj}^{r_L(j+1)-1}\sum_{n=r_Lk}^{r_L(k+1)-1}\bigg(\frac{\Delta t}{r_L}a^{L+1}_{l-1/2,m,n}\Delta y_{L+1} \Delta z_{L+1} F^{L+1,J}_{l-1/2,m,n}\bigg)
\end{split}
\end{align}
where $a^L$ and $a^{L+1}$ denote area fractions at level $L$ and $L+1$, respectively. Steps 3 and 4 are self-explanatory. Steps 2 and 5 are described below.

\subsubsection{Computation of $\delta R$}
Additional corrections to maintain conservation when the fluid-body boundary crosses a coarse-fine boundary are needed because redistribution provides an additional mechanism for numerical transport across a coarse-fine boundary.

To describe this mechanism, we divide the fine cells at a coarse-fine grid interface into two groups: the fine interior cells, i.e., the cells at the grid boundary in the interior of the grid, and the fine exterior cells, i.e., the cells in the boundary region of the grid at the grid boundary. We also divide the coarse cells into two groups: the coarse interior cells, i.e., the coarse cells underlying the fine interior cells, and the coarse exterior cells, i.e., the cells underlying the fine exterior cells. We label these four groups of cells (b), (c), (a), and (d), respectively, for ease of reference; see Figure 4.

As the level L grids advance, mixed cells in the coarse exterior region (region (a)) redistribute mass into the coarse interior region (region (d)), and mixed cells in the interior region redistribute into the exterior region. Similarly, as the level $L+1$ grids advance, mixed cells in in the fine exterior region (region (c)) redistribute mass into the fine interior region (region (b)), and mixed cells in the interior region redistribute into the exterior region. After the coarse (level $L$) and the fine (level $L+1$) grids have been advanced to the same time, the only effects of redistribution that should be present in the solution (in order that conservation be maintained) in the cells at the coarse-fine interface are those that resulted from redistribution from the fine interior cells into the fine exterior cells (i.e., from region (b) into region (c)) and from redistribution from the coarse exterior cells into the coarse interior cells (i.e., from region (a) into region (d)). However, these effects have been lost and must be recovered for the following reasons. The effect due to redistribution from fine interior cells into the fine exterior cells is lost because the fine exterior cells are only in the fine grid boundary region; the coarse grid cells (region (a)) contain the actual solution. The coarse interactions have been lost because the values on the coarse interior cells (region (d)) have been redefined by averaging the values of the overlying fine interior cells (region (b)). In addition, the effects of redistribution from the fine exterior cells (region (c)) into the fine interior cells (region (b)) and from the coarse interior cells (region (d)) into the coarse exterior cells (region (a)) are present in the solution and must be removed. These effects are present simply because the redistribution algorithm (2.27) and (2.29) applies to cells in the boundary region of the grid as well as in the interior of the grid.

In summary, there are four basic coarse-fine re-redistribution quantities that must be computed:

$\delta R^I_{L+1}$: This represents the sum of the values redistributed into the fine interior cells from the fine exterior cells over a single level L time step. Their effect is present but unwanted in the solution, and must be removed.

$\delta R^I_{L}$: This represents the values redistributed from the coarse exterior cells into the coarse interior cells over a single level L time step. Their effect is lost when the coarse values are redefined by averaging the overlying fine values and must be recovered.

$\delta R^E_{L+1}$: This represents the sum of the values redistributed from the fine interior cells into the fine exterior cells over a single level L time step. Their effect must be recovered because the solution in the boundary region is represented by the solution on the underlying coarse exterior cells.

$\delta R^E_{L}$: This represents the redistribution values from the coarse interior cells to the coarse exterior cells over a single level L time step. Their effect is present but unwanted in the solution, and must be removed.

We actually associate all four of the above corrections with the coarse exterior cells. Specifically, we associate $\delta R^I_{L}$ and $\delta R^I_{L+1}$ with the coarse exterior cell from which the redistribution values come and the coarse exterior cell that underlays the fine exterior cells from which the values come, respectively. Further, we associate $\delta R^E_{L}$ and $\delta R^E_{L+1}$ with the coarse exterior cell that receives values and the coarse exterior cell that underlays the fine cells that receive them, respectively. In other words, all four corrections are associated with the appropriate cells $i,j,k$ in region (a) that satisfy the condition mentioned at the beginning of this section: they are coarse cells that are either mixed cells or neighbors of a mixed coarse cell underlying the adjacent fine grid. We define $\delta R^L_{i,j,k}$ by
\begin{equation}\label{eq: del_R_L}
    \delta R^L_{i,j,k} = \big(\delta R^I_{L}\big)_{i,j,k} - \big(\delta R^I_{L+1}\big)_{i,j,k} + \big(\delta R^E_{L+1}\big)_{i,j,k} - \big(\delta R^E_{L}\big)_{i,j,k}.
\end{equation}

Equations (\ref{eq: delta_R_I_L}) - (\ref{eq:delta_R_E_LP1}) below define the four components of $\delta R^L_{i,j,k}$. The following terms and notation, along with Figure 4, facilitate these definitions:

$\delta M^{L}_{\mathcal{I},\mathcal{J},\mathcal{K};\mathcal{L},\mathcal{M},\mathcal{N}}$: This is the amount of "mass" redistributed from level $L$ mixed cell $\mathcal{I},\mathcal{J},\mathcal{K}$ to level $L$ cell $\mathcal{L},\mathcal{M},\mathcal{N}$ while advancing the level $L$ solution one level $L$ time step from $t$ to $t+\Delta t$ as computed by (2.28).

$\delta M^{L+1}_{\mathcal{I},\mathcal{J},\mathcal{K};\mathcal{L},\mathcal{M},\mathcal{N}}$: This is the amount of "mass" redistributed from level $L+1$ mixed cell $\mathcal{I},\mathcal{J},\mathcal{K}$ to level $L+1$ cell $\mathcal{L},\mathcal{M},\mathcal{N}$ while advancing the level $L+1$ solution one level $L+1$ time step from $t+(J-1)\Delta t/r_L$ to $t+J\Delta t/r_L$ as computed by (2.28).

(a-d)$_{\mathcal{I},\mathcal{J},\mathcal{K}}$: This is the set of level $L$ cells $\mathcal{L},\mathcal{M},\mathcal{N}$ in region (d) to which "mass" is redistributed from level $L$ mixed cell $\mathcal{I},\mathcal{J},\mathcal{K}$ in region (a):

\begin{align*}
\text{(a-d)}_{\mathcal{I},\mathcal{J},\mathcal{K}} = nbh(\mathcal{I},\mathcal{J},\mathcal{K}) \cap \{(\mathcal{L},\mathcal{M},\mathcal{N}):(\mathcal{L},\mathcal{M},\mathcal{N})\in \text{(d) and } \Lambda_{\mathcal{L},\mathcal{M},\mathcal{N}}>0\}
\end{align*}

(d-a)$_{\mathcal{I},\mathcal{J},\mathcal{K}}$: This is the set of level $L$ mixed cells $\mathcal{L},\mathcal{M},\mathcal{N}$ in region (d) which redistribute "mass" to level $L$ cell $\mathcal{I},\mathcal{J},\mathcal{K}$ in region (a):

\begin{align*}
\text{(d-a)}_{\mathcal{I},\mathcal{J},\mathcal{K}} = nbh(\mathcal{I},\mathcal{J},\mathcal{K}) \cap \{(\mathcal{L},\mathcal{M},\mathcal{N}):(\mathcal{L},\mathcal{M},\mathcal{N})\in \text{(d) and } 0<\Lambda_{\mathcal{L},\mathcal{M},\mathcal{N}}<1\}
\end{align*}

(b-c)$^{src}_{\mathcal{I},\mathcal{J},\mathcal{K}}$: This is the set of level $L+1$ mixed cells $\mathcal{L},\mathcal{M},\mathcal{N}$ in region (b) which redistribute "mass" to level $L+1$ cells in region (c) overlying level $L$ cell $\mathcal{I},\mathcal{J},\mathcal{K}$ in region (a):

\begin{align*}
\text{(b-c)}^{src}_{\mathcal{I},\mathcal{J},\mathcal{K}} = &\{(\mathcal{L},\mathcal{M},\mathcal{N}):((\mathcal{L},\mathcal{M},\mathcal{N})\in \text{(b) and } 0<\Lambda_{\mathcal{L},\mathcal{M},\mathcal{N}}<1 \text{ and }\\
&(\mathcal{L},\mathcal{M},\mathcal{N})\in nbh(\mathcal{P},\mathcal{Q},\mathcal{R}):((\mathcal{P},\mathcal{Q},\mathcal{R}) \text{ overlays } (\mathcal{I},\mathcal{J},\mathcal{K}) \\
&\text{ and } \Lambda_{\mathcal{P},\mathcal{Q},\mathcal{R}}>0))\}
\end{align*}

(b-c)$^{dst}_{\mathcal{I},\mathcal{J},\mathcal{K}}$: This is the set of level $L+1$ cells $\mathcal{L},\mathcal{M},\mathcal{N}$ in region (c) overlying level $L$ cell $\mathcal{I},\mathcal{J},\mathcal{K}$ in region (a) to which "mass" is redistributed from level $L+1$ mixed cells in region (b):

\begin{align*}
\text{(b-c)}^{dst}_{\mathcal{I},\mathcal{J},\mathcal{K}} = &\{(\mathcal{L},\mathcal{M},\mathcal{N}):((\mathcal{L},\mathcal{M},\mathcal{N}) \text{ overlays } (\mathcal{I},\mathcal{J},\mathcal{K}) \text{ and }\\
&\Lambda_{\mathcal{L},\mathcal{M},\mathcal{N}}>0 \text{ and }\\
&(\mathcal{L},\mathcal{M},\mathcal{N})\in nbh(\mathcal{P},\mathcal{Q},\mathcal{R}):((\mathcal{P},\mathcal{Q},\mathcal{R})\in \text{(b) and } (\mathcal{I},\mathcal{J},\mathcal{K}) \\
&\text{ and } 0<\Lambda_{\mathcal{P},\mathcal{Q},\mathcal{R}}<1))\}
\end{align*}

(c-b)$^{src}_{\mathcal{I},\mathcal{J},\mathcal{K}}$: This is the set of level $L+1$ mixed cells $\mathcal{L},\mathcal{M},\mathcal{N}$ in region (c) which overlay $\mathcal{I},\mathcal{J},\mathcal{K}$ in region (a) and which redistribute "mass" to level $L+1$ cells in region (b):

\begin{align*}
\text{(c-b)}^{src}_{\mathcal{I},\mathcal{J},\mathcal{K}} = &\{(\mathcal{L},\mathcal{M},\mathcal{N}):((\mathcal{L},\mathcal{M},\mathcal{N}) \text{ overlays } (\mathcal{I},\mathcal{J},\mathcal{K}) \\
&\text{ and } 0<\Lambda_{\mathcal{L},\mathcal{M},\mathcal{N}}<1))\}
\end{align*}

(c-b)$^{dst}_{\mathcal{I},\mathcal{J},\mathcal{K}}$: This is the set of level $L+1$ cells $\mathcal{L},\mathcal{M},\mathcal{N}$ in region (b) to which mass is redistributed by mixed cells in region (c) which overlay $\mathcal{I},\mathcal{J},\mathcal{K}$ in region (a):

\begin{align*}
\text{(c-b)}^{dst}_{\mathcal{I},\mathcal{J},\mathcal{K}} = &\{(\mathcal{L},\mathcal{M},\mathcal{N}):((\mathcal{L},\mathcal{M},\mathcal{N})\in \text{(b) and } \Lambda_{\mathcal{L},\mathcal{M},\mathcal{N}}>0\\
&(\mathcal{L},\mathcal{M},\mathcal{N})\in nbh(\mathcal{P},\mathcal{Q},\mathcal{R}):((\mathcal{P},\mathcal{Q},\mathcal{R}) \text{ overlays } (\mathcal{I},\mathcal{J},\mathcal{K}) \\
&\text{ and } 0<\Lambda_{\mathcal{P},\mathcal{Q},\mathcal{R}}<1))\}
\end{align*}

Using the forgoing terminology, we can now define the four components of $\delta R^L_{i,j,k}$ in (\ref{eq: del_R_L}) as follows:

\begin{equation}\label{eq: delta_R_I_L}
    \big(\delta R^I_{L}\big)_{i,j,k} = \sum_{(l,m,n)\in\text{(a-d)}_{i,j,k}} \delta M^L_{i,j,k;l,m,n}
\end{equation}

\begin{equation}
    \big(\delta R^E_{L}\big)_{i,j,k} = \sum_{(l,m,n)\in\text{(d-a)}_{i,j,k}} \delta M^L_{l,m,n;i,j,k}
\end{equation}

\begin{equation}
    \big(\delta R^I_{L+1}\big)_{i,j,k} = \sum_{J=1}^{r_L}\sum_{(l,m,n)\in \text{(c-b)}^{src}_{i,j,k}}\Bigg(\sum_{(I,J,K)\in nbh(l,m,n)\cap\text{(c-b)}^{dst}_{i,j,k}}\delta M^{L+1,J}_{l,m,n:I,J,K}\Bigg)
\end{equation}

\begin{equation}\label{eq:delta_R_E_LP1}
    \big(\delta R^E_{L+1}\big)_{i,j,k} = \sum_{J=1}^{r_L}\sum_{(l,m,n)\in \text{(b-c)}^{dst}_{i,j,k}}\Bigg(\sum_{(I,J,K)\in nbh(l,m,n)\cap\text{(b-c)}^{src}_{i,j,k}}\delta M^{L+1,J}_{I,J,K:l,m,n}\Bigg)
\end{equation}

\subsubsection{Redistribution of $\delta M^{L,c-f}_{i,j,k}$}

If $\Lambda_{i,j,k}<1$, we must redistribute $(1-\Lambda_{i,j,k})\delta M^{L,c-f}_{i,j,k}$ to the level $L$ neighbors of $i,j,k$ using volume weighted versions of (2.27) and (2.29) as follows. We define 

\begin{equation}
    m^{red}_{i,j,k} = \sum_{level\; L \;cells\;(l,m,n)\in nbh(i,j,k)} \Lambda_{l,m,n}\Delta x\Delta y\Delta z.
\end{equation}

Then, for all coarse cells $(l,m,n)\in nbh(i,j,k)$, the value of $U^{n+1}_{l,m,n}$ is modified by the following:

\begin{equation}\label{eq: Unp1}
    U^{n+1}_{l,m,n} = U^{n+1}_{l,m,n} + \frac{(1-\Lambda_{i,j,k})\delta M^{L,c-f}_{i,j,k}}{m^{red}_{i,j,k}}.
\end{equation}

Further, if any of the level $L$ cells $l,m,n$ underlays a level $M$ cell $I,J,K$, where $M>L$, we modify the level M solution in $I,J,K$ by

\begin{equation}
    U^{n+1}_{I,J,K} = U^{n+1}_{I,J,K} + \frac{(1-\Lambda_{i,j,k})\delta M^{L,c-f}_{i,j,k}}{m^{red}_{i,j,k}},
\end{equation}

in other words, using the same equation as (\ref{eq: Unp1}). The same expressions are used  because the total amount of mass redistributed to a cell $\mathcal{I},\mathcal{J},\mathcal{K}$ on level $\mathcal{L}$, $\mathcal{L}=L$ or $M$, with volume fraction $\Lambda_{\mathcal{I},\mathcal{J},\mathcal{K}}$ is 

\begin{equation}
    \Lambda_{\mathcal{I},\mathcal{J},\mathcal{K}}\Delta x_{\mathcal{L}}\Delta y_{\mathcal{L}}\Delta z_{\mathcal{L}}\frac{(1-\Lambda_{i,j,k})\delta M^{L,c-f}_{i,j,k}}{m^{red}_{i,j,k}};
\end{equation}

the update to the solution $U$ is found by dividing this expression by $\Lambda_{\mathcal{I},\mathcal{J},\mathcal{K}}\Delta x_{\mathcal{L}}\Delta y_{\mathcal{L}}\Delta z_{\mathcal{L}}$.

}

%% file: 02_algorithms/frd_matrix.tex
We note that the redistribution operations discussed above can be written in terms of operations with sparse matrices.  We first define an averaging operator, $\mathbf{A}_1$,
by
\[
\mathbf{A}_{1_{R_{i,j,k},C_{ii,jj,kk}} }= 
\begin{cases}
   \frac{ \vfrac_{ii,jj,kk} }{\sum\limits_{ii,jj,kk\in N(i,j,k) } \vfrac_{ii,jj,kk} } & \text{if } 0< \vfrac_{i,j,k} <1 \text{ and } \\
   & ii,jj,kk\in N(i,j,k) \\
   & \\
   0 & \text{otherwise}
\end{cases}
\]
where $R_{i,j,k}$ and $C_{i,j,k}$ map grid indices into rows and columns of a matrix.
Eq. (\ref{eq:init_update}) can then be written as
\begin{equation}
\delta U^{EB,*} = \left[ \mathbf{\vfrac} + (\mathbf{I}- \mathbf{\vfrac}) \mathbf{A}_1  \right ] \delta U^\cons
\label{eq:preup_matrix}
\end{equation}
where $\mathbf{\vfrac}$ is a diagonal matrix of volume fractions and $\delta U^\cons$ is a vector of the conservative updates.
We also note that
\[
\delM = \Delta x \Delta y \Delta z \; \mathbf{\vfrac} (\mathbf{I}- \mathbf{\vfrac})
(\mathbf{I} - \mathbf{A}_1 ) \delta U^\cons
\]

We now define a second averaging operator, $\mathbf{A}_2$ by
\[
\mathbf{A}_{2_{R_{i,j,k},C_{ii,jj,kk}} }= 
\begin{cases}
w^{ii,jj,kk}_{i,j,k}  & \text{if } 0< \vfrac_{ii,jj,kk} <1 \text{ and } \\
   & ii,jj,kk\in N(i,j,k) \\
   & \\
   0 & \text{otherwise}
\end{cases}
\]
\MarginPar{Early version had a sum that wasn't correct}
\commentout{
\[
\mathbf{A}^2_{R(i,j,k),C(ii,jj,kk)} = 
%\begin{cases}
 \sum_{ii,jj,kk \; | \; \ivec \in N(ii,jj,kk) \text{ and } 0 < \vfrac_{ii,jj,kk}<1}   w^{ii,jj,kk}_{i,j,k} 
\]
}
Eq. (\ref{eq:redist_M}) can be written
as
\[
\delta U^{EB} = \delta U^{EB,*} + \mathbf{A}_2 \delM = \delta U^{EB,*} + \mathbf{A}_2 \Delta x \Delta y \Delta z \; \mathbf{\vfrac} (\mathbf{I}- \mathbf{\vfrac})
(\mathbf{I} - \mathbf{A}_1 ) \delta U^\cons
\]
Using Eq. (\ref{eq:preup_matrix}) then gives
\[
\delta U^{EB} = \left [ \left( \mathbf{\vfrac} + (\mathbf{I}- \mathbf{\vfrac}) \mathbf{A}_1  \right ) + \mathbf{A}_2 \Delta x \Delta y \Delta z \; \mathbf{\vfrac} (\mathbf{I}- \mathbf{\vfrac})
(\mathbf{I} - \mathbf{A}_1 ) \right ] \delta U^\cons
\]

%% file: 02_algorithms/frd_matrix_rere.tex
%Let $n$ be the dimension of $\mathbf{A}_2$. 
Here, we show here how to construct the $\delta R$'s from the matrix $\Diag(V)\mathbf{A}_2\Diag(\delta M)$, where $\Diag(V)$ is the diagonal matrix of cell volumes. $V$ and $\Diag(\delta M)$ is the diagonal matrix whose diagonal entries are the $\delta M$s.%\\\MarginPar{I will fix all at a later time}

First, we consider the coarse grid. Let $\mathbf{C}^\ell$ be the set of covered coarse cells at level $\ell$. 
We denote the row corresponding to cell $(\IVEC)$ as $R_\IVEC$ and the column corresponding to cell $(\IVEC)$ as $C_\IVEC$. Then we can observe that the entry $(\Diag(V)\mathbf{A}_2\Diag(\delta M))_{R_\IVEC,C_{R,S,T}}$ contains the contribution of cell $(R,S,T)$ to cell $(\IVEC)$. 
Referring to the definitions in \ref{sec:FRRD}, we observe that
\[
\delta R_{\IVEC}^{l,u\rightarrow c}=\sum_{(\IVEC)|(\IVEC)\in N(R,S,T) \text{ and } (R,S,T)\in \mathbf{C}^\ell} (\Diag(V)\mathbf{A}_2\Diag(\delta M))_{R_{R,S,T},C_\IVEC}
\]
and
\[
\delta R_{\IVEC}^{l,c\rightarrow u}=\sum_{(R,S,T)\in N(\IVEC) \cap \mathbf{C}^\ell} (\Diag(V)\mathbf{A}_2\Diag(\delta M))_{R_\IVEC,C_{R,S,T}}
\]

Regarding the $\delta R$'s corresponding to the fine grid, and denoting $\mathbf{V}^{\ell+1}$ to be the set of valid cells at level $\ell+1$, we have that 
\[
\delta R_{\ivec}^{l,g\rightarrow v}=\sum_{(\ivec)|(\ivec)\in N(r,s,t) \text{ and } (r,s,t)\in \mathbf{V}^{\ell+1}} (\Diag(V)\mathbf{A}_2\Diag(\delta M))_{R_{r,s,t},C_\ivec}
\]
and
\[
\delta R_{\ivec}^{l,v\rightarrow g}=\sum_{(r,s,t)\in N(\ivec) \cap \mathbf{V}^{\ell+1}} (\Diag(V)\mathbf{A}_2\Diag(\delta M))_{R_\ivec,C_{r,s,t}}
\]

%% file: 02_algorithms/ebgdnv.tex
We use a piecewise linear Godunov method based on the unsplit discretization of Colella \cite{colella1990multidimensional} and Saltzman \cite{SALTZMAN1994}.  The method uses a fourth-order reconstruction of slopes within a cell.  The solution is predicted to cell faces and half-time levels using a characteristic tracing procedure.  Derivatives transverse to the face are treated in conservation form and act as sources in the characteristic tracing.  

In the version we use with cut cells, we view $U$ as being defined at cell centers (as opposed to centroids) even for cut cells.  When we use these values to reconstruct slopes within each (regular or cut) cell, we use the original fourth-order slopes if the fourth-order stencil does not cross any faces with zero area.   If the fourth-order computation is disallowed, we reduce to second-order slopes in that direction if the second-order stencil does not cross any faces with zero area.  
If the second-order computation is disallowed the slopes in that direction are zeroed.

We then compute fluxes only at faces with nonzero area.  If, for a given cell, in trying to compute the flux on a face, the area fractions of all of the faces on which we need values to compute the transverse derivatives that contribute to that flux are not nonzero, then we do not include any transverse terms.  For example, if we are predicting from a cell center to an $x$-face, unless both of the $y$ faces and both of $z$ faces have nonzero area fractions, we do not include transverse derivatives in the computation of the flux on the $x$-face.

Once we have predicted values on either side of faces with nonzero area fractions, we then compute the flux at each face using a two-shock approximate Riemann solver.  We view the flux as being defined at the center of the face, but when we compute the divergence of fluxes needed in Eq. (\ref{eq:cons_up}) to update the solution, we use bilinear interpolation to interpolate the fluxes from face centers to face centroids.

%% file: main.bbl
\begin{thebibliography}{10}
\expandafter\ifx\csname url\endcsname\relax
  \def\url#1{\texttt{#1}}\fi
\expandafter\ifx\csname urlprefix\endcsname\relax\def\urlprefix{URL }\fi
\expandafter\ifx\csname href\endcsname\relax
  \def\href#1#2{#2} \def\path#1{#1}\fi

\bibitem{gokhaleNikosKlein:2018}
N.~Gokhale, N.~Nikiforakis, R.~Klein, A dimensionally split {C}artesian cut
  cell method for hyperbolic conservation laws, J. Comput. Phys. 364 (2018)
  186--208.

\bibitem{mjb-hel-rjl:hbox}
M.~J. Berger, C.~Helzel, R.~J. LeVeque, {H}-box methods for the approximation
  of one-dimensional conservation laws on irregular grids, SIAM J. Numer. Anal.
  41 (2003) 893--918.

\bibitem{MURALIDHARAN2016}
B.~Muralidharan, S.~Menon,
  \href{https://www.sciencedirect.com/science/article/pii/S0021999116301954}{A
  high-order adaptive cartesian cut-cell method for simulation of compressible
  viscous flow over immersed bodies}, Journal of Computational Physics 321
  (2016) 342--368.
\newblock \href {https://doi.org/https://doi.org/10.1016/j.jcp.2016.05.050}
  {\path{doi:https://doi.org/10.1016/j.jcp.2016.05.050}}.
\newline\urlprefix\url{https://www.sciencedirect.com/science/article/pii/S0021999116301954}

\bibitem{saye2017implicit}
R.~Saye, Implicit mesh discontinuous galerkin methods and interfacial gauge
  methods for high-order accurate interface dynamics, with applications to
  surface tension dynamics, rigid body fluid--structure interaction, and free
  surface flow: Part i, Journal of Computational Physics 344 (2017) 647--682.

\bibitem{saye2017part2}
R.~Saye,
  \href{https://www.sciencedirect.com/science/article/pii/S0021999117303728}{Implicit
  mesh discontinuous galerkin methods and interfacial gauge methods for
  high-order accurate interface dynamics, with applications to surface tension
  dynamics, rigid body fluid–structure interaction, and free surface flow:
  Part ii}, Journal of Computational Physics 344 (2017) 683--723.
\newblock \href {https://doi.org/https://doi.org/10.1016/j.jcp.2017.05.003}
  {\path{doi:https://doi.org/10.1016/j.jcp.2017.05.003}}.
\newline\urlprefix\url{https://www.sciencedirect.com/science/article/pii/S0021999117303728}

\bibitem{GULIZZI2022}
V.~Gulizzi, A.~S. Almgren, J.~B. Bell, A coupled discontinuous galerkin-finite
  volume framework for solving gas dynamics over embedded geometries, Journal
  of Computational Physics 450 (2022) 110861.

\bibitem{pember1995adaptive}
R.~B. Pember, J.~B. Bell, P.~Colella, W.~Y. Curtchfield, M.~L. Welcome, An
  adaptive cartesian grid method for unsteady compressible flow in irregular
  regions, Journal of computational Physics 120~(2) (1995) 278--304.

\bibitem{colella2006cartesian}
P.~Colella, D.~T. Graves, B.~J. Keen, D.~Modiano, A cartesian grid embedded
  boundary method for hyperbolic conservation laws, Journal of Computational
  Physics 211~(1) (2006) 347--366.

\bibitem{hu2006conservative}
X.~Y. Hu, B.~Khoo, N.~A. Adams, F.~Huang, A conservative interface method for
  compressible flows, Journal of Computational Physics 219~(2) (2006) 553--578.

\bibitem{klein2009well}
R.~Klein, K.~Bates, N.~Nikiforakis, Well-balanced compressible cut-cell
  simulation of atmospheric flow, Philosophical Transactions of the Royal
  Society A: Mathematical, Physical and Engineering Sciences 367~(1907) (2009)
  4559--4575.

\bibitem{graves2013cartesian}
D.~Graves, P.~Colella, D.~Modiano, J.~Johnson, B.~Sjogreen, X.~Gao, A cartesian
  grid embedded boundary method for the compressible navier--stokes equations,
  Communications in Applied Mathematics and Computational Science 8~(1) (2013)
  99--122.

\bibitem{SCHNEIDERS2013786}
L.~Schneiders, D.~Hartmann, M.~Meinke, W.~Schröder,
  \href{https://www.sciencedirect.com/science/article/pii/S0021999112005839}{An
  accurate moving boundary formulation in cut-cell methods}, Journal of
  Computational Physics 235 (2013) 786--809.
\newblock \href {https://doi.org/https://doi.org/10.1016/j.jcp.2012.09.038}
  {\path{doi:https://doi.org/10.1016/j.jcp.2012.09.038}}.
\newline\urlprefix\url{https://www.sciencedirect.com/science/article/pii/S0021999112005839}

\bibitem{almgren1997cartesian}
A.~S. Almgren, J.~B. Bell, P.~Colella, T.~Marthaler, A cartesian grid
  projection method for the incompressible euler equations in complex
  geometries, SIAM Journal on Scientific Computing 18~(5) (1997) 1289--1309.

\bibitem{trebotich:2015}
D.~Trebotich, D.~Graves, An adpative finite volume method for the
  incompressible {N}avier-{S}tokes equations in complex geometries,
  Communications in Applied Mathematics and Computational Science (2015)
  43--82.

\bibitem{berger2021state}
M.~Berger, A.~Giuliani, A state redistribution algorithm for finite volume
  schemes on cut cell meshes, Journal of Computational Physics 428 (2021)
  109820.

\bibitem{giuliani2022dgsrd}
A.~Giuliani, A two-dimensional stabilized discontinuous galerkin method on
  curvilinear embedded boundary grids, SIAM Journal on Scientific Computing
  44~(1) (2022) A389--A415.

\bibitem{giuliani2022weighted}
A.~Giuliani, A.~S. Almgren, J.~B. Bell, M.~J. Berger, M.~H. de~Frahan,
  D.~Rangarajan, A weighted state redistribution algorithm for embedded
  boundary grids, Journal of Computational Physics (2022) 111305.

\bibitem{berger2023new}
M.~Berger, A.~Giuliani, A new provably stable weighted state redistribution
  algorithm (2023).
\newblock \href {http://arxiv.org/abs/2308.16332} {\path{arXiv:2308.16332}}.

\bibitem{berger1989local}
M.~J. Berger, P.~Colella, Local adaptive mesh refinement for shock
  hydrodynamics, Journal of computational Physics 82~(1) (1989) 64--84.

\bibitem{bell1994three}
J.~Bell, M.~Berger, J.~Saltzman, M.~Welcome, Three-dimensional adaptive mesh
  refinement for hyperbolic conservation laws, SIAM Journal on Scientific
  Computing 15~(1) (1994) 127--138.

\bibitem{AMREX:IJHPCA}
W.~Zhang, A.~Myers, K.~Gott, A.~Almgren, J.~Bell, {AMReX}: Block-structured
  adaptive mesh refinement for multiphysics applications, The International
  Journal of High Performance Computing Applications 35~(6) (2021) 508--526.
\newblock \href {https://doi.org/10.1177/10943420211022811}
  {\path{doi:10.1177/10943420211022811}}.

\bibitem{bell2011adaptive}
J.~Bell, M.~Day, Adaptive methods for simulation of turbulent combustion, in:
  Turbulent Combustion Modeling: Advances, New Trends and perspectives,
  Springer, 2011, pp. 301--329.

\bibitem{colella1990multidimensional}
P.~Colella, Multidimensional upwind methods for hyperbolic conservation laws,
  Journal of Computational Physics 87~(1) (1990) 171--200.

\bibitem{SALTZMAN1994}
J.~Saltzman,
  \href{https://www.sciencedirect.com/science/article/pii/S0021999184711843}{An
  unsplit 3d upwind method for hyperbolic conservation laws}, Journal of
  Computational Physics 115~(1) (1994) 153--168.
\newblock \href {https://doi.org/https://doi.org/10.1006/jcph.1994.1184}
  {\path{doi:https://doi.org/10.1006/jcph.1994.1184}}.
\newline\urlprefix\url{https://www.sciencedirect.com/science/article/pii/S0021999184711843}

\bibitem{EmmettETAL2019}
M.~Emmett, E.~Motheau, W.~Zhang, M.~Minion, J.~B. Bell, A fourth-order adaptive
  mesh refinement algorithm for the multicomponent, reacting compressible
  navier–stokes equations, Combustion Theory and Modelling 23~(4) (2019)
  592--625.

\bibitem{url:CAMR}
\href{https://github.com/AMReX-Codes/CAMR}{{CAMR} code}.
\newline\urlprefix\url{https://github.com/AMReX-Codes/CAMR}

\bibitem{PeleC_IJHPCA}
M.~T. {Henry de Frahan}, J.~S. Rood, M.~S. Day, H.~Sitaraman, S.~Yellapantula,
  B.~A. Perry, R.~W. Grout, A.~Almgren, W.~Zhang, J.~B. Bell, J.~H. Chen,
  \href{https://doi.org/10.1177/10943420221121151}{{PeleC: An adaptive mesh
  refinement solver for compressible reacting flows}}, The International
  Journal of High Performance Computing Applications 37~(2) (2022) 115--131.
\newblock \href {https://doi.org/10.1177/10943420221121151}
  {\path{doi:10.1177/10943420221121151}}.
\newline\urlprefix\url{https://doi.org/10.1177/10943420221121151}

\bibitem{PeleViz}
H.~Sitaraman, N.~Brunhart-Lupo, M.~H. de~Frahan, S.~Yellapantula, B.~Perry,
  J.~Rood, R.~Grout, M.~Day, R.~Binyahib, K.~Gruchalla,
  \href{https://link.aps.org/doi/10.1103/PhysRevFluids.6.110504}{Visualizations
  of direct fuel injection effects in a supersonic cavity flameholder}, Phys.
  Rev. Fluids 6 (2021) 110504.
\newblock \href {https://doi.org/10.1103/PhysRevFluids.6.110504}
  {\path{doi:10.1103/PhysRevFluids.6.110504}}.
\newline\urlprefix\url{https://link.aps.org/doi/10.1103/PhysRevFluids.6.110504}

\bibitem{Sitaraman2021}
H.~Sitaraman, S.~Yellapantula, M.~T. {Henry de Frahan}, B.~Perry, J.~Rood,
  R.~Grout, M.~Day,
  \href{https://www.sciencedirect.com/science/article/pii/S0010218021002741}{Adaptive
  mesh based combustion simulations of direct fuel injection effects in a
  supersonic cavity flame-holder}, Combustion and Flame 232 (2021) 111531.
\newblock \href
  {https://doi.org/https://doi.org/10.1016/j.combustflame.2021.111531}
  {\path{doi:https://doi.org/10.1016/j.combustflame.2021.111531}}.
\newline\urlprefix\url{https://www.sciencedirect.com/science/article/pii/S0010218021002741}

\bibitem{ni2022immersed}
R.~Ni, J.~Li, X.~Zhang, X.~Zhou, X.~Cui, An immersed boundary-material point
  method for shock-structure interaction and dynamic fracture, Journal of
  Computational Physics 470 (2022) 111558.

\bibitem{jiang2022development}
C.~Jiang, J.~Pan, Y.~Zhu, J.~Li, E.~K. Quaye, Development and verification of a
  high-speed compressible reactive flow solver in openfoam, Journal of
  Computational Science 63 (2022) 101780.

\bibitem{yang1987computation}
J.~Yang, Y.~Liu, H.~Lomax, Computation of shock wave reflection by circular
  cylinders, AIAA journal 25~(5) (1987) 683--689.

\bibitem{chaudhuri2011use}
A.~Chaudhuri, A.~Hadjadj, A.~Chinnayya, On the use of immersed boundary methods
  for shock/obstacle interactions, Journal of Computational Physics 230~(5)
  (2011) 1731--1748.

\bibitem{zoltak1998hybrid}
J.~Z{\'o}{\l}tak, D.~Drikakis, Hybrid upwind methods for the simulation of
  unsteady shock-wave diffraction over a cylinder, Computer Methods in Applied
  Mechanics and Engineering 162~(1-4) (1998) 165--185.

\bibitem{whitham1957new}
G.~Whitham, A new approach to problems of shock dynamics part i two-dimensional
  problems, Journal of Fluid Mechanics 2~(2) (1957) 145--171.

\bibitem{bryson1961diffraction}
A.~Bryson, R.~Gross, Diffraction of strong shocks by cones, cylinders, and
  spheres, Journal of Fluid Mechanics 10~(1) (1961) 1--16.

\bibitem{toro2013riemann}
E.~F. Toro, Riemann solvers and numerical methods for fluid dynamics: a
  practical introduction, Springer Science \& Business Media, 2013.

\bibitem{eagle-ref}
{National Renewable Energy Laboratory},
  \href{https://www.nrel.gov/hpc/eagle-system.html}{Eagle, {National Renewable
  Energy Laboratory}} (2022).
\newline\urlprefix\url{https://www.nrel.gov/hpc/eagle-system.html}

\end{thebibliography}
